\documentclass[a4paper,11pt]{article}
\usepackage{amsmath,amssymb,amsthm,graphicx}
\usepackage{hyperref,xurl,enumitem,bm,esint}
\usepackage{color}
\usepackage{comment}
\usepackage{algorithm}
\usepackage{algpseudocode}
\usepackage[T1]{fontenc}
\usepackage[normalem]{ulem}
\usepackage{pifont}
\usepackage{xcolor,enumitem}
\usepackage{graphicx}
\usepackage{multirow}
\usepackage{makecell} 
\usepackage{subcaption}
\usepackage{mwe}
\usepackage{overpic}
\usepackage[margin=30mm]{geometry}
\usepackage{booktabs}

\newcommand{\be}{\mathbf e}
\newcommand{\bd}{\mathbf d}
\newcommand{\bc}{\mathbf c}
\newcommand{\bb}{\mathbf b}

\newcommand{\bx}{\bm{x}}
\newcommand{\bn}{\mathbf n}
\newcommand{\br}{\mathbf r}
\newcommand{\bA}{\mathbf A}
\newcommand{\bp}{\mathbf p}
\newcommand{\bH}{\mathbf H}

\newcommand{\anna}[2][cyan]{\textcolor{#1}{#2}}
\newcommand{\ming}[2][magenta]{\textcolor{#1}{#2}}

\newcommand{\baoli}[2][red]{\textcolor{#1}{#2}}
\theoremstyle{plain}

\newtheorem{remark}{Remark}[section]

\numberwithin{equation}{section}

\DeclareMathOperator*{\argmin}{arg\,min}

\def\cl {\nonumber \\}
\def\el {\nonumber }

\title{Crime hotspot dynamics in residential burglary models with police response}

\author{Baoli Hao$^1$, Kamrun Mily$^2$, Annalisa Quaini$^2$, and Ming Zhong$^2$}

\date{}
\begin{document}
\maketitle

\begin{center}
$^1$ Department of Applied Mathematics, Illinois Institute of Technology, 10 West 35th Street
Chicago, IL 60616 \\ {\tt bhao2@hawk.illinoistech.edu}

\noindent $^2$ Department of Mathematics, University of Houston, 3551 Cullen Blvd, Houston TX 77204, USA \\
{\tt kmily@cougarnet.uh.edu, \{aquaini, mzhong3\}@central.uh.edu}
\end{center}

\begin{abstract}

We develop and analyze mathematical models for residential burglary that incorporates police deployment through a delayed feedback mechanism. Motivated by empirical observations from publicly available crime and policing data, we extend a well-known agent-based model by introducing a dynamic police response driven by crime information that becomes available only after a finite delay. 
Taking the mean-field limit, we derive a coupled continuum system consisting of three partial differential equations and one ordinary differential equation describing the interactions among criminal density, environmental attractiveness, delayed crime signal, and police deployment. Linear stability analysis of homogeneous steady states reveals that response delays can destabilize otherwise stable equilibria through Hopf bifurcations. As a result, the model predicts sustained temporal oscillations and dynamically evolving crime hotspots.
Numerical simulations of both the agent-based and continuum models confirm the theoretical analysis and uncover rich spatio-temporal behaviors, including moving, splitting, and merging hotspots.
Through a parametric study, we investigate the roles of police density, crime information delay, and neighborhood effects in controlling stability, hotspot size, and oscillatory behavior. 
Our results indicate that timely access to crime data plays a more important role than police density in stabilizing crime levels. 
\end{abstract}


\vskip .2cm
\noindent \emph{AMSC}: 35Q91, 65M22, 91D10.

\section{Introduction}\label{sec:intro}
%
Predictability and control of urban crime remain fundamental challenges for public safety. 
Crime patterns are known to exhibit complex spatial and temporal variability, often characterized by the emergence, evolution, and dissipation of hotspots. 
See, e.g., ~\cite{Brantinghambook,Eck1995,Chaineybook,Weisburd2015,mohler2019reducing}.
Continuum models derived from agent-based descriptions (see, e.g.,~\cite{short2008statistical,short2010dissipation,pan2018,Hao2026}) have provided a mathematically tractable framework to study such phenomena in the case of residential burglary, which is simpler than other crimes since mobile offenders target stationary sites. These models, which describe the coupled evolution of criminal density and environmental attractiveness, exhibit mechanisms analogous to chemotactic aggregation in biological systems \cite{KELLER1970399,Byrne2004,doi:10.1137/050637923,BellomoWinkler2022}, leading to the formation of localized crime hotspots. Although mathematically interesting and computationally challenging, several of these models do not incorporate law enforcement.

The integration of police intervention in continuum models has been explored through deterrence terms with
fixed police deployment ~\cite{short2010nonlinear} (i.e., the criminal activity is allowed to evolve with no policing till
a given time at which time 
police density is computed and held fixed thereafter) or optimal control strategies \cite{PITCHER_2010,zipkin2014}. 
{For an explicit treatment of enforcement dynamics in agent-based models with instantaneously responding police, we refer to \cite{jones2010}.}
These studies demonstrate that police presence can significantly alter the spatial structure of crime, suppressing or redistributing hotspots depending on deployment strategies. However, most existing agent-based or continuum models assume that police response is either static or instantaneous 
and globally optimal, failing to fully capture the operational constraints of real-world law enforcement. 
Other kinds of models that include police response are ODE models \cite{10.3389/fams.2022.1086745,https://doi.org/10.1155/jom/1372780}
and, more recently, hybrid agent-based/machine learning models~\cite{Nurhan_Short_2026}.

\begin{figure}[htb!]
    \centering
    \includegraphics[width=0.9\textwidth]{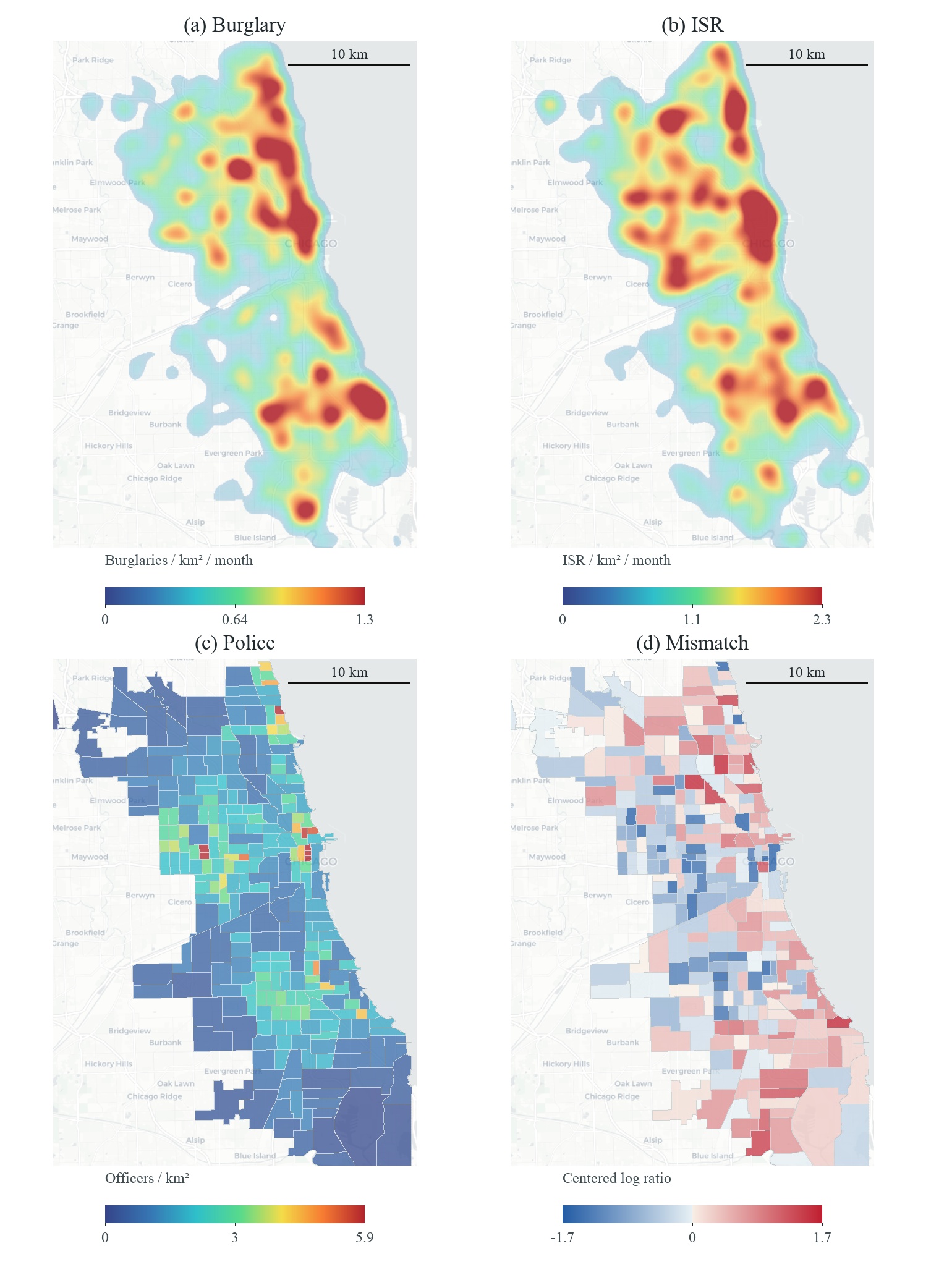}
    \caption{
    Burglaries per
\(\mathrm{km}^2\) per month \cite{Chicago_crime} (top left), filtered 
ISRs per \(\mathrm{km}^2\) per month \cite{ISR_Chicago} (top right), beat-wise constant police density in officers per \(\mathrm{km}^2\) \cite{OIG_Chicago} (bottom left), 
and crime-to-staffing mismatch (bottom right) in Chicago
between July and September 2025.
    }
    \label{fig:empirical_correspondence_q3_2025}
\end{figure}

\begin{figure}[htb!]
    \centering
    \includegraphics[width=0.9\textwidth]{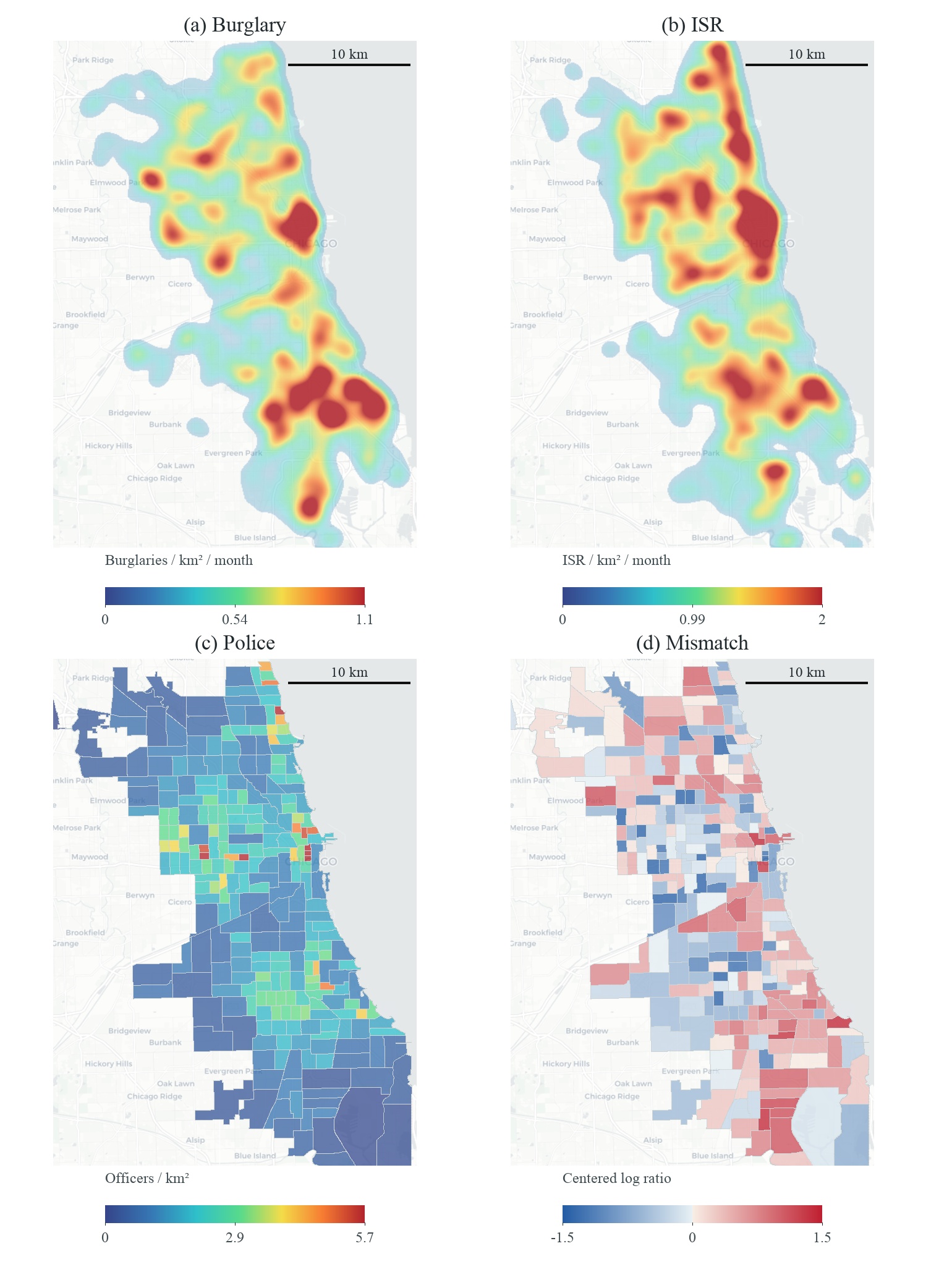}
    \caption{
   Burglaries per
\(\mathrm{km}^2\) per month \cite{Chicago_crime} (top left), filtered 
ISRs per \(\mathrm{km}^2\) per month \cite{ISR_Chicago} (top right), beat-wise constant police density in officers per \(\mathrm{km}^2\) \cite{OIG_Chicago} (bottom left), 
and crime-to-staffing mismatch (bottom right) in Chicago
between October and December 2025.
    }
    \label{fig:empirical_correspondence_q4_2025}
\end{figure}

A key feature of practical policing is the presence of delays in response to evolving crime conditions. Police deployment decisions rely on information processing, forecasting, and resource allocation, all of which introduce temporal gaps between the observation of crime and the intervention. This can be seen in data made publicly available by the city of Chicago, which appears to be the best U.S. case for spatio-temporal coupling of patrol presence and crime. 
See Appendix A to learn about the data sources and processing.
For context, the Chicago Police Department divides the city into 22 police districts, which are further divided into roughly 279 police beats. These small geographic areas are assigned dedicated patrol units.
Figs.~\ref{fig:empirical_correspondence_q3_2025} and 
\ref{fig:empirical_correspondence_q4_2025} report
burglaries per
\(\mathrm{km}^2\) per month, filtered Investigatory Stop Reports
(ISRs) per \(\mathrm{km}^2\) per month, beat-wise constant police density in officers per \(\mathrm{km}^2\), 
and crime-to-staffing mismatch for two consecutive quarters in 2025 (July-September and October-December, respectively). 
The crime-to-staffing mismatch is computed at the beat level as
\[
    M_i
    =
    \log\!\left(
        \frac{C_i+\alpha}{O_i+\beta}
    \right)
    -
    \operatorname{median}_j
    \log\!\left(
        \frac{C_j+\alpha}{O_j+\beta}
    \right),
\]
where \(C_i\) is average monthly burglary density in beat \(i\),
\(O_i\) is average officer density in beat \(i\), and
\(\alpha,\beta>0\) are small constants used to avoid instability in
low-count or low-staffing beats. A positive value of $M_i$ indicates higher burglary
activity relative to staffing than the citywide median, while a negative value
indicates lower burglary activity relative to staffing. Because the four panels in Figs.~\ref{fig:empirical_correspondence_q3_2025} and 
\ref{fig:empirical_correspondence_q4_2025}
represent different empirical quantities, color scales are normalized separately
within each panel and should be interpreted as relative spatial intensity.  

Figs.~\ref{fig:empirical_correspondence_q3_2025} and 
\ref{fig:empirical_correspondence_q4_2025} show burglary
hotspots that exhibit temporal evolution and spatial migration. This evolution cannot be captured by the models for burglary in \cite{short2008statistical,short2010dissipation,short2010nonlinear,rodriguez2010,zipkin2014}. Indeed, these models
give rise to static equilibria with 
hotspots, motivating the need for models with a time-dependent mechanism.
Filtered ISR activity exhibits partial spatial
overlap with burglary intensity, but it is more diffuse, meaning that there are
areas where police contacts are frequent even when burglary intensity is contained. 
This is expected since ISR activity is a proxy for a
police-observed or police-acted-upon field \cite{Chen2025}, rather than a direct measurement of reported crime. Moreover, it supports the modeling distinction between actual crime intensity
and perceived or processed crime signal available to police.
The bottom left panels in Figs.~\ref{fig:empirical_correspondence_q3_2025} and 
\ref{fig:empirical_correspondence_q4_2025} show that 
staffing density does not vary substantially 
from one quarter to the other, most likely due to rigidity in staffing allocation. 
In the mismatch panels, across both quarters we see
some high-burglary areas that appear under-resourced and
some well-staffed areas that experience relatively low burglary. Such misalignment is dynamic. We would like to stress that Figs.~\ref{fig:empirical_correspondence_q3_2025} and 
\ref{fig:empirical_correspondence_q4_2025}
should be read as qualitative empirical motivation for the models presented in this work, with no advocacy to increase staffing with the goal of reducing crime. 

Motivated by the need for more realistic models for police response and in light of the considerations drawn from Figs.~\ref{fig:empirical_correspondence_q3_2025} and 
\ref{fig:empirical_correspondence_q4_2025}, in this paper we extend the agent-based model 
for residential burglary in~\cite{short2008statistical,Hao2026} by incorporating a delayed police response. 
Specifically, we introduce an additional state variable representing delayed crime data and couple it to the equation for police deployment. 
We will show that delays in police response introduce memory effects that can fundamentally alter the dynamics of the system.
By taking the  mean-field limit of the agent-based model, we obtain the associated
continuum model, which
consists of a coupled system (three PDEs and one ODE) with highly nonlinear interactions.

To elucidate the mechanisms underlying hotspot formation and temporal variability in the proposed model, we conduct a linear stability analysis of spatially homogeneous steady states. 
By applying the Routh-Hurwitz criterion, we identify parameter regimes in which the homogeneous state is stable, as well as conditions under which delay induces oscillatory instabilities. The analysis reveals the emergence of Hopf bifurcations driven by the police response delay for given parameters and police population, leading to time-periodic average crime levels and dynamically evolving hotspots. 
To confirm the theoretical findings, 
we run a series of 
simulations with both the PDE and agent-based models, revealing a rich range of spatio-temporal dynamics. Finally, we
compare the proposed delayed-response mechanism with two alternative strategies 
considered in the literature, i.e., fixed police deployment~\cite{short2010nonlinear}
and globally optimal police patrolling 
\cite{zipkin2014}.
These comparisons reveal qualitative differences in the resulting crime patterns and highlight the trade-offs between responsiveness, stability, and resource allocation.

Before continuing any further, we would like to mention that we are aware that the use of mathematical models for predictive policing, i.e., to suggest law enforcement where crime is most likely to occur, has attracted criticism \cite{Haskins2019,Castelvecchi2020}.
Our intention is not to predict where the next crimes are likely to happen, 
but to develop a mathematical model 
and computational framework that can help assess the effectiveness of policing strategies.

The rest of the paper is organized as follows.  Sec.~\ref{sec:agent_model} and ~\ref{sec:pde_model} discuss the agent-based model with delayed police response and associated PDE model, with linear stability analysis. 
In Sec.~\ref{sec:num_met}, we present the time and space discretization of the PDE model and introduce an iterative partitioned algorithm to solve the discretized problem. Sec.~\ref{sec:num_res} reports the numerical results obtained with both models.  
Conclusions are drawn in Sec.~\ref{sec:concl}.

\section{The agent-based model }\label{sec:agent_model}

This section presents our discrete
statistical model for crime dynamics with police intervention. The basis for our model is the 
statistical model from \cite{short2008statistical}.

Let us consider a smooth and bounded  domain $\Omega\subset \mathbb{R}^2$ where there are three kinds of ``agents'': the policemen, the burglars, and the houses where the burglaries can occur. The houses are located on a two-dimensional lattice. 
For simplicity, we will consider a rectangular domain and 
consider a Cartesian lattice with uniform
spacing $h = \Delta x = \Delta y$.
Similarly, time is discretized with a time step $\Delta t$. Houses are located at the lattice
sites $s$ and are associated with a time-dependent
level of attractiveness $A_s(t)$, which measures the burglars’
perception of the attractiveness of the home at site $s$.
In general, we will use subindex $s$ to denote the variables in the discrete model that refer to site $s$.
Variable $A_s(t)$ is written as the sum of a static (possibly
spatially varying \cite{Hao2026}) component $A_s^{st}$ and a dynamic component $B_s(t)$:
\begin{equation}\label{eq:attr}
    A_s(t) = A_s^{st} + B_s(t).
\end{equation}
The dynamic component evolves based on repeat
and near-repeat victimization, broken-window effect, and the presence of law enforcement.

Let us denote by $m_s(t)$ the number of policemen at site $s$. 
We assume there is a fixed police population in $\Omega$, i.e., $\sum_s m_s = M$.
At a given site, criminals can burglarize
it or get stopped by police. Otherwise,
they move to a neighboring site with a random walk biased toward areas of high attractiveness.
The probability of each crime occurring at site $s$ during time interval $[t, t+\Delta t]$ is given by:
\begin{equation}\label{eq:prob}
    p^c_s(t)= 1-e^{-A_s(t)\exp{(-\beta m_s/h^2)}\Delta t},
\end{equation}
in accordance with a standard Poisson process. Notice that the presence of police makes crimes less likely to happen through a factor $\beta > 0$, which multiplies local police density $m_s/h^2$. 
When a site $s$ is burglarized, the corresponding criminal agent is removed from the lattice to represent
the tendency of burglars to flee the crime location. 
Obviously, removed burglars can become active again, which can be seen as a regeneration of criminals. To represent this in the model, a parameter $\Gamma$ is introduced at each lattice site $s$.
However, if a criminal has been arrested, we assume they will
not become active again for the entire time interval under consideration.

To capture repeat victimization, the attractiveness $B_s(t)$ at location $s$ 
is increased by $\theta$ for each crime 
committed at $s$ in the preceding time period \cite{short2008statistical,short2010dissipation}. Let parameter $\omega$ be 
a time scale over which repeat victimizations are most likely to occur
and $E_s(t)$ the number of burglaries occurred at site $s$ during the time
interval $(t - \Delta t, t)$. 
Thus, to account for repeat victimization
the equation for $B_s(t + \Delta t)$ needs to include term:
\begin{equation*}
    B_s(t)(1-\omega\Delta t)+\theta E_s(t).
\end{equation*}
To model near-repeat victimization and the broken-window effect, 
$B_s(t)$ has to spread spatially from a house to its neighbors. So, $B_s(t + \Delta t)$ needs to include a term of this kind:
\begin{equation*}
     \eta\bigg[\frac{1}{4}\sum_{r \sim s} B_{r}(t) - B_s(t) \bigg](1-\omega\Delta t),
\end{equation*}
where parameter $\eta \in [0, 1]$ characterizes the strength of the 
near-repeat victimization and the broken-window effect
and $r \sim s$ denotes the set of houses neighboring site $s$.
{For a Cartesian lattice, the cardinality of this set is four, 
which explains the factor of 4 in the formula above.} See, e.g., \cite{Hao2026} for more details.



Putting the three terms discussed above together, we get
\begin{equation}\label{eq:discrete_bertozzi-1}
    B_{s}(t+\Delta t) = \bigg[B_s(t)+\frac{\eta h^2}{4} \Delta_h B_{s}(t)\bigg](1-\omega\Delta t)+\theta E_s(t), 
\end{equation}
where
\begin{equation}\label{eq:Delta_Bst}
    \Delta_h B_{s}(t) = \frac{1}{h^2} \left(\sum_{s' \sim s} B_{s'}(t) - 4 B_{s}(t) \right)
\end{equation}
is the discrete spatial Laplacian operator. 

One way to estimate $E_s(t)$ in \eqref{eq:discrete_bertozzi-1} is given by:
\begin{equation}\label{eq:Es}
    E_s(t) = n_s(t) p^c_s(t),
\end{equation}
where $n_s (t)$ is the number of criminals at site $s$ at time $t$ and $p^c_s(t)$ is defined in \eqref{eq:prob}.

Next, we need to describe the evolution of $n_s (t)$.
Recall that criminals can burglarize (and thus be removed from the system)
or move to a neighboring site following a biased
random walk, with probability of moving to site $r$:
\begin{equation}\label{eq:move}
    q^c_{s \rightarrow r} = \frac{A_r (t)}{\displaystyle{\sum_{s' \sim s} A_{s'}(t)}}.
\end{equation}
Thus, at the following time step 
at site $s$ criminals arrive after either failing to burglarize a neighboring site or being generated at site $s$ at rate $(1 - \Sigma)\Gamma e^{-\beta \frac{m_s}{h^2}}$, where
$\Sigma$ is the probability that a burglary leads to an arrest. Note that the presence of police dampens the regeneration of criminals in the same way it dampens attractiveness in their probability to strike
\eqref{eq:prob} and the multiplication by factor $(1 - \Sigma)$ accounts for the fact that arrested burglars will not become active again in the time interval under consideration.
Putting together the terms that give
movement and generation of criminals,
we obtain the following: 
\begin{equation}\label{eq:n_s_P0}
    n_s(t+\Delta t) =  \sum_{r\sim s} n_r(t)\Big(1-p_r^c(t)\Big)q_{r \rightarrow s}^c(t) +e^{-\beta \frac{m_s}{h^2}}\Gamma (1 - \Sigma) \Delta t,
\end{equation}
which can be rewritten as
\begin{equation}\label{eq:n_s_P}
      n_s(t+\Delta t) = A_s \sum_{r\sim s} \frac{n_r(t)[1-p^c_r(t)]}{T_{r}(t)}+ e^{-\beta \frac{m_s}{h^2}}\Gamma (1 - \Sigma) \Delta t, \quad 
    T_{r}(t) = \sum_{s'\sim r}A_{s'}(t).
\end{equation}
To prevent the denominator of $q_{s \rightarrow r}$ from becoming zero,
we assume $A_s(t) \geq a_0 > 0$ at all times $t$ under consideration and for all sites $s$ in the lattice.

In summary, the effect of police presence in our model is twofold. 
First, it makes it less likely for a criminal to strike, though 
\eqref{eq:prob}. Second, it reduces the number of criminals
circulating in the region, through the last term in \eqref{eq:n_s_P}.

It is left to describe how policemen patrol. Following the same idea
used for criminals, we assume that police are either busy at site $s$, e.g., arresting someone, or they move according to 
a biased random walk towards areas with higher crime density. 
The expected number of crimes per unit area per unit time at site $s$ and time $t$ is given by
\begin{align}
      S_s(t) &= \frac{E_s(t)}{{h^2} \Delta t} = \frac{n_s(t) p^c_s(t)}{{h^2} \Delta t}, \label{eq:S_dis}
\end{align}
Obviously, we cannot expect police to obtain and share this kind of data in real time to decide where to patrol. Instead, it is reasonable to assume these data
will become available and used for decision-making with some delay, denoted with $\tau$. Hence, we introduce a first-order lag:
\begin{align}
      H_s(t+\Delta t) &= \left(1-\frac{\Delta t}{\tau} \right)H_s(t) + \frac{\Delta t}{\tau} S_s(t). \label{eq:H}
\end{align}
It is $H_s(t)$ that will determine where a policeman moves to at time $t$. Thus, we define the per-step transition probability from site $s$ to site $r$ for a policeman as
\begin{equation}\label{eq:moveP}
   q_{s \rightarrow r}^p = \frac{H_r (t)}{\displaystyle{\sum_{s' \sim s} H_{s'}(t)}}. 
\end{equation}
Since $H_s$ is simply $S_s$ in \eqref{eq:S_dis} lagged in time
and $S_s$ is strictly positive, 
the denominator of \eqref{eq:moveP}
is guaranteed to be strictly positive.




To write the equivalent of \eqref{eq:n_s_P0} for police, we
introduce the probability an officer makes at least one arrest in time interval of length $\Delta t$, denoted with $p^p_s(t)$. If incidents that could lead to arrests occur as a Poisson process with rate $n_s(t)p_s^c(t)/\Delta t$ (incidents per unit time) and each incident leads to an arrest independently with probability $\Sigma$, then arrests themselves form a Poisson process with rate 
$\Sigma n_s(t)p_s^c(t)/\Delta t$. So, we have:
\begin{equation}\label{eq:psP}
    p_s^p(t) = 1 - e^{- \Sigma n_s(t)p_s^c(t)},
\end{equation}
and police move according to the following equation:
\begin{equation}\label{eq:ms2}
    m_s(t+\Delta t) =  \sum_{r\sim s} m_r(t)[1-p_r^p(t)]q_{r \rightarrow s}^p(t) +m_s(t)p_s^p(t),
\end{equation}
where $m_s(t)p^p_s(t)$ is the number of policemen busy arresting criminals in time interval of length $\Delta t$.
Notice that eq.~\eqref{eq:ms2} resembles \eqref{eq:n_s_P0}. However, 
unlike \eqref{eq:n_s_P0}, there is no generation or removal term and so the total police population is conserved. In fact, we have:
\begin{align*}
        \sum_s m_s(t + \Delta t) &= \sum_s\sum_{r\sim s} m_r(t)[1-p^p_r(t)]q^p_{r\rightarrow s}(t)+\sum_s m_s(t) p^p_s(t)\\
        &= \sum_r m_r(t) [1-p^p_r(t)]\sum_{s\sim r} q^p_{r\rightarrow s}(t)+\sum_s m_s(t) p_s^p(t)\\
        & = \sum_r m_r(t) [1-p_r^p(t)] +\sum_s m_s(t) p^p_s(t) =\sum_s m_s(t).
\end{align*}
By plugging \eqref{eq:moveP} into eq.~\eqref{eq:ms2}, we get: 
\begin{equation}\label{eq:ms3}
      m_s(t+\Delta t) =  H_s \sum_{r\sim s} \frac{m_r(t)[1-p^p_r(t)]}{V_{r}(t)}+ {m_s(t)p^p_s(t)}, \quad 
    V_{r}(t) = \sum_{s'\sim r}H_{s'}(t).
\end{equation}

Problem \eqref{eq:discrete_bertozzi-1},\eqref{eq:n_s_P},\eqref{eq:H},\eqref{eq:ms3}, with $E_s(t)$ as in \eqref{eq:Es}, $p^c_s(t)$ as in \eqref{eq:prob}, $S_s(t)$ as in \eqref{eq:S_dis}, and $p^p_s(t)$ as in \eqref{eq:psP}, represents the probabilistic agent-based model. 
We consider a finite domain and natural boundary conditions, meaning that at the boundary the normal derivatives of the system variables {are} set to zero. In the agent-based model, we enforce such boundary conditions
{for the level of attractiveness, the number of criminals, and the number of policemen. 
See \cite{Hao2026} for details on how to impose these boundary conditions.


Let $N_s$ be the total number of sites in the lattice. 
When all sites have the same level of attractiveness $\bar{A} = A^{st}+\bar{B}$ and, on average, the same number of criminals $\bar{n}$ and the same number of policemen $\bar{m}$, the system has a homogeneous equilibrium
solution with:
\begin{equation}\label{eq:equil_sol}
    \bar{B} =\frac{\theta \Gamma}{\omega} (1 - \Sigma)e^{-\beta\frac{\bar{m}}{h^2}}, \quad 
    \bar{n} = \frac{{\Delta t} \Gamma }{\bar{p}_c} (1 - \Sigma)e^{-\beta\frac{\bar{m}}{h^2}}, \quad
    \bar{m} = \frac{M}{N_s}, \quad
    \quad \bar{H} = \frac{\bar{n}\bar{p}^c}{h^2\Delta t},
\end{equation}
where
\begin{align*}
    \bar{p}_c = 1-e^{-\bar{A}\exp{(-\beta \bar{m}/h^2)}\Delta t}.
\end{align*}


\section{The PDE model}\label{sec:pde_model}

Let $\rho = n_s(t)/h^2$ be the density
of criminals and $\pi = m_s(t)/h^2$ be the density of policemen.
By manipulating eq.~\eqref{eq:discrete_bertozzi-1}, with $E_s(t)$ as in \eqref{eq:Es}, and taking its limit for $h, \Delta t \rightarrow 0$, while keeping the ratio $\frac{h^2}{\Delta t}$ fixed to $D$ and quantity $\theta \Delta t$ fixed to $\epsilon$, we obtain:
\begin{equation}\label{eq:contB_dim}
    \frac{\partial B}{\partial t}  =  \frac{\eta D}{4} \Delta B - \omega B + \epsilon D \rho A e^{-\beta \pi}, \quad A(\bx, t) = A^{st}(\bx) + B(\bx, t),
\end{equation}
where $\bx$ gives the spatial coordinates. In \eqref{eq:contB_dim},
we have removed the subindex $s$ {from} the variables since the model is not discrete anymore.

Similarly, by manipulating eq.~\eqref{eq:n_s_P} and taking its limit for $h, \Delta t \rightarrow 0$, while keeping the ratio $\frac{\Gamma}{h^2}$ fixed to $\gamma$, we get:
\begin{equation}\label{eq:rho_P2}
     \frac{\partial \rho}{\partial t}= \frac{D}{4} \nabla \cdot \Big(\nabla\rho - \frac{2 \nabla A}{A} \rho \Big)-\rho A e^{-\beta \pi} +\gamma (1 - \Sigma) e^{-\beta \pi}.
\end{equation}
To reflect the strict positivity requirement on $A_s(t)$ stated after eq.~\eqref{eq:n_s_P}, we assume 
$A(\bx, t) \geq a_0 > 0$ for $\bx \in \Omega$ and $t$ in the time interval of interest.

If we take the limit for $h, \Delta t \rightarrow 0$ of eq.~\eqref{eq:ms2}, we get:
\begin{align}
\frac{\partial \pi}{\partial t} = \frac{D}{4}  \nabla \cdot \Big(\nabla \pi - \frac{2 \nabla H}{H} \pi \Big).\label{eq:p_P}
\end{align}
Next, we complete the PDE model and
comment on the positivity of $H$.

We need to take the hydrodynamic limit of \eqref{eq:H}. For this, we note that as $h, \Delta t \rightarrow 0$:
\begin{align}
      S_s(t) 
      \rightarrow S = \rho A e^{-\beta \pi}. \label{eq:S}
\end{align}
Thus, in the limit eq.~\eqref{eq:H}
becomes
 \begin{align}
      \frac{\partial H}{\partial t}= \frac{1}{\tau} \left(\rho A e^{-\beta \pi} - H \right). \label{eq:H2}
\end{align}
Since $H$ is $S$ lagged in time, the strict positivity of $A$ and $\rho$ imply the 
the strict positivity of $H$.

To simplify problem \eqref{eq:contB_dim}-\eqref{eq:p_P},\eqref{eq:H2}, we introduce the following non-dimensional quantities:
\begin{equation}\label{eq:non_dim}
    \tilde{\bx} = 2 \frac{\bx}{L}, \quad \tilde{t} = \omega t, \quad \tilde{A} = \frac{A}{\omega} = \frac{A^{st}}{\omega} + \frac{B}{\omega},\quad \tilde{\rho} = \epsilon L^2 \rho, \quad \tilde{\pi} = \beta \pi, \quad
    \tilde{H} = \frac{\epsilon L^2}{\omega} H, \quad \tilde{\tau} = \omega \tau,
\end{equation}
where 
$L = \sqrt{{D}/{\omega}}$ is a characteristic length. 
We note that the 2 in the first equation in \eqref{eq:non_dim}
comes from $\sqrt{4}$, where 4 is the number of neighbors of a given site. 
For simplicity of notation, we will omit the tilde on the non-dimensional variables and from now on every variable in the continuum (or PDE) model will be non-dimensional.

The non-dimensional continuum model can be written as: find $A(\bx, t)$, $\rho(\bx, t)$, $\pi(\bx, t)$, and $H(\bx, t)$, such that 
\begin{align}
     & \frac{\partial A}{\partial t} - \eta \Delta A + A - \rho A e^{-\pi}= - \eta \Delta A^{st} + A^{st}, \label{eq:continous_1}\\
     &\frac{\partial \rho}{\partial t} - \nabla \cdot \Big(\nabla\rho - \frac{2 \nabla A}{A} \rho\Big)+\rho A e^{-\pi}  - \frac{\Gamma \theta}{\omega^2} (1 - \Sigma) e^{-\pi}= 0,\label{eq:continous_2} \\
     &\frac{\partial \pi}{\partial t} - \nabla \cdot\Big(\nabla \pi -\frac{2 \nabla H}{H} \pi \Big)=0, \label{eq:continous_3} \\
     &\frac{\partial H}{\partial t} - \frac{1}{\tau} \Big(\rho A e^{- \pi} - H \Big) = 0. \label{eq:continous_4}
\end{align}
We note that the above system
features 3 parameters: $\eta$, $\Gamma \theta (1-\Sigma)/{\omega^2} $, and $\tau$. 

Problem \eqref{eq:continous_1}-\eqref{eq:continous_4} needs to be supplemented with appropriate boundary and initial conditions. In line with our
previous work \cite{Hao2026}, we enforce 
natural boundary conditions
\begin{align}
   \nabla A {\cdot} \bn &= 0,  \quad \text{on } \partial\Omega,\; t>0, \label{eq:BC1} \\
   \nabla \rho {\cdot} \bn &= 0, \quad \text{on } \partial\Omega,\; t>0, \label{eq:BC2} \\
   \nabla \pi {\cdot} \bn &= 0, \quad \text{on } \partial\Omega,\; t>0, \label{eq:BC3}
\end{align}
where $\bn$ is the outward normal unit vector. 
As for the initial conditions, we set:
\begin{align}
     A(\bm{x}, 0) = A^{st}(\bm{x}) +  B_0(\bm{x}),\quad
\rho(\bm{x}, 0) = \rho_0(\bm{x}), \quad
\pi(\bm{x}, 0) =
\pi_0(\bm{x}), \label{eq:ic} 
\end{align}
where 
$A^{st}(\bm{x})+B_0(\bm{x})$ is the prescribed initial level of attractiveness,
$\rho_0(\bm{x})$ is the initial criminal density, and $\pi_0(\bm{x})$ is the initial police density. As for the initial condition for \eqref{eq:continous_4}, since $H$ is delayed $S$, we set
\begin{align}
    H(\bm{x}, 0) = \rho(\bm{x}, 0)A(\bm{x}, 0) e^{- \pi(\bm{x}, 0)}.\label{eq:ic4}
\end{align}

\begin{remark}\label{rem:pol_strategies}



If one chooses to prescribe $\pi$ in the above model, 
eq.~\eqref{eq:continous_3} and \eqref{eq:continous_4}
are discarded from the model. In particular, if one sets $\Sigma = 0$ and prescribes $\pi$ at and after a given time $t_s > 0$ (i.e., police are deployed after hotspots are formed) as
\begin{align}
    \pi (\bm{x}) = -\log\left( \frac{1}{2} \left( 1 - \tanh{\mu \left( A(\bm{x}, t_s) - A_c\right) } \right) \right), \label{eq:fixed_pol}
\end{align}
for given positive parameters $\mu$ and $A_c$, model \eqref{eq:continous_1}-\eqref{eq:continous_2} coincides with the model from \cite{short2010nonlinear,short2010dissipation}. 
Instead, if one prescribes $\pi$ according to the following minimization problem: 
\begin{align*}
    \pi(\bm{x},t)  = \argmin \left\{ \int_\Omega e^{-\pi} \rho A d\bm{x} : \pi \in L^1(\Omega),  \pi \geq 0,  \int_\Omega {\pi} d\bm{x} = M \right\},
\end{align*}
model \eqref{eq:continous_1}-\eqref{eq:continous_2} coincides with the model from \cite{zipkin2014}, which assumes that police are deployed to minimize the total crime occurring instantaneously. Thus, the models in 
\cite{short2010nonlinear,short2010dissipation,zipkin2014} can be seen as particular cases of model
\eqref{eq:continous_1}-\eqref{eq:continous_4}. 
In the model in \cite{jones2010}, police agents perform a biased random walk toward high-attractiveness sites, however without 
information delay. This is structurally related to, though distinct from, 
the $\tau \to 0^+$ limit of our model.
\end{remark}

\begin{remark}
Local existence and uniqueness of solutions of the model in 
\cite{short2008statistical} (i.e., the original model with no police) is demonstrated in \cite{rodriguez2010}. As mentioned in the previous remark, the model with police in \cite{short2010nonlinear,short2010dissipation} is a simple modification of the model in \cite{short2008statistical}, hence the local existence and uniqueness result holds for that model too. We believe that, leveraging the work in 
\cite{rodriguez2010}, one could prove local existence and uniqueness of solutions of model 
\eqref{eq:continous_1}-\eqref{eq:continous_4} too. However, this is beyond the scope of this work.
\end{remark}

\subsection{Homogeneous equilibrium solution and linear stability}\label{sec:lin_stab}
We denote by $(\bar{A}, \bar{\rho}, \bar{\pi}, \bar{H})$ a spatially homogeneous equilibrium solution of
\eqref{eq:continous_1}-\eqref{eq:continous_4}.
For homogeneous  $H$, eq.~\eqref{eq:continous_3} is satisfied by any $\pi$ homogeneous in space and time. Also, since eq.~\eqref{eq:continous_3} has no source/sink term, given an initial homogeneous in space $\bar{\pi}$ the steady state equilibrium will not move away from $\bar{\pi}$. Then, 
assuming $\bar{\pi}$ is given and
$A^{st}$ is homogeneous in space, solution $(\bar{A},\bar{\rho},\bar{H})$ satisfies:
\begin{align}
     & \bar{A} - \bar{\rho} \bar{A} e^{- \bar{\pi}} = A^{st}, \label{eq:continous_1ss} \\
     &\bar{\rho} \bar{A} e^{- \bar{\pi}} - \frac{\Gamma \theta}{\omega^2} (1-\Sigma)  e^{- \bar{\pi}} = 0, \label{eq:continous_2ss} \\
     &\bar{\rho} \bar{A} e^{- \bar{\pi}} - \bar{H} = 0. \label{eq:continous_3ss}
\end{align}
By manipulating \eqref{eq:continous_1ss}-\eqref{eq:continous_3ss}, one obtains:
\begin{align}
\bar{A} = A^{st} + \frac{\Gamma \theta}{\omega^2} (1-\Sigma)  e^{- \bar{\pi}}, \quad \bar{\rho} =  \frac{\Gamma \theta}{\omega^2}  \frac{(1-\Sigma)}{\bar{A}}, \quad
\bar{H} = \bar{\rho} \bar{A} e^{- \bar{\pi}}, \label{eq:hom_eq_sol}
\end{align}
which can be easily computed once $\frac{\Gamma \theta}{\omega^2} (1 - \Sigma)$ is set and $\bar{\pi}$ is given. We remark that this homogeneous equilibrium solution is nondimensionalized according to \eqref{eq:non_dim}, while the 
equilibrium solution in \eqref{eq:equil_sol} is dimensional.

Since we impose homogeneous Neumann boundary conditions \eqref{eq:BC1}-\eqref{eq:BC3}, we consider perturbations in the eigenfunctions 
$\{\varphi_j\}_{j\ge0}$ of the Neumann Laplacian: 
\begin{equation}\label{eq:neu_lapl}
    -\Delta\varphi_j=\mu_j\varphi_j,
\qquad
0=\mu_0<\mu_1\le\mu_2\le\cdots.
\end{equation}
For each mode, we perturb the spatially homogeneous equilibrium as follows: 
\begin{align}
    A(\bm{x}, t) &= \bar{A} + \delta_A e^{\lambda t} \varphi_j(x), \label{eq:pert_A} \\
    \rho(\bm{x}, t) &= \bar{\rho} + \delta_\rho e^{\lambda t} \varphi_j(x), \label{eq:pert_rho} \\
    \pi(\bm{x}, t) &= \bar{\pi} + \delta_\pi e^{\lambda t} \varphi_j(x), \label{eq:pert_pi} \\
    H(\bm{x}, t) &= \bar{H} + \delta_H e^{\lambda t} \varphi_j(x). \label{eq:pert_H}
\end{align}

By plugging \eqref{eq:pert_A}-\eqref{eq:pert_H} into \eqref{eq:continous_1}-\eqref{eq:continous_4}, we obtain the following linearized system:
\begin{align}\label{eq:eig}
   & J
\begin{bmatrix}
    \delta_A \\
    \delta_\rho \\
    \delta_\pi \\
    \delta_H
\end{bmatrix} 
= 
\lambda \begin{bmatrix}
    \delta_A \\
    \delta_\rho \\
    \delta_\pi \\
    \delta_H
\end{bmatrix}, \quad
J = \begin{bmatrix}
  -\eta \mu_j - 1 + \bar{\rho} e^{- \bar{\pi}} & \bar{A} e^{- \bar{\pi}} & -\bar{\rho}\bar{A}e^{-\bar{\pi}} & 0 \\
  2 \mu_j \frac{\bar{\rho}}{\bar{A}} -  \bar{\rho} e^{- \bar{\pi}} &  - \mu_j  -  \bar{A} e^{- \bar{\pi}} & 0 & 0 \\
  0 & 0 & - \mu_j & 2 \mu_j \frac{\bar{\pi}}{\bar{H}} \\
  \frac{1}{\tau} \bar{\rho} e^{- \bar{\pi}} & \frac{1}{\tau} \bar{A} e^{- \bar{\pi}} & -\frac{1}{\tau}\bar{\rho}\bar{A}e^{-\bar{\pi}} & -\frac{1}{\tau}
\end{bmatrix}.
\end{align}
Using the fact that $\bar\rho\bar{A}e^{-\bar\pi} = \bar{H}$, we can derive the characteristic polynomial for the matrix $J$ in \eqref{eq:eig} by co-factor expansion along the third row and set:
\[
\begin{aligned}
\det(\lambda I - J)= &(\mu_j + \lambda)\left(\frac{1}{\tau} + \lambda\right)\Big[(-\eta\mu_j - 1 + \bar\rho e^{-\bar\pi} - \lambda)(-\mu_j - \bar{A}e^{-\bar\pi} - \lambda) - e^{-\bar\pi}(2\mu_j\bar\rho - \bar{H})\Big] \\
&+ 2\mu_j\frac{\bar\pi}{\tau}(-\eta\mu_j - 1 - \lambda)(-\mu_j - \bar{A}e^{-\bar\pi} - \lambda) = 0.
\end{aligned}
\]
Finding the analytic dependence of the eigenvalues on $\bar\rho, \bar{A}, \bar\pi, \eta$ and $\mu_j$ is too involved from the algebraic point of view
and beyond the scope of this paper. Instead, we will proceed with the Routh-Hurwitz criterion.

The characteristic polynomial of $J$ can be written as
\begin{equation}\label{eq:char_poly}
P(\lambda) = \lambda^4 + a_3\,\lambda^3 + a_2\,\lambda^2
           + a_1\,\lambda + a_0,
\end{equation}
where the coefficients depend on $\mu_j$, $\tau$, and the equilibrium
values in~\eqref{eq:hom_eq_sol}.
See Appendix B.
 
By the Routh--Hurwitz criterion, all roots of~\eqref{eq:char_poly}
have negative real part if and only if $a_i > 0$ for $i=0,\ldots,3$
and the Hurwitz determinant
\begin{equation}\label{eq:hurwitz}
H_3 = a_1\,a_2\,a_3 - a_1^2 - a_0\,a_3^2
\end{equation}
is positive. A Hopf bifurcation occurs when $H_3 = 0$ with all
$a_i > 0$, at which point a complex-conjugate pair
$\lambda = \pm\,i\omega_0$ crosses the imaginary axis with
non-zero speed. 
Substituting $\lambda = i\omega_0$ into~\eqref{eq:char_poly}
and separating real and imaginary parts yields the oscillation
frequency
\begin{equation}\label{eq:hopf_freq}
\omega_0 = \sqrt{a_1/a_3}.
\end{equation}
 
On a square domain $\Omega = [0,L] \times [0,L]$, for given length $L$, with homogeneous Neumann boundary conditions, the eigenvalues of 
\eqref{eq:neu_lapl} are $\mu_j = (m^2+n^2)\pi^2/L^2$ for non-negative integers $m,n$.
Thus, 
for each given $\mu_j$ 
and for given equilibrium
values in~\eqref{eq:hom_eq_sol},
we can locate the critical value of $\tau$, denoted with $\tau_c(\mu_j)$, at which the Hopf bifurcation occurs by solving
$H_3(\mu_j,\tau)=0$ numerically.
The transversality condition
$\partial_\tau \mathrm{Re}(\lambda)\big|_{\tau=\tau_c} \neq 0$
can also be verified numerically and we will do so for one case in Sec.~\ref{sec:num_res}. 
We define
the critical delay as 
\begin{equation}\label{eq:tau_c}
\tau_c^* = \min_j\,\tau_c(\mu_j).
\end{equation}
We will denote with $\mu^*$ the value of $\mu_j$ for which the minimum is attained.
 
Fig.~\ref{fig:phase_diagrams} maps $\tau_c^*$ as a function of
$\eta$, $\Gamma\theta(1-\Sigma)/\omega^2$, and~$\pi_0$. The left
panel shows that $\tau_c^*$ increases monotonically with~$\eta$:
stronger 
neighborhood effects (i.e., larger $\eta$) stabilize the system
by smoothing local crime concentrations.
The center panel shows that $\tau_c^*$ is less sensitive to
$\Gamma\theta(1-\Sigma)/\omega^2$.
The right panel reveals that $\tau_c^* \to \infty$ as $\pi_0 \to 0$.
Fig.~\ref{fig:phase_diagrams} suggests that the Hopf bifurcation
mechanism is governed primarily by spatial diffusion (i.e., $\eta$) and police
density (i.e., $\pi_0$), rather than by the criminal ``regeneration''/suppression rate (i.e., $\Gamma\theta(1-\Sigma)/\omega^2$).

\begin{figure}[htb!]
\centering
\includegraphics[width=\textwidth]{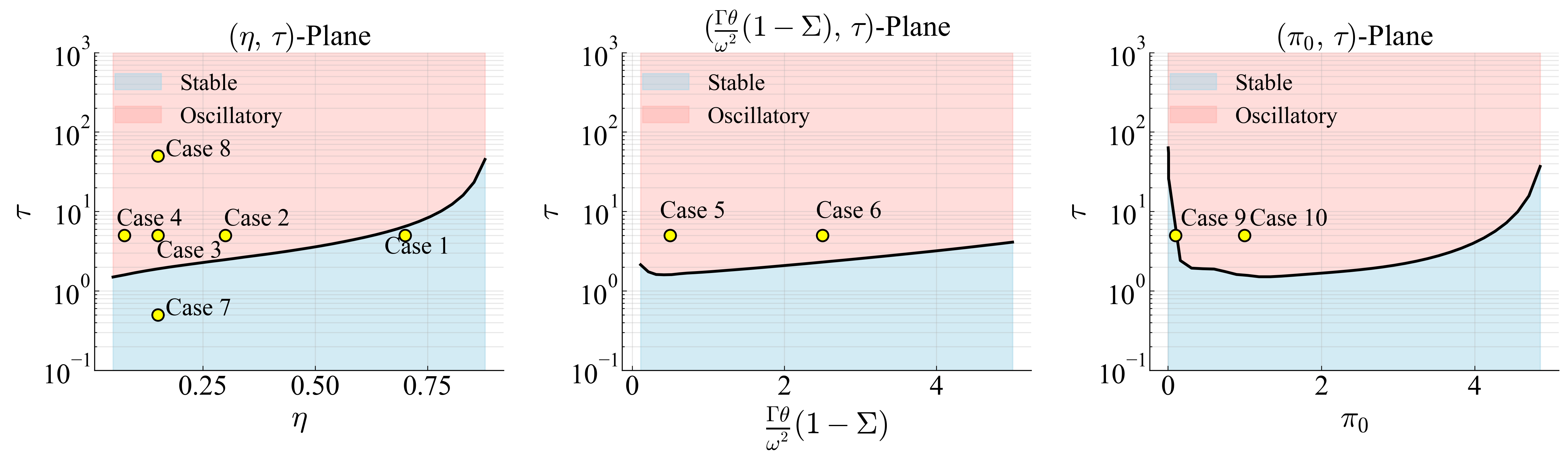}
\caption{Stability phase diagrams from the Routh-Hurwitz criterion with the black
curve giving the Hopf threshold $\tau_c^*$ as a function of
$\eta$ for $\Gamma\theta(1-\Sigma)/\omega^2 = 1.5, \pi_0 =0.5$ (left), $\tau_c^*$ as a function of $\Gamma\theta(1-\Sigma)/\omega^2$ for $\eta = 0.15, \pi_0 =0.5$ (center), and $\tau_c^*$ as a function of $\pi_0$ for $\eta =0.15, \Gamma\theta(1-\Sigma)/\omega^2 = 1.5$ (right). The stable part of the plane is colored in blue, while the oscillatory part is colored in pink.
The yellow dots represent the cases discussed in Sec.~\ref{sec:num_res}.}
\label{fig:phase_diagrams}
\end{figure}

For an explanation of the $\pi_0 \to 0$ limit, we refer to the Appendix B. 
The U-shaped profile in the right panel of Fig.~\ref{fig:phase_diagrams} has a local minimum near $\pi_0\approx 0.3$,
identifying intermediate police density as the parameter regime
most susceptible to delay-induced instability.

\subsection{Variational formulation}
For the numerical method that will be discussed in Sec.~\ref{sec:num_met}, 
we need a variational formulation of problem
\eqref{eq:continous_1}-\eqref{eq:continous_3}.
We leave \eqref{eq:continous_4}
as it is because it is an ODE, simply stating that $H$ is $S = \rho A e^{-\beta \pi}$ lagged in time.
To devise the variational problem, one multiplies \eqref{eq:continous_1} by $u\in H^1(\Omega)$, \eqref{eq:continous_2} by $v\in H^1(\Omega)$, and \eqref{eq:continous_3} by $w\in H^1(\Omega)$, 
integrates over $\Omega$ and employs the integration by parts identity. 
Taking into account boundary 
conditions \eqref{eq:BC1}-\eqref{eq:BC3}, and the fact that 
at the boundary $\partial \Omega$ we have  $\nabla H {\cdot} \bn = \nabla (\rho A e^{-\pi}) {\cdot} \bn = 0$,
we obtain the following weak problem: Find $A(\bx, t)\in H^1(\Omega) {\cap
L^\infty(\Omega)}$ and $\rho(\bx, t), \pi(\bx, t) \in H^1(\Omega)$ 
such that
\begin{align}
        &\int_\Omega \frac{\partial A}{\partial t} u\;d\Omega + \int_\Omega \eta\nabla A \cdot \nabla u\;d\Omega 
        +\int_\Omega Au \;d\Omega  
        - \int_\Omega \rho A e^{-\pi} u \;d\Omega  = \int_\Omega (A^{st} - \eta \Delta A^{st}) u \;d\Omega , \label{eq:weak_1} \\
        & \int_\Omega \frac{\partial \rho}{\partial t} v \;d\Omega +\int_\Omega  \Big(\nabla\rho - \frac{2 \nabla A}{A} \rho \Big) \nabla v\;d\Omega 
        +\int_\Omega \rho A e^{-\pi}  v\;d\Omega - \int_\Omega \frac{\Gamma \theta}{\omega^2} (1 - \Sigma) e^{-\pi} v\;d\Omega=  0,  \label{eq:weak_2} \\
        &\int_\Omega \frac{\partial \pi}{\partial t} w \;d\Omega +\int_\Omega  \Big(\nabla\pi - \frac{2 \nabla H}{H} \pi \Big) \nabla w\;d\Omega 
         = 0, \label{eq:weak_3}
\end{align}
and \eqref{eq:continous_4} hold
for all $(u,v,w)\in H^1(\Omega)\times H^1(\Omega)\times H^1(\Omega)$.

\section{Numerical method for the PDE model}\label{sec:num_met}


Let us start from the time discretization. For this, 
we introduce a time step $\Delta t = T/N_t$ for given $N_t$, set $t^n =n\Delta t $ for $n = 1, \dots, N_t$
and denote with $y^n$ the approximation of quantity 
$y$ at time $t^n$. 
For simplicity, we apply Backward Differentiation Formula
of order 1 to \eqref{eq:continous_4}, \eqref{eq:weak_1}-\eqref{eq:weak_3} and get the following problem at the generic time $t^{n+1}$: for $n\geq 0$, 
given $A^n$, $\rho^n$, $\pi^n$, and $H^n$, find $A^{n+1}\in H^1(\Omega){\cap
L^\infty(\Omega)}$, $\rho^{n+1}, \pi^{n+1} \in H^1(\Omega) $, and $H^{n+1}\in H^1(\Omega){\cap
L^\infty(\Omega)}$ 
such that
\begin{align}
           & \frac{1}{\Delta t}\int_\Omega A^{n+1} u\;d\Omega + 
           \int_\Omega \eta\nabla A^{n+1} \cdot \nabla u\;d\Omega +\int_\Omega A^{n+1}u \;d\Omega \cl
           & \quad \quad -\int_\Omega \rho^{n+1} A^{n+1} e^{- \pi^{n+1}} u \;d\Omega  = \int_\Omega (A^{st} - \eta \Delta A^{st} ) u \;d\Omega + \frac{1}{\Delta t}\int_\Omega A^{n} u\;d\Omega, \label{eq:implicit-1} \\
         & \frac{1}{\Delta t} \int_\Omega \rho^{n+1} v \;d\Omega +\int_\Omega  \Big(\nabla\rho^{n+1} - \frac{2\nabla A^{n+1}}{A^{n+1}} \rho^{n+1}\Big) \nabla v\;d\Omega  \cl
         & \quad \quad +\int_\Omega \rho^{n+1} A^{n+1} e^{-\pi^{n+1}}  v\;d\Omega - \int_\Omega \frac{\Gamma \theta}{\omega^2} (1 - \Sigma) e^{-\pi^{n+1}} v\;d\Omega = \frac{1}{\Delta t} \int_\Omega \rho^{n} v \;d\Omega, \label{eq:implicit-2} \\
         & \frac{1}{\Delta t} \int_\Omega \pi^{n+1} w \;d\Omega +\int_\Omega  \Big(\nabla\pi^{n+1} - \frac{2\nabla H^{n+1}}{H^{n+1}} \pi^{n+1}\Big) \nabla w\;d\Omega =  \frac{1}{\Delta t} \int_\Omega \pi^{n} w \;d\Omega, \label{eq:implicit-3} \\
         & \frac{1}{\Delta t} H^{n+1}   -  \frac{1}{\tau} \Big(\rho^{n+1} A^{n+1} e^{- \pi^{n+1}} - H^{n+1} \Big)   = \frac{1}{\Delta t}  H^{n}, \label{eq:implicit-4}
\end{align}
for all $(u,v,w)\in H^1(\Omega)\times H^1(\Omega)\times H^1(\Omega)$.
System \eqref{eq:implicit-1}-\eqref{eq:implicit-4} is
implicit and coupled. Hence, if one adopts a fine mesh for space discretization, it could lead to 
a large nonlinear system that needs to be solved
at each time step through an
iterative solver like, e.g., Newton’s method. Since 
problem \eqref{eq:implicit-1}-\eqref{eq:implicit-4} features several strongly nonlinear terms, the convergence of such 
iterative method could be slow, especially for large-scale simulations \cite{quarteroni1994numerical,saad2003iterative}. 

One way to contain the computational cost is to decouple system \eqref{eq:implicit-1}-\eqref{eq:implicit-4}, while linearizing it, so that one solves two smaller linear systems, instead of one large nonlinear system, per time step. 
While successfully containing the computational time, such a scheme may suffer
from instability for certain combination of parameter values. 
The following algorithm is 
a way to keep the advantages of a decoupled scheme without possibly incurring into instabilities.
At every time step $t^{n+1}$, given $A^n$, $\rho^n$, $\pi^n$, and $H^n$, perform at iteration $k+1$, $k \ge 0$:
\begin{itemize}
    \item[-] \emph{Step 1}:  find $A^{k+1}\in H^1(\Omega){\cap
L^\infty(\Omega)}$ such that
    \begin{align}
           &\frac{1}{\Delta t}\int_\Omega A^{k+1} u\;d\Omega + 
           \int_\Omega \eta\nabla A^{k+1} \cdot \nabla u\;d\Omega +\int_\Omega A^{k+1}u \;d\Omega \cl
           & \quad \quad -\int_\Omega \rho^{k} {e^{- \pi^{k}}} A^{k+1} u \;d\Omega  = \int_\Omega (A^{st} - \eta \Delta A^{st}) u \;d\Omega + \frac{1}{\Delta t}\int_\Omega A^{n} u\;d\Omega, \label{eq:iterative-1} 
\end{align}
for all $u \in H^1(\Omega)$. For $k = 0$, 
$\rho^{k} = \rho^{n}$.
 \item[-] \emph{Step 2}:  find $\rho^{k+1}\in H^1(\Omega)$ such that
 \begin{align}
          & \frac{1}{\Delta t} \int_\Omega \rho^{k+1} v \;d\Omega +\int_\Omega  \Big(\nabla\rho^{k+1} - \frac{2\nabla A^{k+1}}{A^{k+1}} \rho^{k+1}\Big) \nabla v\;d\Omega  \cl
         & \quad \quad +\int_\Omega \rho^{k+1}  A^{k+1} {e^{-\pi^{k}}}  v\;d\Omega = \int_\Omega \frac{\Gamma \theta}{\omega^2} (1 - \Sigma) {e^{-\pi^{k}}} v\;d\Omega + \frac{1}{\Delta t} \int_\Omega \rho^{n} v \;d\Omega, \label{eq:iterative-2}
\end{align}
for all $v \in H^1(\Omega)$.
\item[-] \emph{Step 3}:  find $H^{k+1}\in L^2(\Omega){\cap
L^\infty(\Omega)}$ such that
\begin{align}
 \left( \frac{1}{\Delta t} + \frac{1}{\tau} \right)  H^{k+1}    =  \frac{1}{\tau} \Big(\rho^{k+1} A^{k+1} e^{- \pi^{k}} \Big)  +  \frac{1}{\Delta t}  H^{n}. \label{eq:iterative-3}
\end{align}
\item[-] \emph{Step 4}:  find $\pi^{k+1}\in H^1(\Omega)$ such that
 \begin{align}
         & \frac{1}{\Delta t} \int_\Omega \pi^{k+1} w \;d\Omega +\int_\Omega  \Big(\nabla\pi^{k+1} - \frac{2\nabla H^{k+1}}{H^{k+1}} \pi^{k+1}\Big) \nabla w\;d\Omega  =  \frac{1}{\Delta t} \int_\Omega \pi^{n} w \;d\Omega, \label{eq:iterative-4}
\end{align}
for all $w \in H^1(\Omega)$.
\item[-] \emph{Step 5}: Check the stopping criterion, e.g.
\begin{equation}\label{eq:crit}
    \frac{\|A^{k+1}-A^{k}\|_{L^2(\Omega)}}{\|A^{k}\|_{L^2(\Omega)}}< \texttt{tol}_1, \quad \frac{\|\rho^{k+1}-\rho^{k}\|_{L^2(\Omega)}}{\|\rho^{k}\|_{L^2(\Omega)}}< \texttt{tol}_2, 
    \quad \frac{\|\pi^{k+1}-\pi^{k}\|_{L^2(\Omega)}}{\|\pi^{k}\|_{L^2(\Omega)}}< \texttt{tol}_3,
\end{equation}
where $\texttt{tol}_1$, $\texttt{tol}_2$, and $\texttt{tol}_3$ are given stopping tolerances. If not satisfied, repeat steps 1-5. If satisfied, set $A^{n+1} = A^{k+1}$, $\rho^{n+1} = \rho^{k+1}$, $\pi^{n+1} = \pi^{k+1}$, 
$H^{n+1} = H^{k+1}$.
\end{itemize}
At a given time step, algorithm \eqref{eq:iterative-1}-\eqref{eq:crit} solves
four smaller linear problems, i.e., \eqref{eq:iterative-1}-\eqref{eq:iterative-4}, as many times as needed to satisfy
criterion \eqref{eq:crit}. 

Next, we present the space discretization of \eqref{eq:iterative-1}-\eqref{eq:iterative-4}.
We introduce a partition $\mathcal{T}_h$ of $\Omega$ into $N_e$
quadrilaterals $Q_k$ such that
\begin{equation*}
    \overline{\Omega} = \bigcup_{k = 1}^{N_e} \overline{Q}_{k},
\end{equation*}
with 
\begin{equation*}
    h = \max_{Q_k \in \mathcal{T}_h} h_k, \quad h_k = \text{diam}(Q_k), \quad k= 1, \dots, N_e.
\end{equation*}
Let $\mathbb{Q}_N (Q_k)$ be the set of algebraic polynomials, defined on $Q_k$, of degree less than or equal to $N$ in each space variable, and let $V_h (\Omega)$ be the space of global continuous functions on $\overline{\Omega}$ that are
polynomials of degree $N$ on each $Q_k \in \mathcal{T}_h $:
\begin{equation*}
    V_h (\Omega) = \{ v_h \in C^0 (\overline{\Omega})~:~v_h |_{Q_k} \in \mathbb{Q}_N, \forall Q_k \in \mathcal{T}_h  \}.
\end{equation*}
We denote with $\{ \phi_i \}_{i = 1}^{N_Q}$ a basis for $V_h$, where $N_Q$ the total number of degrees of freedom in $\Omega$.

The space-discrete counterpart of \eqref{eq:iterative-1} reads: find $A^{k+1}_h \in V_h (\Omega)$ such that
    \begin{align}
           &\frac{1}{\Delta t}\int_\Omega A^{k+1}_h u\;d\Omega + 
           \int_\Omega \eta\nabla A^{k+1}_h \cdot \nabla u\;d\Omega +\int_\Omega A^{k+1}_h u \;d\Omega \cl
           & \quad \quad -\int_\Omega \rho^{k}_h {e^{- \pi^{k}_h}} A^{k+1}_h u \;d\Omega  = \int_\Omega (A_{0,h} - \eta \Delta A_{0,h}) u \;d\Omega + \frac{1}{\Delta t}\int_\Omega A^{n}_h u\;d\Omega, \label{eq:iterative-1-sd} 
\end{align}
for all $u \in V_h (\Omega)$.
The space-discrete counterpart of 
\eqref{eq:iterative-2} reads: 
find $\rho^{k+1}_h\in V_h (\Omega)$ such that
 \begin{align}
          & \frac{1}{\Delta t} \int_\Omega \rho^{k+1}_h v \;d\Omega +\int_\Omega  \Big(\nabla\rho^{k+1}_h - \frac{2\nabla A^{k+1}_h}{A^{k+1}_h} \rho^{k+1}_h\Big) \nabla v\;d\Omega  \cl
         & \quad \quad +\int_\Omega \rho^{k+1}_h A^{k+1}_h {e^{-\pi^{k}_h}} v\;d\Omega = \int_\Omega \frac{\Gamma \theta}{\omega^2} (1 - \Sigma) {e^{-\pi^{k}_h}} v\;d\Omega + \frac{1}{\Delta t} \int_\Omega \rho^{n}_h v \;d\Omega, \label{eq:iterative-2-sd}
\end{align}
for all $v \in V_h (\Omega)$. In a space-discrete setting, \eqref{eq:iterative-3} simply becomes: set
\begin{align}
 \left( \frac{1}{\Delta t} + \frac{1}{\tau} \right)  H^{k+1}_h    =  \frac{1}{\tau} \Big(\rho^{k+1}_h A^{k+1}_h e^{- \pi^{k}_h} \Big)  +  \frac{1}{\Delta t}  H^{n}_h. \label{eq:iterative-3-sd}
\end{align}
The space discretization of~\eqref{eq:iterative-4} yields: find $\pi^{k+1}_h\in V_h(\Omega)$ such that
 \begin{align}
         & \frac{1}{\Delta t} \int_\Omega \pi^{k+1}_h w \;d\Omega +\int_\Omega  \Big(\nabla\pi^{k+1}_h - \frac{2\nabla H^{k+1}_h}{H^{k+1}_h} \pi^{k+1}_h\Big) \nabla w\;d\Omega  =  \frac{1}{\Delta t} \int_\Omega \pi^{n}_h w \;d\Omega, \label{eq:iterative-4-sd}
\end{align}
for all $w \in V_h (\Omega)$. 

Obviously, stopping criterion \eqref{eq:crit} becomes
\begin{equation}\label{eq:crit-sd}
    \frac{\|A^{k+1}_h - A^{k}_h\|_{L^2(\Omega)}}{\|A^{k}_h\|_{L^2(\Omega)}}< \texttt{tol}_1, \quad  \frac{\|\rho^{k+1}_h -\rho^{k}_h \|_{L^2(\Omega)}}{\|\rho^{k}_h \|_{L^2(\Omega)}}< \texttt{tol}_2,
    \quad \frac{\|\pi^{k+1}_h-\pi^{k}_h\|_{L^2(\Omega)}}{\|\pi^{k}_h\|_{L^2(\Omega)}}< \texttt{tol}_3.
\end{equation}

\section{Numerical experiments}\label{sec:num_res}


We implemented the iterative partitioned scheme described in the previous section in FEniCSx \cite{BarattaEtal2023,BasixJoss,AlnaesEtal2014}, a popular open-source computing platform for solving PDE problems with the finite element method. 
For the agent-based model, we developed  an in-house solver in MATLAB. 
Both the agent-based and the PDE solvers are available on GitHub \cite{OurCode}.  

All the tests presented below involve
computational domain $\Omega = [0,10]\times [0,10]$. For the PDE model, the time step is set to $\Delta t = 1/50$,
$\Omega$ is meshed with a structured 
grid of size $h = 1/10$, and we use first-order Lagrange finite elements on quadrilateral cells
$\mathbb{Q}_1$.
The PDE solver uses the iterative partitioned algorithm, with stopping
tolerances $\texttt{tol}_1 = \texttt{tol}_2 = \texttt{tol}_3 = 10^{-6}.$
In all the tests, we set $A^{st} = 1/50$, $\rho_0 = 0.6$, and $B_0 = \Gamma\theta(1-\Sigma)/\omega^2$. The remaining parameters for each case considered in the following are summarized in Tab.~\ref{tab:case_summary}.

\begin{table}[htb!]
\begin{center}
\begin{tabular}{c c c c c c l}
\toprule
Case & $\eta$ & $\Gamma\theta(1-\Sigma)/\omega^2$ & $\tau$ & $\pi_0$ & Regime & Figures \\
\midrule
1  & $0.7$   & $1.5$ & $5$    & $0.5$ & S & Figs.~\ref{fig:case1_pde}, \ref{fig:case1_ab} \\
2  & $0.3$   & $1.5$ & $5$    & $0.5$ & O & Figs.~\ref{fig:case2_ABM}--\ref{fig:case2_PDE} \\
3  & $0.15$  & $1.5$ & $5$    & $0.5$ & O & Figs.~\ref{fig:eta},  \ref{fig:eta015_snap} \\
4  & $0.075$ & $1.5$ & $5$    & $0.5$ & O & Figs.~\ref{fig:eta}, \ref{fig:eta0075_snap} \\
5  & $0.15$  & $0.5$ & $5$    & $0.5$ & O & Fig.~\ref{fig:bbar} \\
6  & $0.15$  & $2.5$ & $5$    & $0.5$ & O & Figs.~\ref{fig:bbar}, \ref{fig:bbar25_snap} \\
7  & $0.15$  & $1.5$ & $0.5$  & $0.5$ & S & Fig.~\ref{fig:tau} \\
8  & $0.15$  & $1.5$ & $50$   & $0.5$ & O & Figs.~\ref{fig:tau}, \ref{fig:tau50_snap}, \ref{fig:tau50_snapABM} \\
9  & $0.15$  & $1.5$ & $5$    & $0.1$ & O & Figs.~\ref{fig:pi}, \ref{fig:pi01_snap}, \ref{fig:pi01_snapABM} \\
10 & $0.15$  & $1.5$ & $5$    & $1.0$ & O & Figs.~\ref{fig:pi}, \ref{fig:pi1_snap} \\
\bottomrule
\end{tabular}\caption{
Parameter values for each test case considered in this section and figures where the corresponding results are presented.
The regime column indicates whether the case lies in the stable (S) or oscillatory (O) part of the stability diagrams in Fig.~\ref{fig:phase_diagrams}.}
\label{tab:case_summary}
\end{center}
\end{table}

\noindent
{\bf Spatially homogeneous equilibrium}.
We start by validating both PDE and agent-based solver using the homogeneous equilibrium solutions
in \eqref{eq:equil_sol} and \eqref{eq:hom_eq_sol}. For ease of comparison, we will report the non-dimensional solution given by both solvers. Thus, 
all the quantities mentioned in the following are non-dimensional. 
The time interval of interest is $[0, 150]$.
We consider case 1. 
With the given parameters, 
the equilibrium solution is $\bar{A} = 0.9298$, 
$\bar{\rho} = 1.6133$, and $\bar{H} = 0.9098$. As shown in the left panel of Fig.~\ref{fig:phase_diagrams}, 
case 1 belongs to the stable part of the plane.
Figs.~\ref{fig:case1_pde} and 
\ref{fig:case1_ab} show the level of attractiveness, time-lagged crime map, densities of criminals and policemen  given by the both solvers
at the end of the simulation.
We chose different color bars for Figs.~\ref{fig:case1_pde} and 
\ref{fig:case1_ab} for visualization purposes.
In both figures, we see that the system approaches a steady-state
solution that is uniform in space (i.e., no hotspots form), approximating well the expected equilibrium solution.


\begin{figure}[htb!]
     \centering
         \begin{overpic}[percent,width=0.24\textwidth]{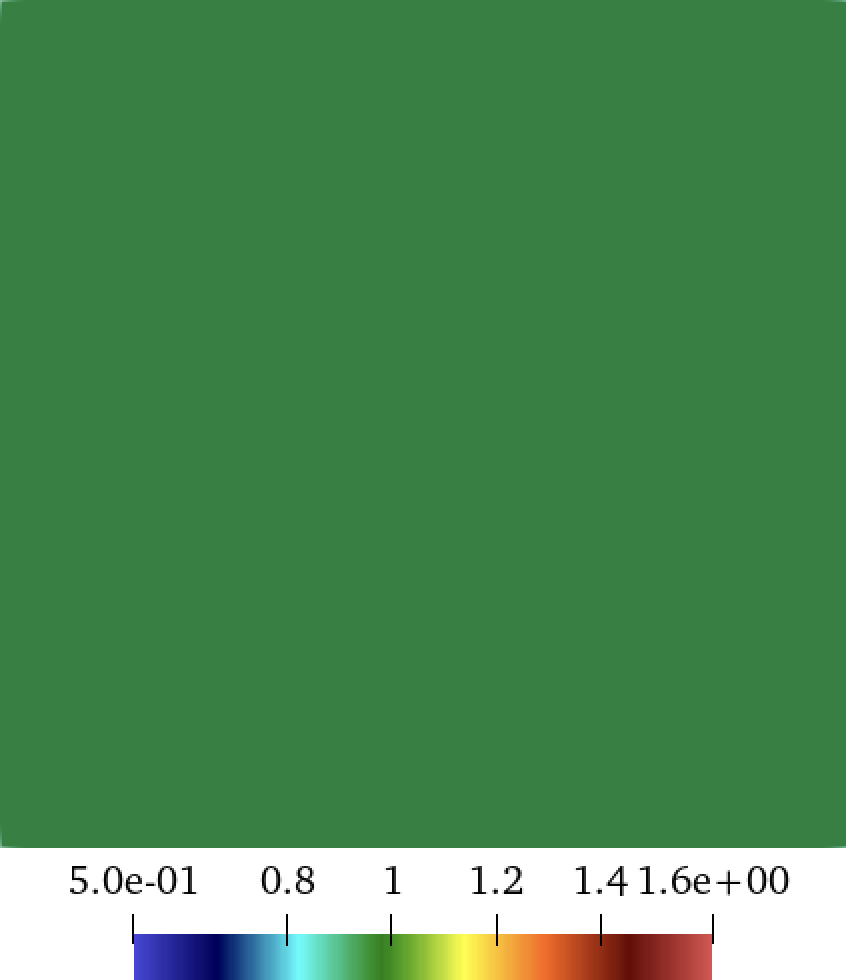}
     \put(43,102){\footnotesize{$A$}}
    \end{overpic} 
     \begin{overpic}[percent,width=0.24\textwidth, grid=false]{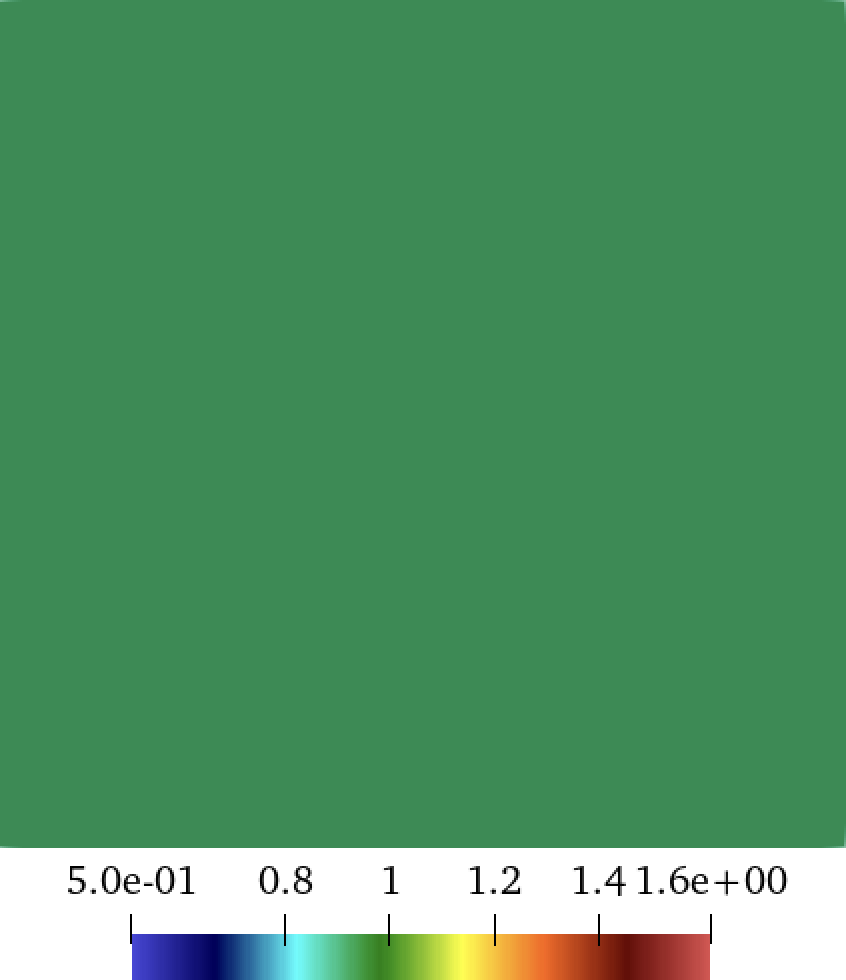}
     \put(43,102){\footnotesize{$H$}}
    \end{overpic}
    \begin{overpic}[percent,width=0.24\textwidth]{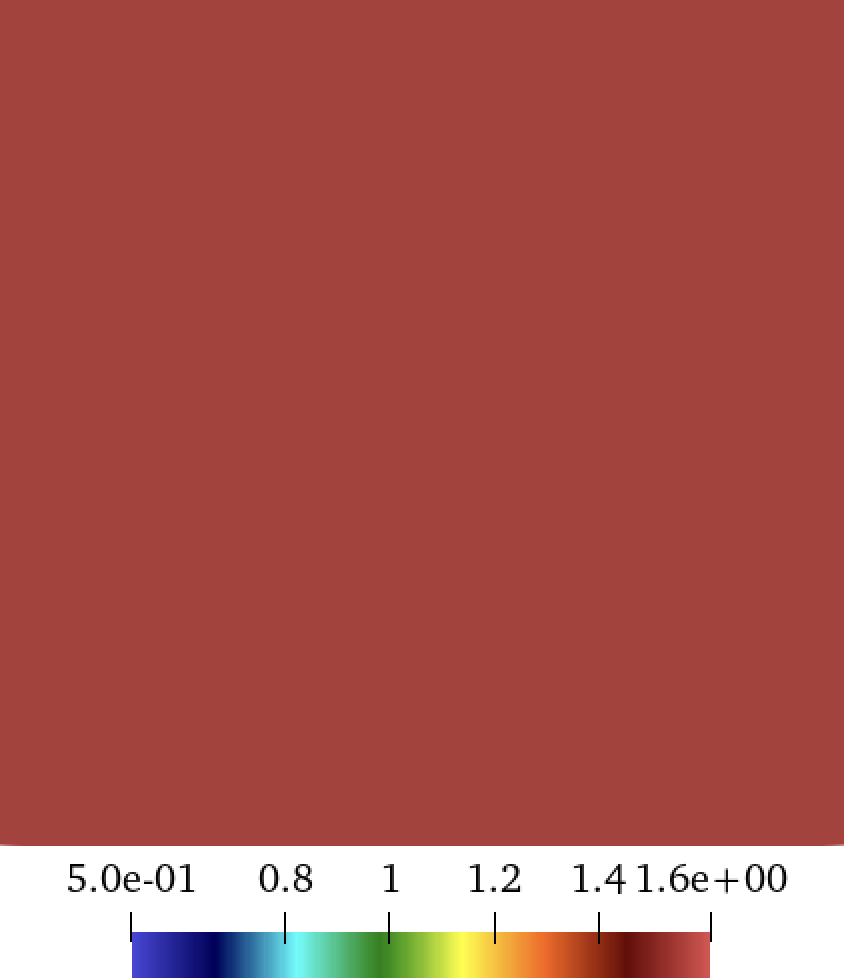}
    \put(43,102){\footnotesize{$\rho$}}
    \end{overpic} 
        \begin{overpic}[percent,width=0.24\textwidth]{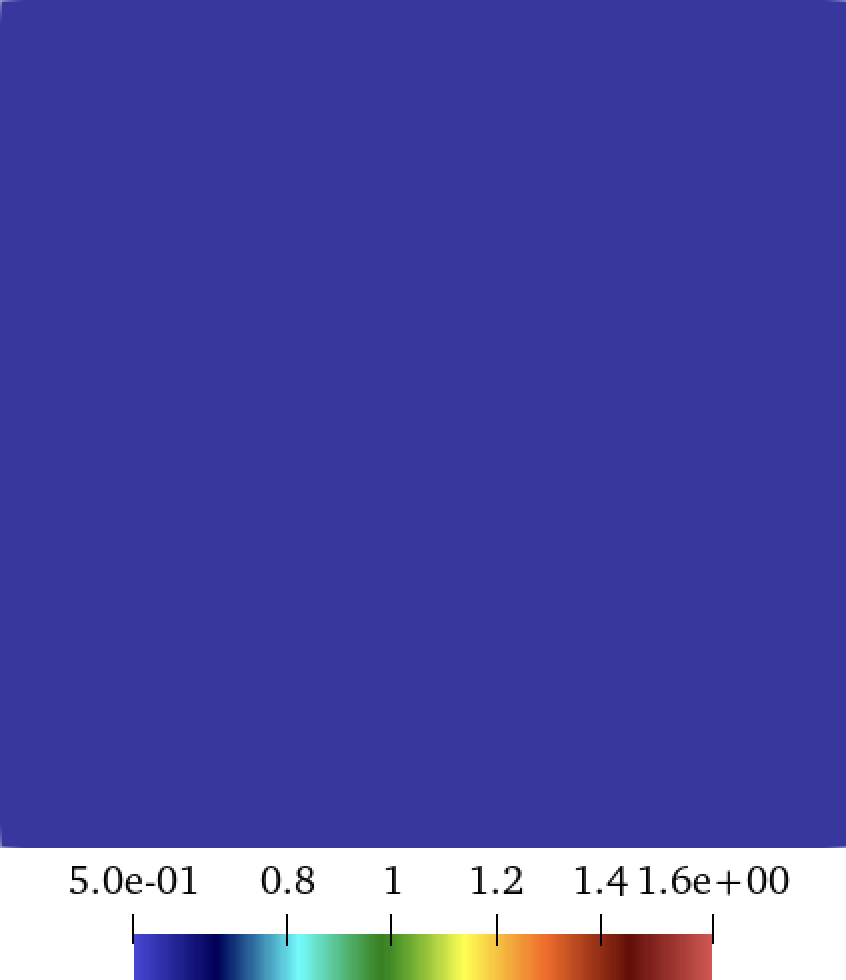}
        \put(43,102){\footnotesize{$\pi$}}
    \end{overpic} 
    
    \caption{Case 1 - PDE model: level of attractiveness $A$, delayed crime density $H$, density of criminals $\rho$, and density of policemen $\pi$
    at the end of the simulation ($T = 150$).}
    \label{fig:case1_pde}
\end{figure}

\begin{figure}[htb!]
     \centering
         \begin{overpic}[percent,width=0.24\textwidth]{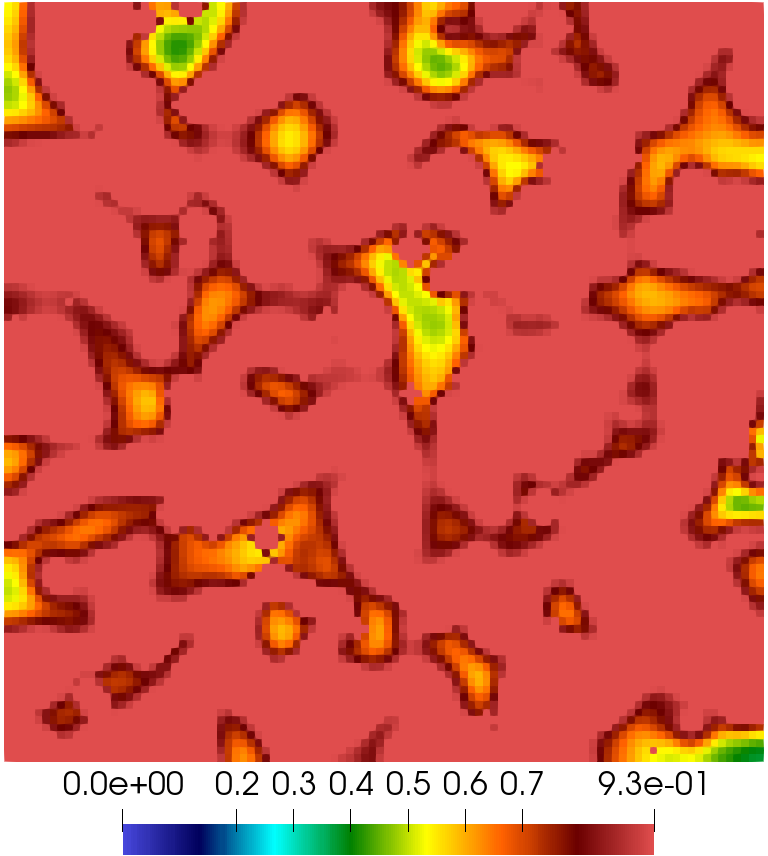}
     \put(43,102){\footnotesize{$A$}}
    \end{overpic} 
     \begin{overpic}[percent,width=0.24\textwidth, grid=false]{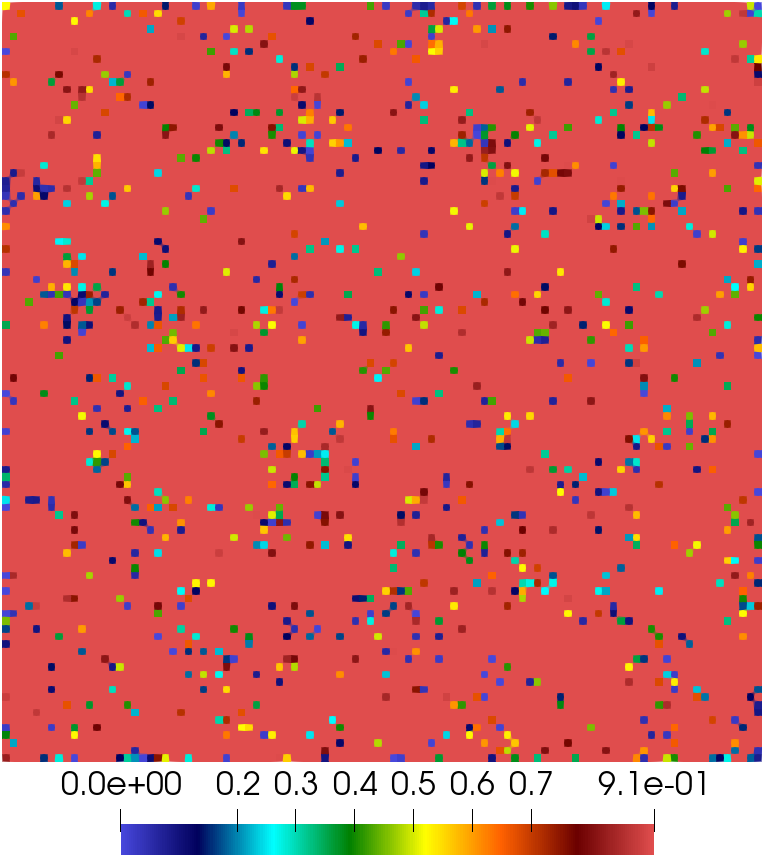}
     \put(43,102){\footnotesize{$H$}}
    \end{overpic}
    \begin{overpic}[percent,width=0.24\textwidth]{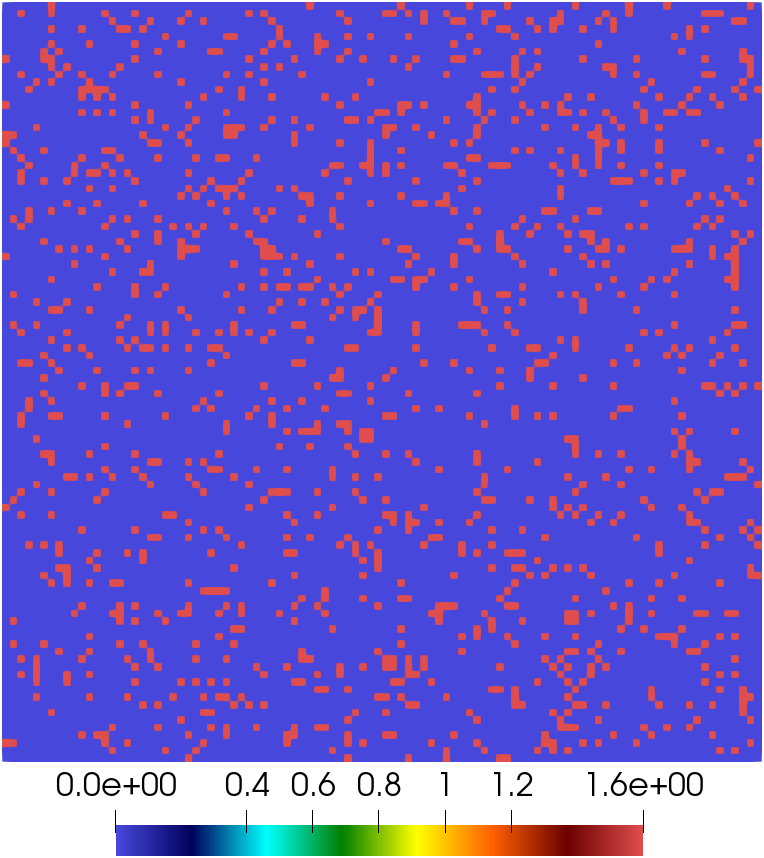}
    \put(43,102){\footnotesize{$\rho$}}
    \end{overpic} 
        \begin{overpic}[percent,width=0.24\textwidth]{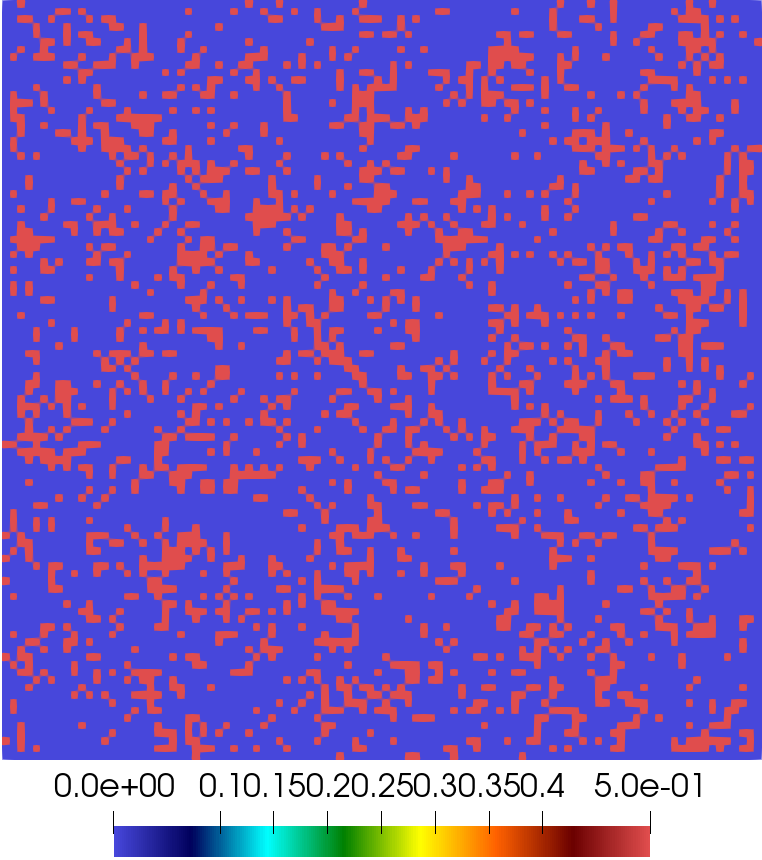}
        \put(43,102){\footnotesize{$\pi$}}
    \end{overpic} 
     
    
    
        
    
    \caption{Case 1 - agent-based model: level of attractiveness $A$, delayed crime density $H$, density of criminals $\rho$, and density of policemen $\pi$
    at the end of the simulation ($T = 150$). 
    }
    \label{fig:case1_ab}

\end{figure}

\noindent {\bf Sample solutions}.
Next, we investigate the effect of the different parameters on the solution.
Recall that model \eqref{eq:continous_1}-\eqref{eq:continous_4}
features 3 parameters: $\eta$, $\Gamma \theta(1-\Sigma)/\omega^2$, and $\tau$. 
We start by lowering the value of $\eta$ in case 1 while leaving all other parameters unchanged, as in case 2, which
belongs to the oscillatory part of the plane as shown in the left panel of Fig.~\ref{fig:phase_diagrams}.


To investigate this oscillatory nature,
we consider spatially averaged quantities
\[
\langle f\rangle(t) := \frac{1}{|\Omega|}\int_\Omega f(x,t)\,dx,
\]
where $f$ can be a system variable or derived quantity like, e.g., $S$ \eqref{eq:S}.
Figs.~\ref{fig:case2_ABM} and \ref{fig:case2_PDE} show the evolution of 
$\langle A \rangle$, $\langle\rho \rangle$, $\langle H \rangle$, and $\langle \pi \rangle$
for $t \in [0, 500]$ given by the agent-based solver and PDE solver, respectively.
Since the total number of policemen is conserved, we see that 
$\langle \pi \rangle$ is constant in time in both figures.
The other three variables settle into persistent, bounded oscillations,
with a well-defined amplitude
after an initial transient.
The mean values of the oscillations
given by the agent-based and PDE models match. Additionally, 
we note that, in both Figs.~\ref{fig:case2_ABM} and \ref{fig:case2_PDE}, $\langle\rho \rangle$
 and $\langle A \rangle$
oscillate with comparable amplitude, 
while  
$\langle H \rangle$ follows with smaller amplitude.

\begin{figure}[htb!]
    \centering
    \includegraphics[width=.8\linewidth]{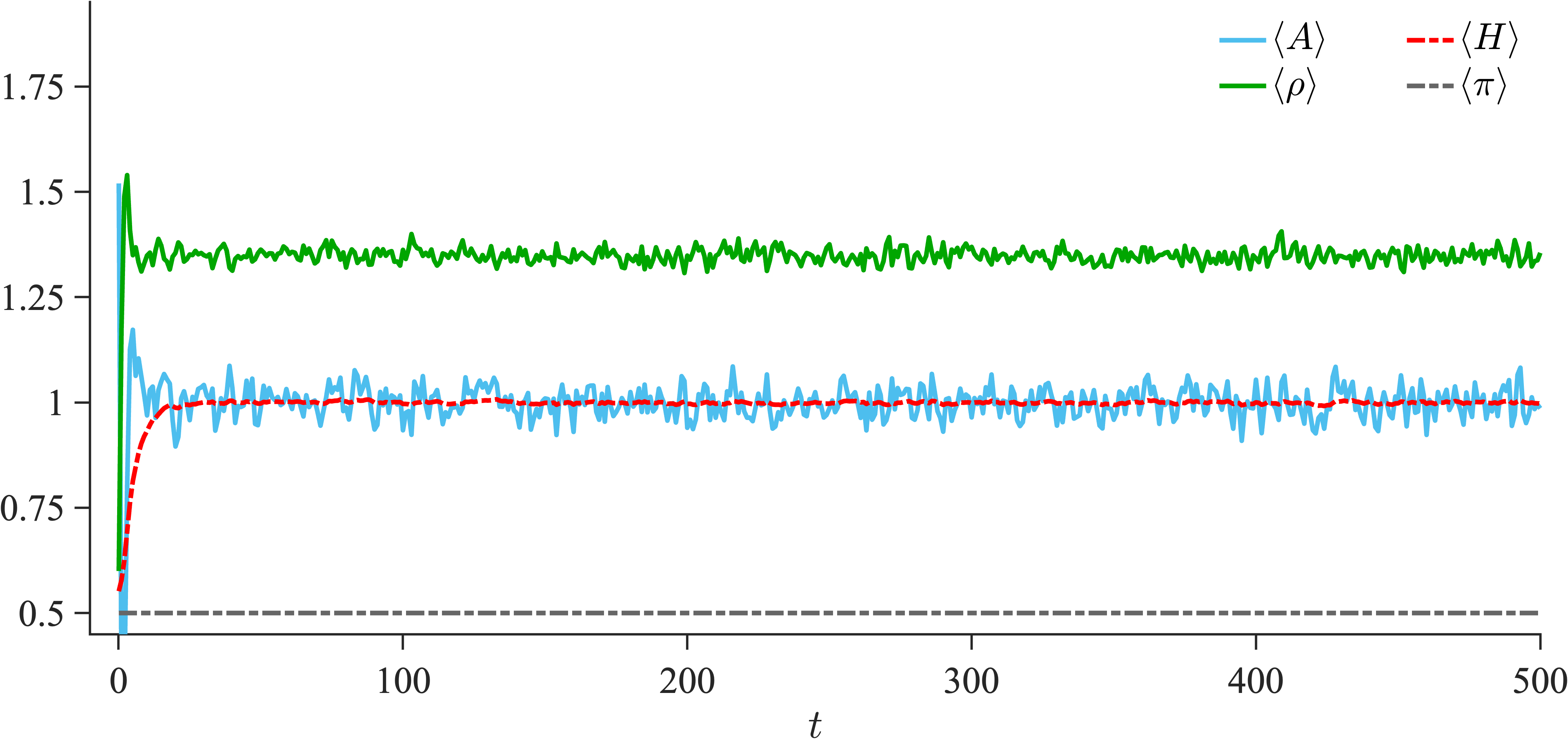}
    \caption{Case 2 - agent-based model: time evolution of spatially averaged quantities $\langle A \rangle$, $\langle \rho \rangle$, $\langle H \rangle$, and $\langle \pi \rangle$.}
    \label{fig:case2_ABM}
\end{figure}

\begin{figure}[htb!]
    \centering
    \includegraphics[width=.8\linewidth]{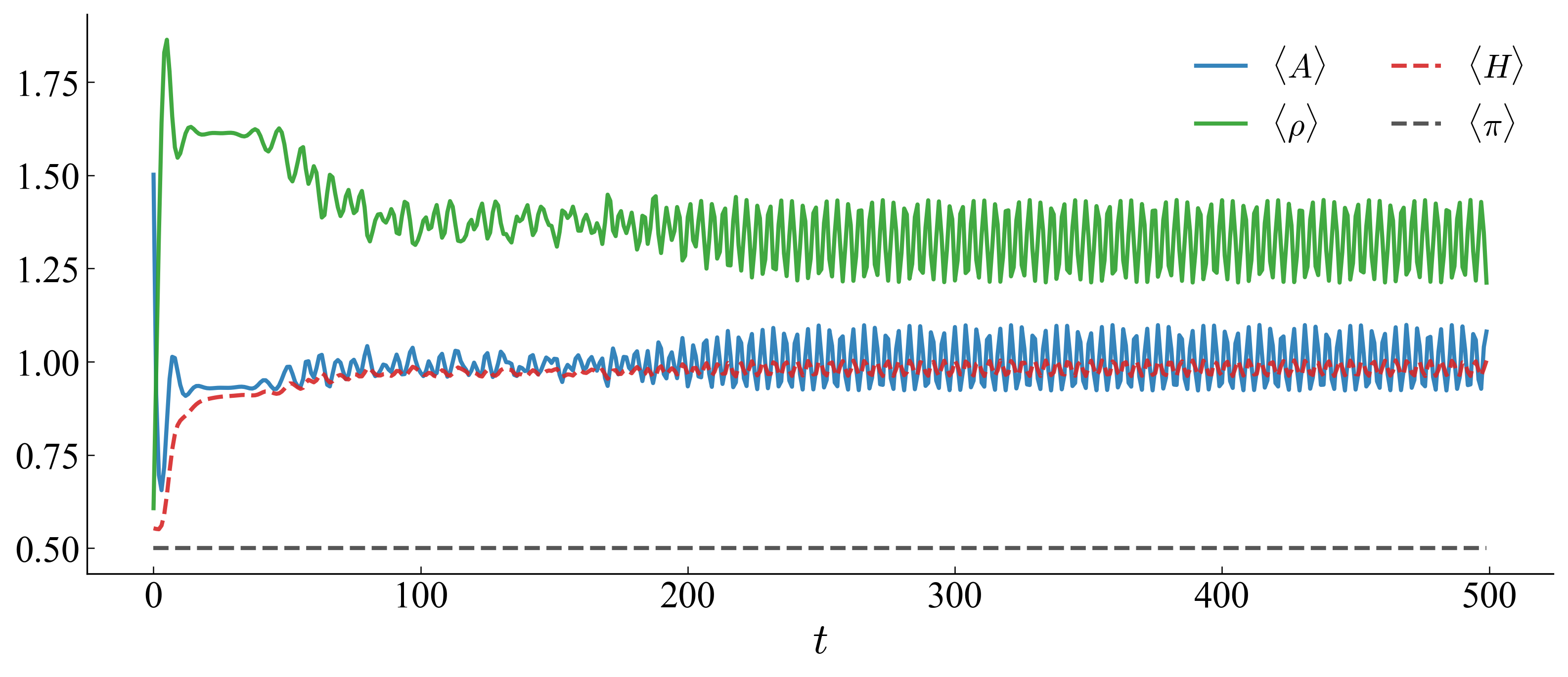}
    \caption{Case 2 - PDE model: time evolution of spatially averaged quantities $\langle A \rangle$, $\langle \rho \rangle$, $\langle H \rangle$, and $\langle \pi \rangle$.}
    \label{fig:case2_PDE}
\end{figure}

Since the oscillations in 
Fig.~\ref{fig:case2_PDE} have a well-defined
period and are synchronized, indicating that the long-time dynamics are governed by a global feedback mechanism, we
further analyze them. 
Fig.~\ref{fig:case2_fft} displays the power spectra of $\langle A \rangle$, $\langle\rho \rangle$, $\langle H \rangle$, and $\langle \pi \rangle$ given by the PDE model
for $t \in [500, 1200]$, i.e., an interval of time when the oscillatory signals show periodic oscillations around a stable mean value. There is a dominant sharp peak around frequency 0.2357
for $\langle A \rangle$, $\langle\rho \rangle$, and $\langle H \rangle$. We also see a smaller peak for a higher frequency, indicating that the signals are not perfectly monochromatic, as
one could already tell from Fig.~\ref{fig:case2_PDE}. The dominant frequency indicates a period of about 4.24 time unites.
In Fig.~\ref{fig:case2_snap}, we report the plots of $S$ computed by the PDE solver
at selected times during over one 
fundamental cycle ($T_f\approx 17$ time units). The cycle visits four symmetry-related configurations: 
symmetric hotspots ($t = 739$), horizontal splitting ($t = 741$), 
near-symmetric ($t = 743$), and vertical splitting ($t = 746$). Finally, at $t=756$ we recover 
a solution similar to the one at $t=739$.
The spatially averaged signal $\langle S\rangle$ takes similar values 
at each sub-state, producing a dominant FFT peak at $4/T_f \approx 0.236$ 
rather than at the fundamental frequency $1/T_f \approx 0.059$.
Note that no part of the domain is safe
from crime. 

\begin{figure}[htb!]
    \centering
    \includegraphics[width=.8\linewidth]{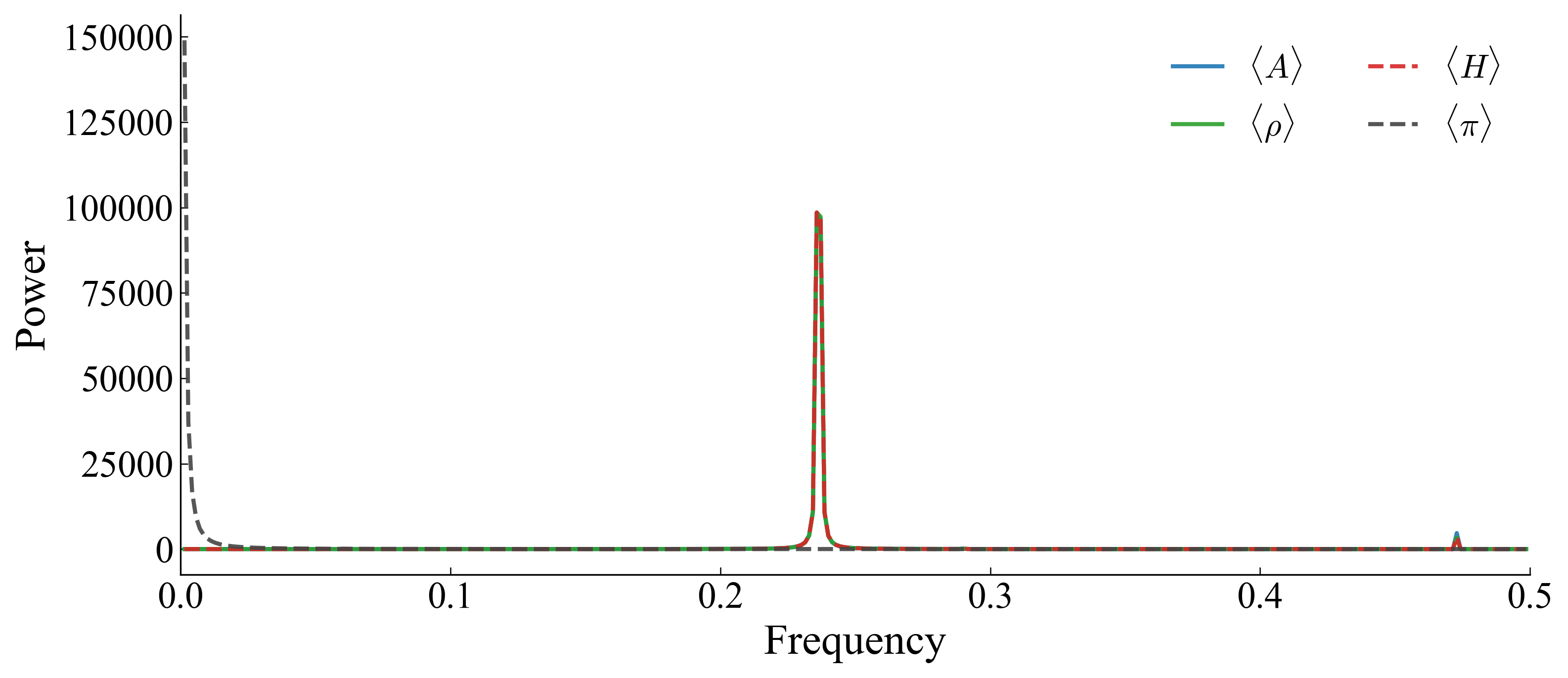}
    \caption{Case 2 - PDE model: power spectra of spatially averaged quantities $\langle A \rangle$, $\langle\rho \rangle$, $\langle H \rangle$, and $\langle \pi \rangle$
for $t \in [500, 1200]$.}
    \label{fig:case2_fft}
\end{figure}

\begin{figure}[htb!]
     \centering
         \begin{overpic}[percent,width=0.16\textwidth]{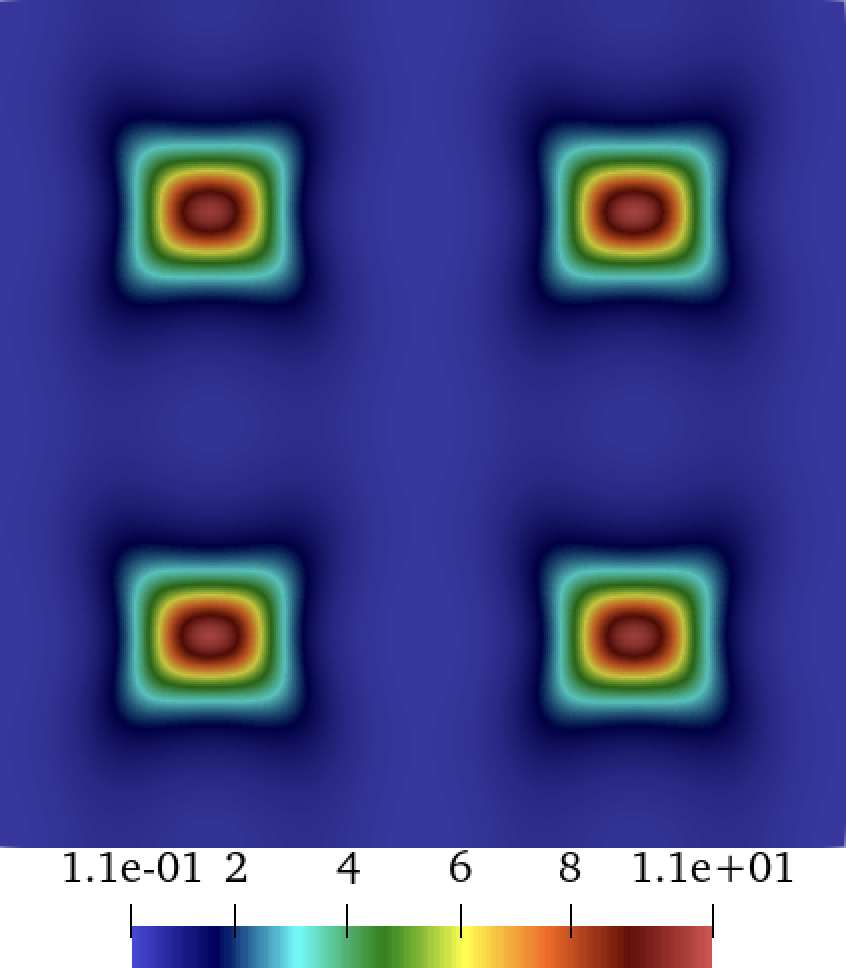}
     \put(23,102){\footnotesize{$t=739$}}
    \end{overpic} 
     \begin{overpic}[percent,width=0.16\textwidth, grid=false]{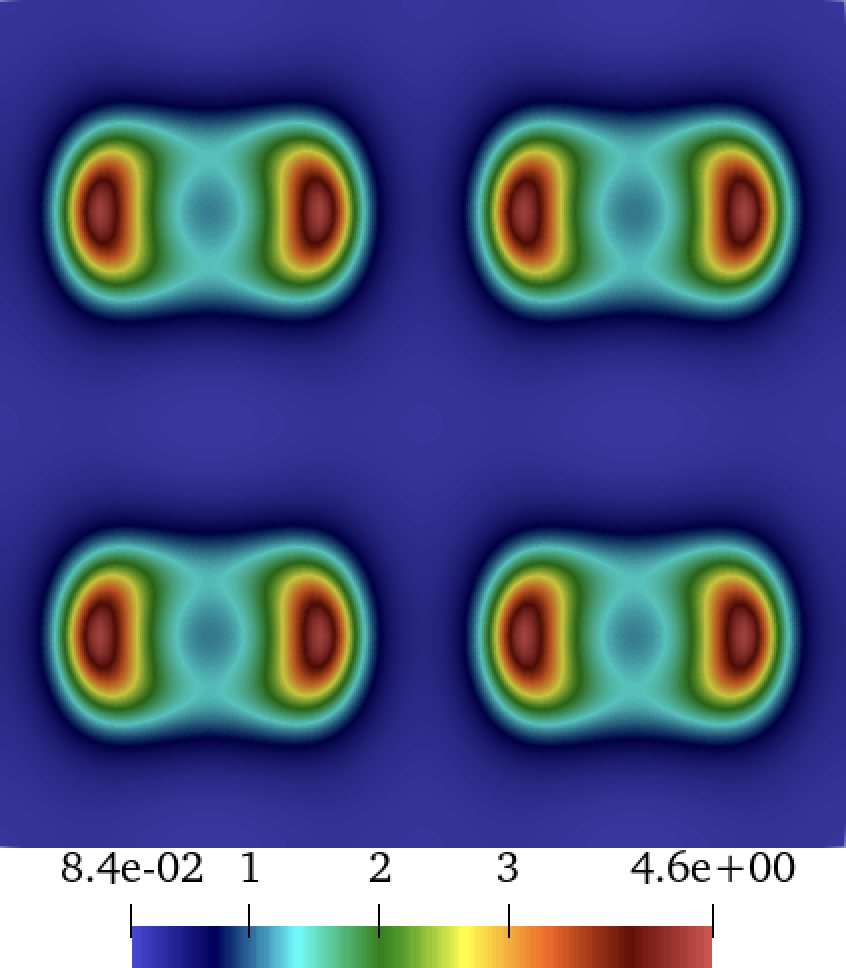}
     \put(23,102){\footnotesize{$t=741$}}
    \end{overpic}
    \begin{overpic}[percent,width=0.16\textwidth]{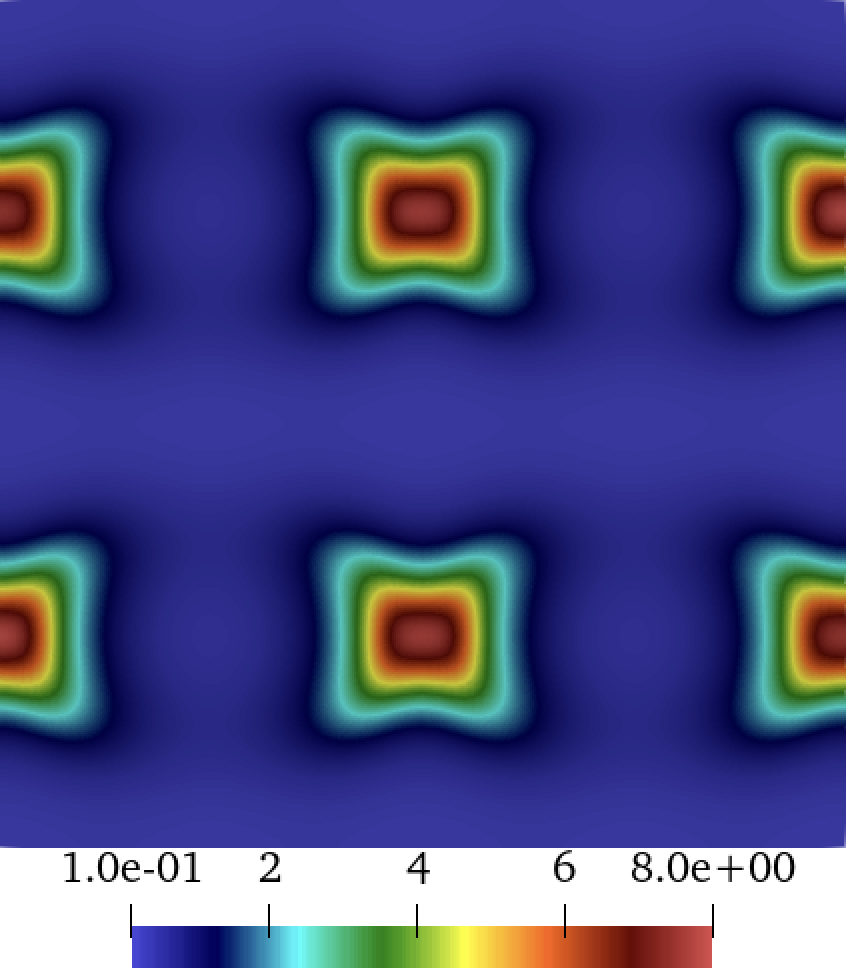}
  \put(23,102){\footnotesize{$t=743$}}
    \end{overpic} 
        \begin{overpic}[percent,width=0.16\textwidth]{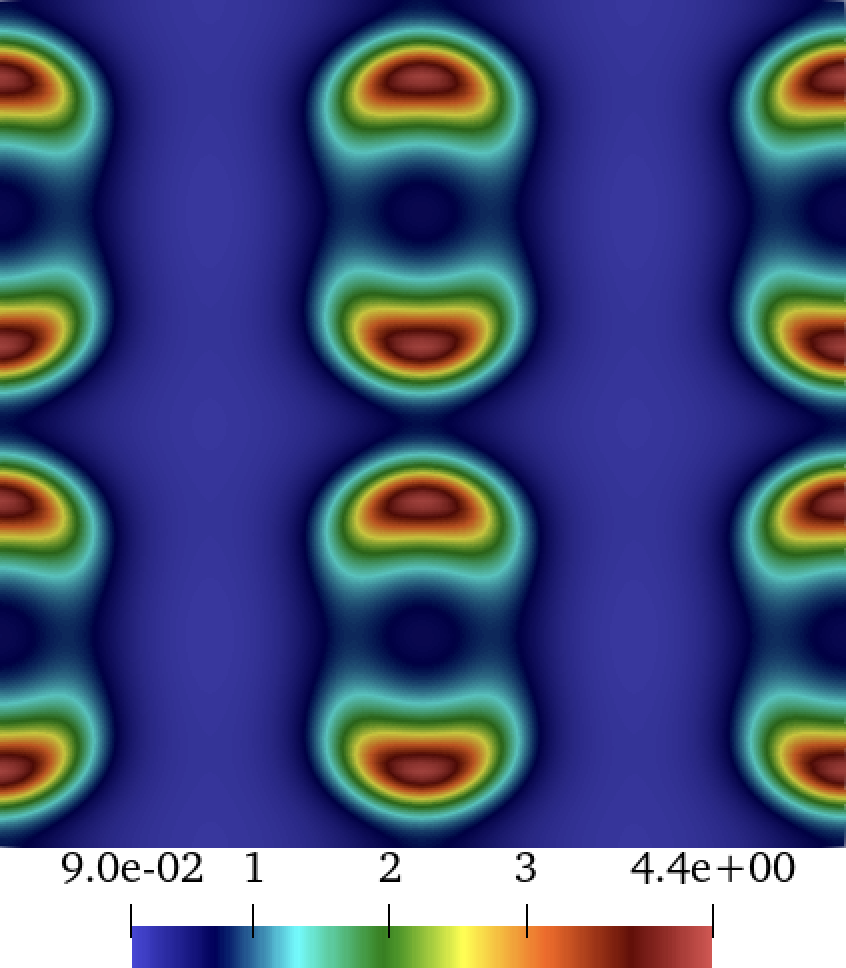}
       \put(23,102){\footnotesize{$t=746$}}
    \end{overpic} 
      \begin{overpic}[percent,width=0.16\textwidth]{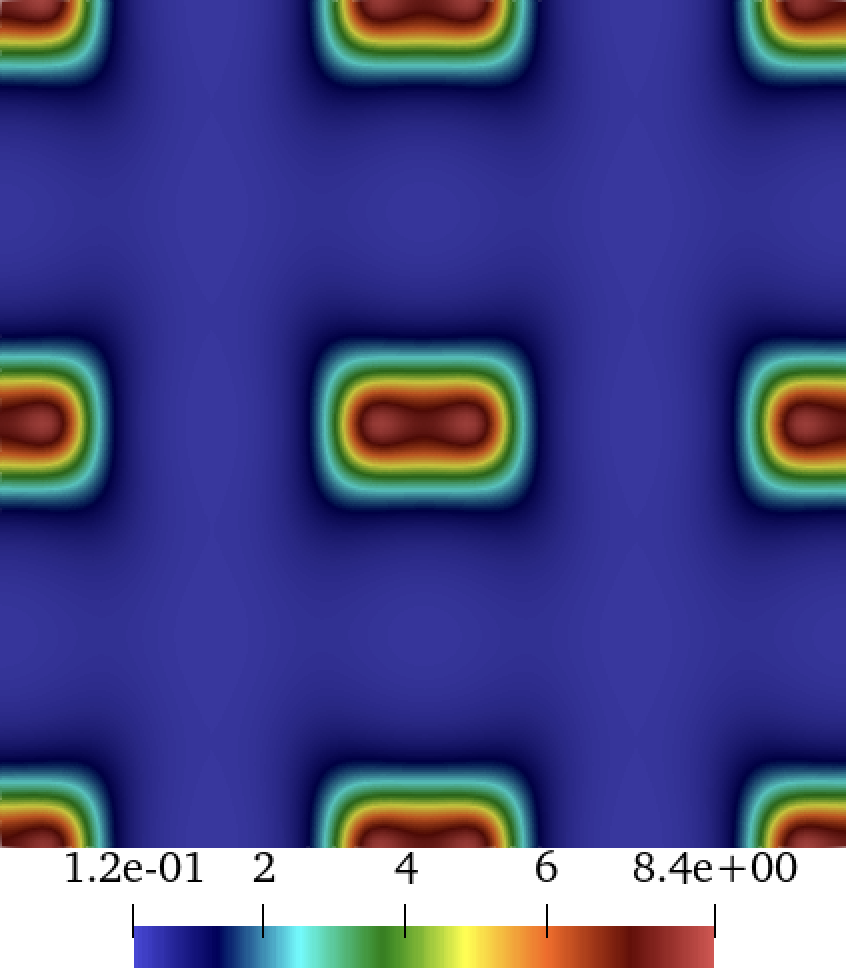}
       \put(23,102){\footnotesize{$t=748$}}
    \end{overpic} 
  \begin{overpic}[percent,width=0.16\textwidth]{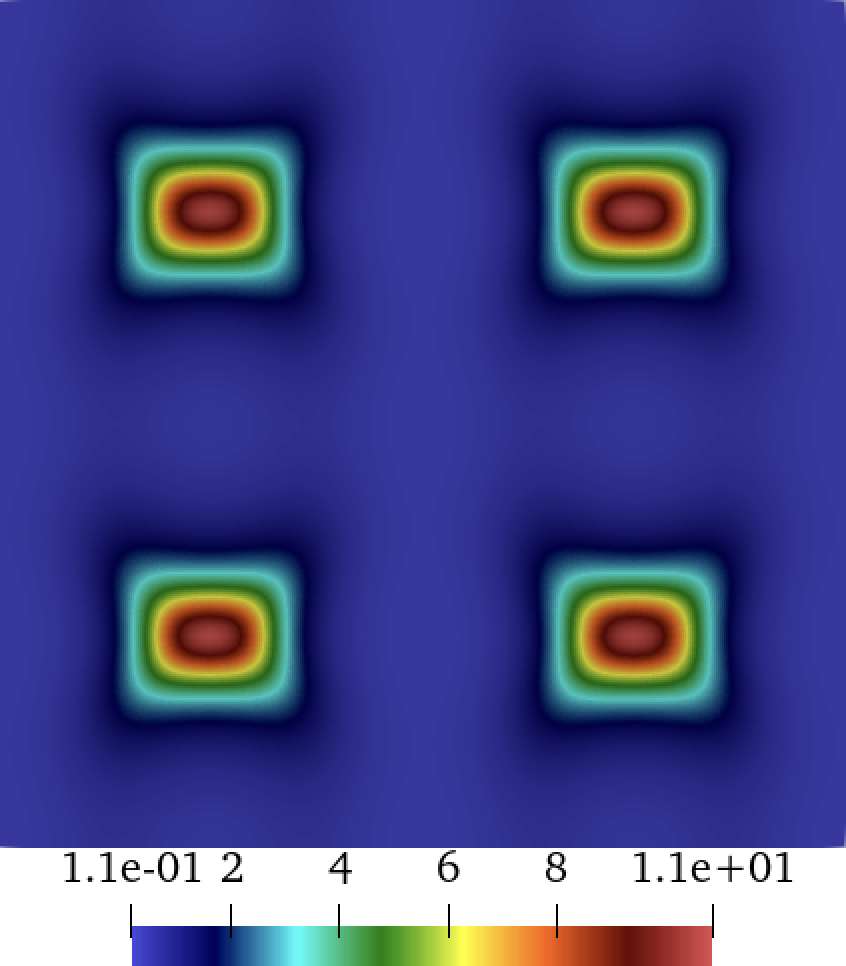}
      \put(23,102){\footnotesize{$t=756$}}
    \end{overpic} 

    \caption{Case 2 - PDE model: expected number of crimes per unit area $S$ from $t = 739$ to $t = 756$.}
    \label{fig:case2_snap}
\end{figure}



If one was to add noise to the initial conditions \eqref{eq:ic}-\eqref{eq:ic4} in the fashion of \cite{short2008statistical,Hao2026}, the effect of the noise would not alter the shape or position of the hotspots in Fig.~\ref{fig:case2_snap}.

Fig.~\ref{fig:portrait} displays phase portraits of  $\langle A\rangle$ versus $\langle\rho\rangle$ and $\langle H\rangle$ versus $\langle S \rangle$ given by the PDE solver for $t \in [500, 1200]$. 
In the $(\langle A\rangle,\langle\rho\rangle)$ plane, the trajectory forms a closed, smooth, elliptical loop, indicating a near single-frequency periodic orbit in the given nonlinear system. The ellipse is slightly tilted, which
means that $\langle A\rangle$ and $\langle\rho\rangle$ oscillate out of phase. 
This is due to a phase-lagged interaction between attractiveness and crime,
driven by the nonlinear production term $\rho A e^{-\pi}$. We also note that the thickness of the loop is small, which means that the oscillation amplitude is small. 
In the $(\langle H \rangle,\langle S \rangle)$ plane, the portrait is similar: a smooth, elongated ellipse.
The nearly vertical loop structure highlights the delayed response of $H$
to variations in the crime intensity, consistent with the relaxation dynamics
in \eqref{eq:continous_4}.
Together, the results in Figs.~\ref{fig:case2_PDE}, \ref{fig:case2_fft}, and \ref{fig:portrait}
 demonstrate that the oscillatory behavior is an emergent
property of the coupled system, characterized by a dominant global frequency that
persists despite significant local variability.
This is an important difference from the model with no police \cite{short2008statistical, Hao2026}, which gives solutions tending towards a steady state (i.e., no oscillations) with hotspots as $\eta$ is lowered.

\begin{figure}[htb!]
    \centering
    \includegraphics[width=.8\linewidth]{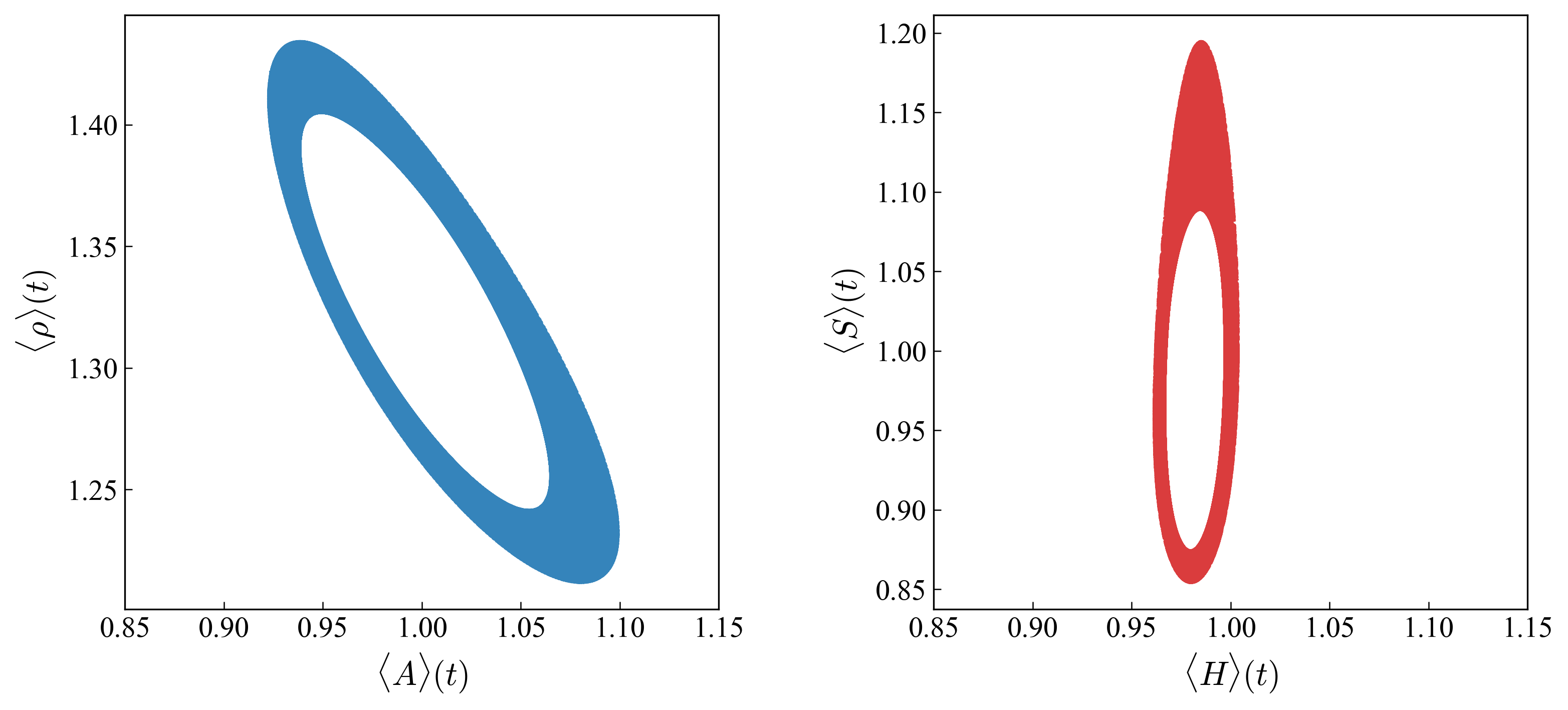}
    \caption{Case 2 - PDE model: phase portraits of  $\langle A\rangle$ versus $\langle\rho\rangle$ (left) and $\langle H\rangle$ versus $\langle S \rangle$ (right) for $t \in [500, 1200]$.}
    \label{fig:portrait}
\end{figure}

In the rest of the paper, we will focus mostly on the results given by the PDE model. 
So, if a figure does not explicitly mention the model,
it means the results therein are obtained with the PDE model.
However, in a couple of cases, we will present results given by the agent-based model too and that will be specified. 

Let us further lower $\eta$ with cases 3 and 4, which lie in the oscillatory part of the plane in Fig.~\ref{fig:phase_diagrams} (left)
and see how it affects the solution. Fig.~\ref{fig:eta}
shows the evolution of $\langle S\rangle$ for cases 2, 3, and 4, and associated power spectra. 
A first observation is that $\eta$ controls
the length of the initial transient: the larger
the value of $\eta$, the longer the initial
transient. Given that after non-dimensionalization $S = \rho A e^{-\pi}$, $\langle S\rangle$ exhibits the same
peaks seen in Fig.~\ref{fig:case2_fft} for $\eta = 0.3$. If $\eta$ is decreased to 
$0.15$, the dominant sharp peak
occurs at a slightly higher frequency. That is the case for the secondary peak too. In addition, we observe the emergence of a third peak, near the dominant peak. 
Upon a further reduction of $\eta$ to $0.075$, the dominant sharp peak moves to a slightly smaller frequency and several more frequencies get excited. 
Figs.~\ref{fig:eta015_snap} and 
\ref{fig:eta0075_snap} report the plots of $S$
at selected times 
for case 3 and 4, respectively.
Just like for the model with no police
\cite{short2008statistical, Hao2026}, the hotspot size decreases with $\eta$.
We observe that the hotspot dynamics in Figs.~\ref{fig:case2_snap} and 
\ref{fig:eta015_snap} are qualitatively similar. In both cases, the system exhibits a regular spatio-temporal evolution characterized by the emergence, deformation, and reorganization of hotspots. The main difference between these two regimes
concern only quantitative features, namely the number of hotspots (increasing as $\eta$ decreases) and their distance (decreasing as $\eta$ decreases). Fig.~\ref{fig:eta0075_snap}
shows a different behavior: hotspots become less regular in shape, frequently split or merge, and display
a more erratic motion, indicating a loss
of spatio-temporal coherence. This is consistent with the power spectra in Fig.~\ref{fig:eta} (right), where case 4 ($\eta = 0.075$) exhibits a broader spectrum with multiple excited frequencies.

\begin{figure}[htb!]
    \centering
    \includegraphics[width=1.0\linewidth]{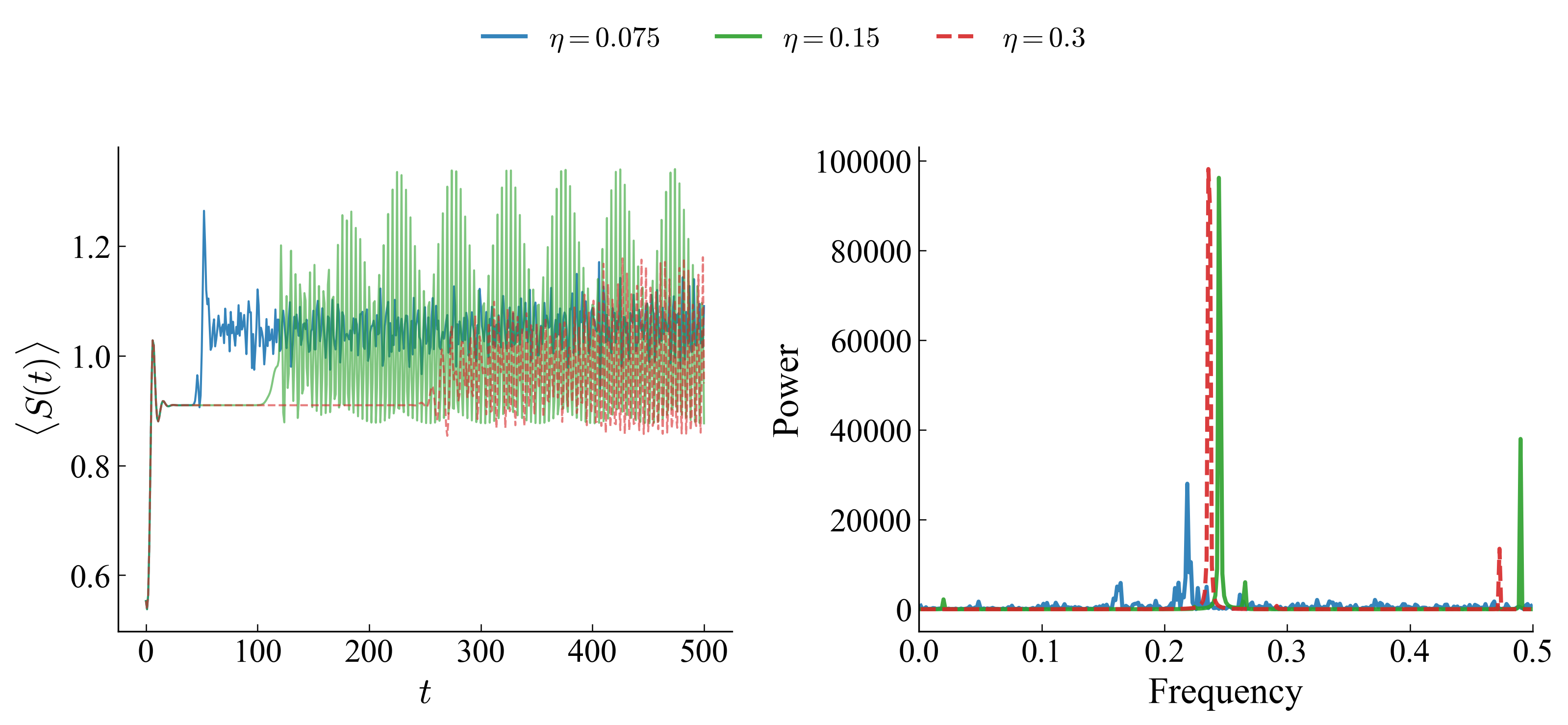}
    \caption{Left: time evolutions of spatially averaged quantities $\langle S \rangle$ for case 2 ($\eta = 0.3$), case 3 ($\eta = 0.15$), and case 4 ($\eta = 0.075$). Right: associated power spectra. }
    \label{fig:eta}
\end{figure}

\begin{figure}[htb!]
     \centering
         \begin{overpic}[percent,width=0.16\textwidth]{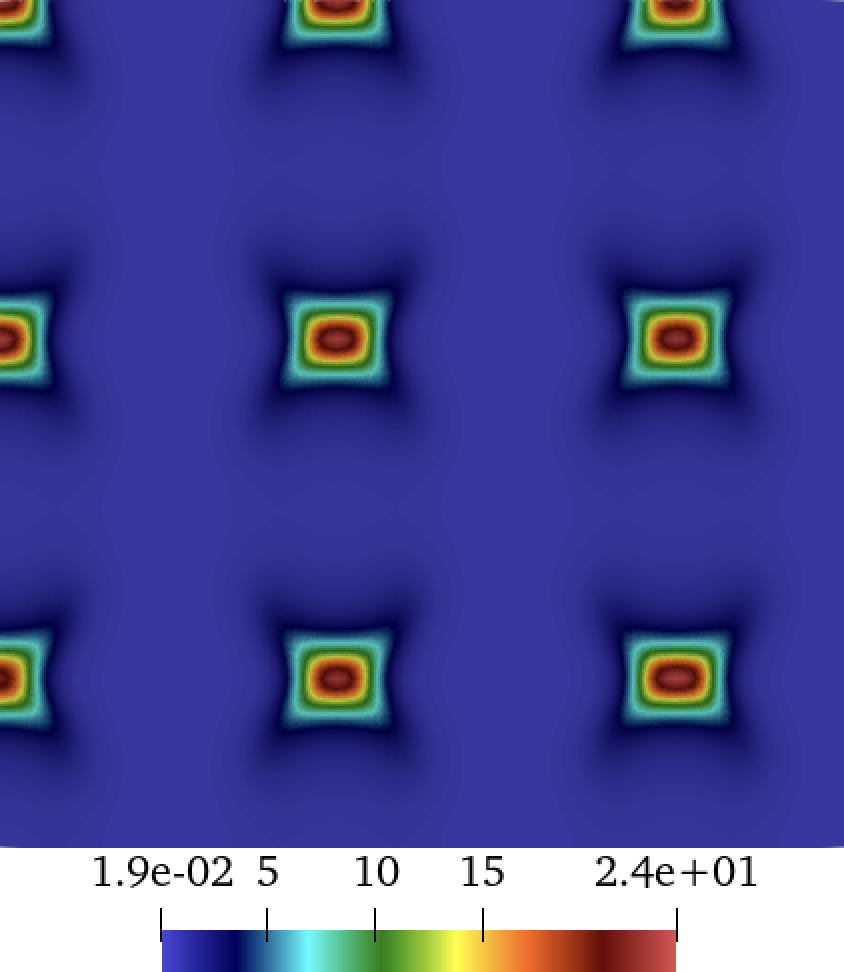}
     \put(23,102){\footnotesize{$t=854$}}
    \end{overpic} 
     \begin{overpic}[percent,width=0.16\textwidth, grid=false]{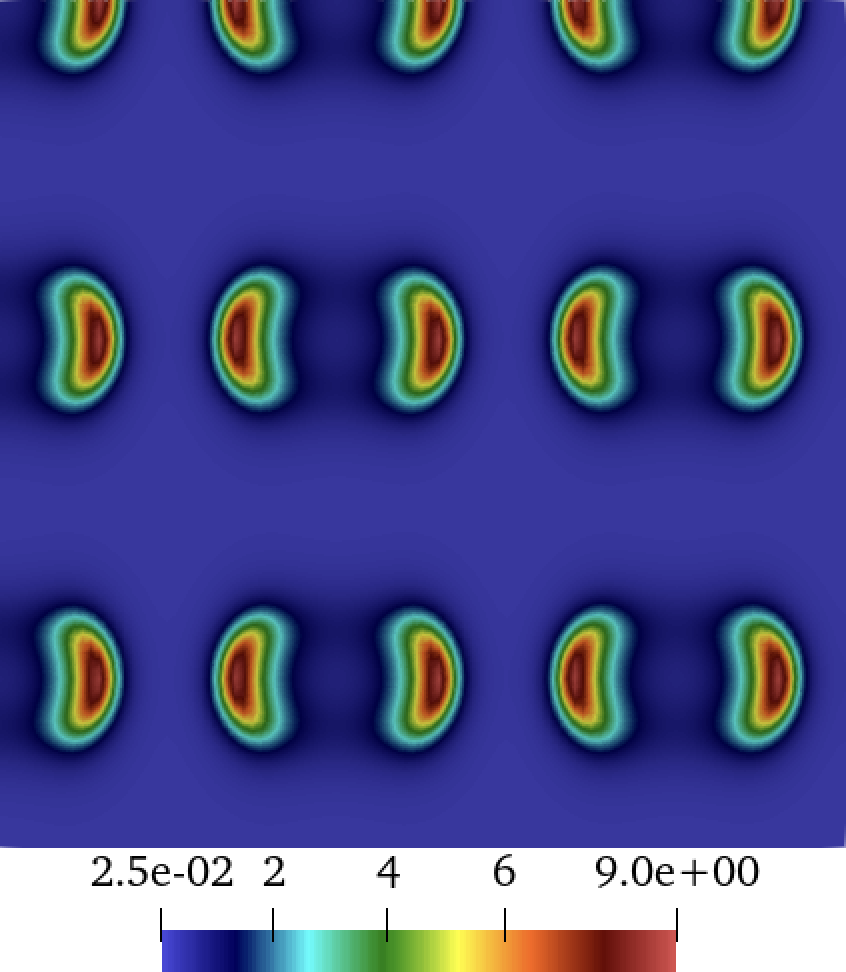}
     \put(23,102){\footnotesize{$t=856$}}
    \end{overpic}
    \begin{overpic}[percent,width=0.16\textwidth]{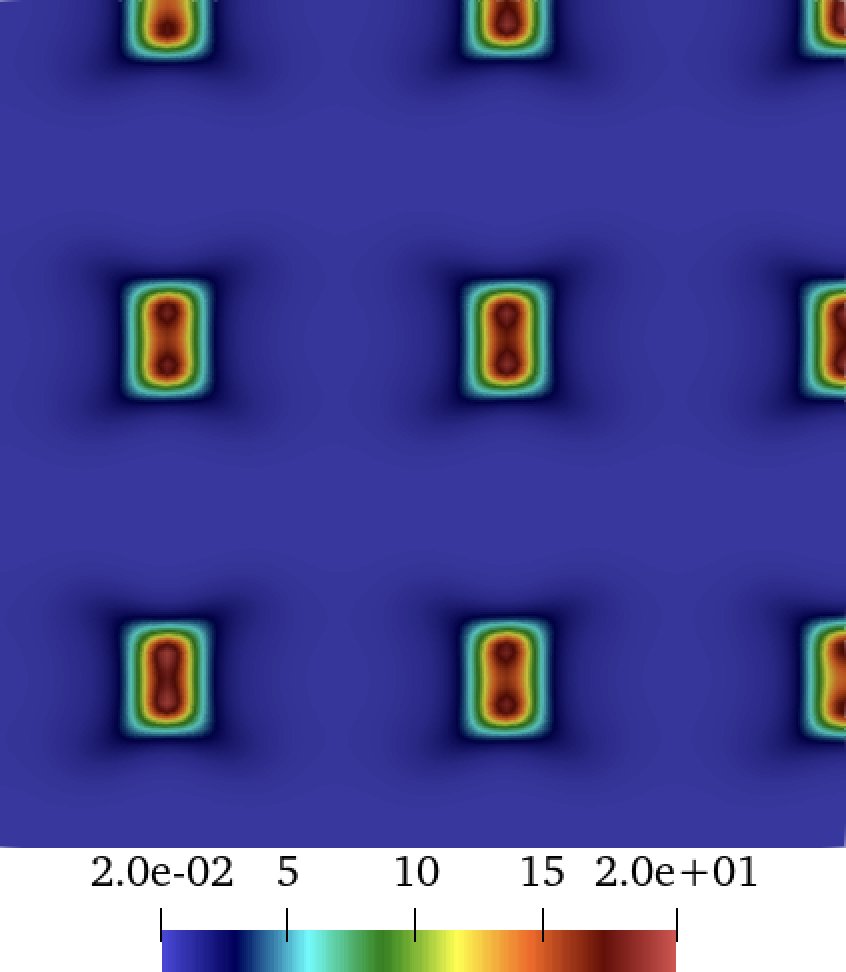}
  \put(23,102){\footnotesize{$t=858$}}
    \end{overpic} 
        \begin{overpic}[percent,width=0.16\textwidth]{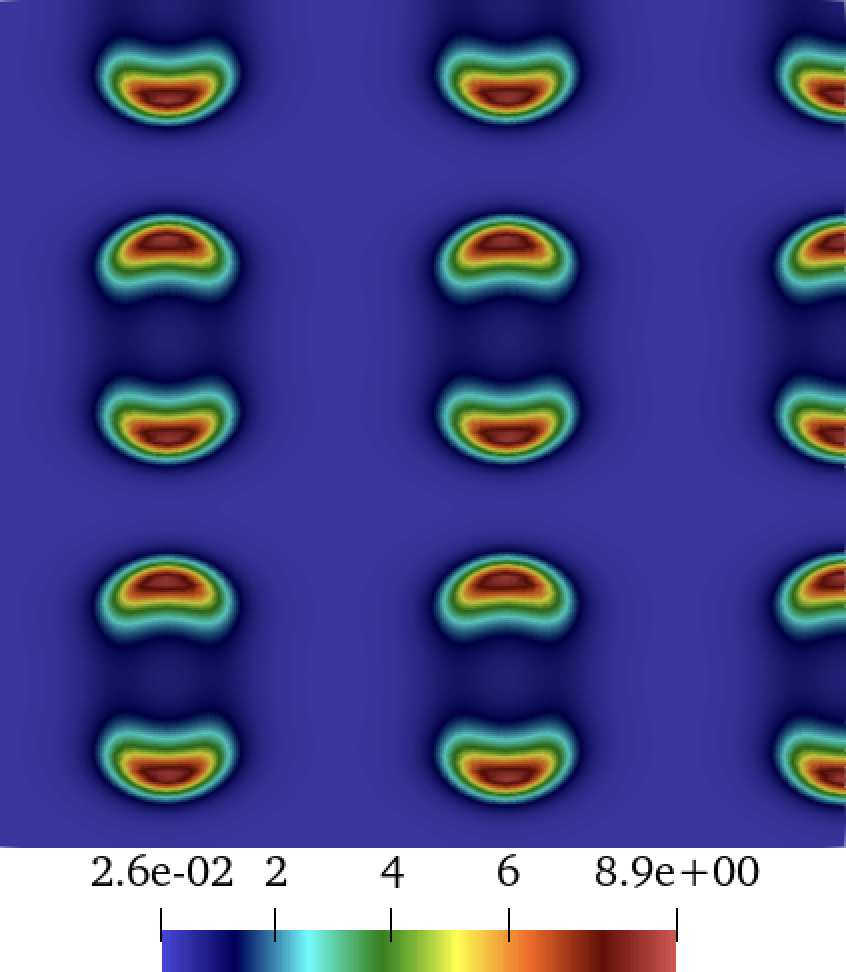}
       \put(23,102){\footnotesize{$t=860$}}
    \end{overpic} 
      \begin{overpic}[percent,width=0.16\textwidth]{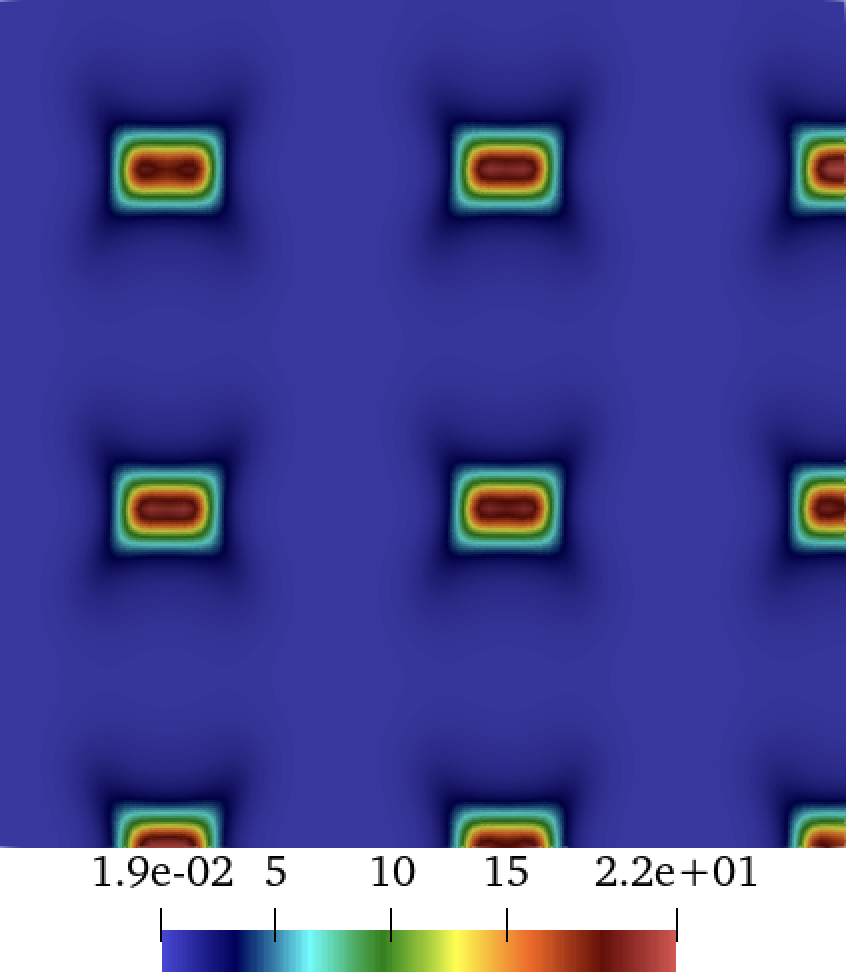}
       \put(23,102){\footnotesize{$t=862$}}
    \end{overpic} 
  \begin{overpic}[percent,width=0.16\textwidth]{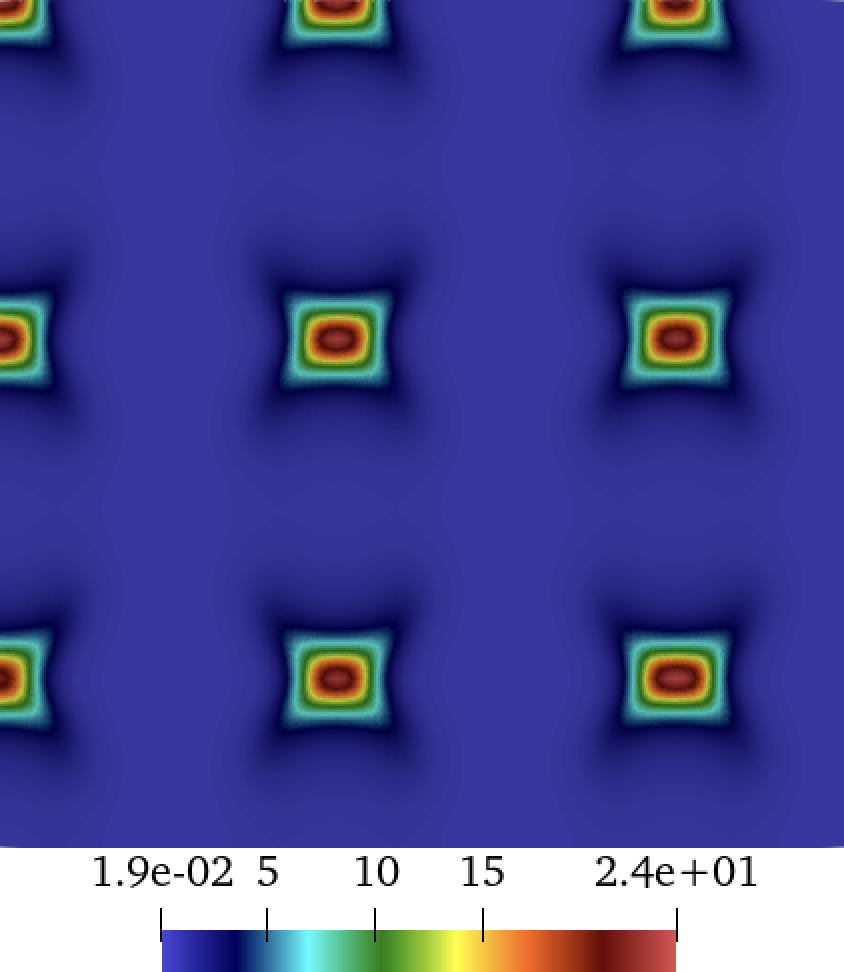}
      \put(23,102){\footnotesize{$t=871$}}
    \end{overpic} 

    \caption{Case 3: expected number of crimes per unit area $S$ from $t = 854$ to $t = 871$. }
    \label{fig:eta015_snap}
\end{figure}

\begin{figure}[htb!]
     \centering
         \begin{overpic}[percent,width=0.16\textwidth]{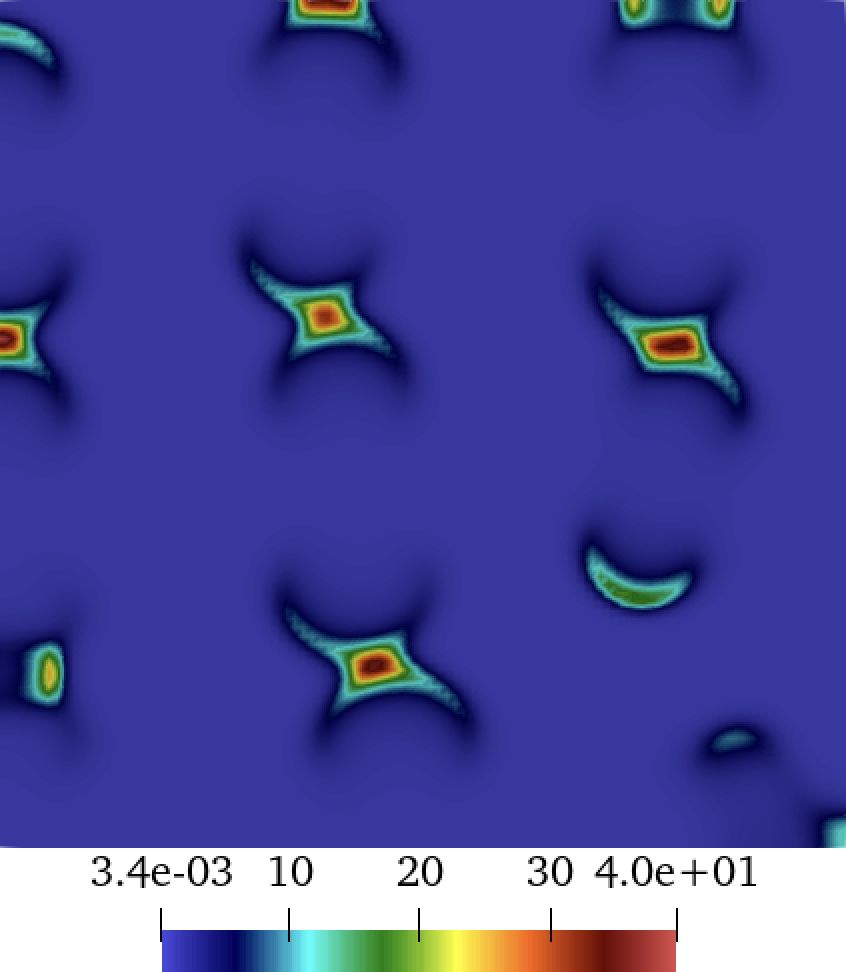}
     \put(23,102){\footnotesize{$t=1360$}}
    \end{overpic} 
     \begin{overpic}[percent,width=0.16\textwidth, grid=false]{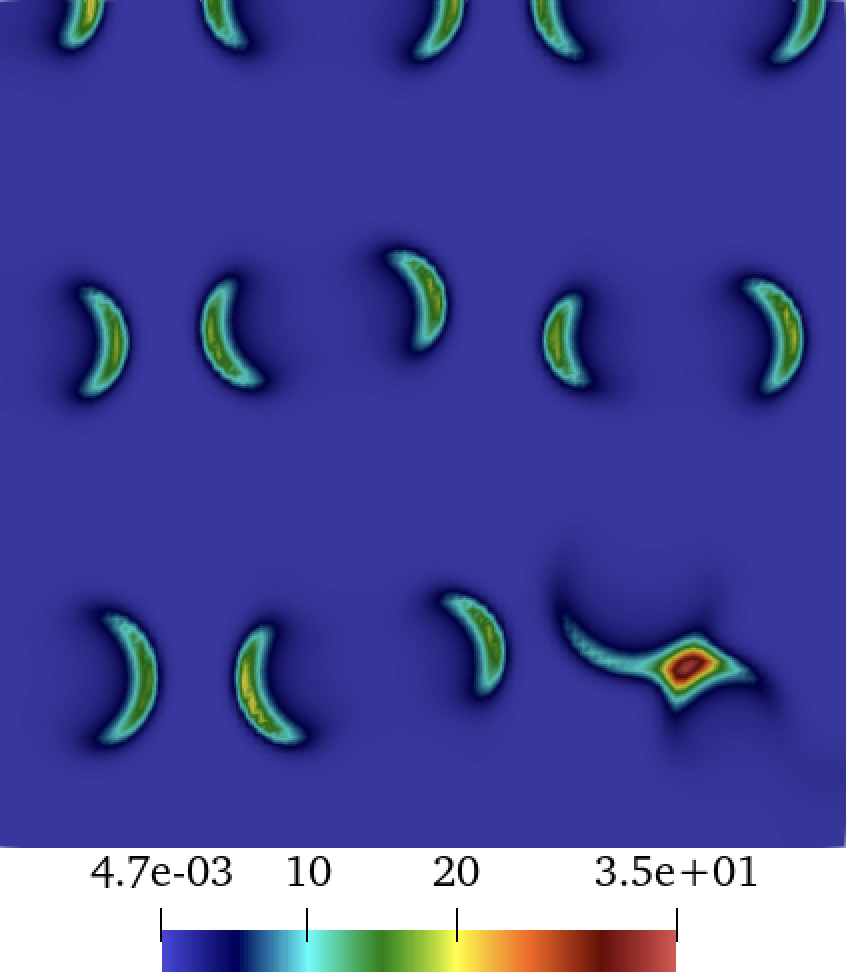}
     \put(23,102){\footnotesize{$t=1363$}}
    \end{overpic}
    \begin{overpic}[percent,width=0.16\textwidth]{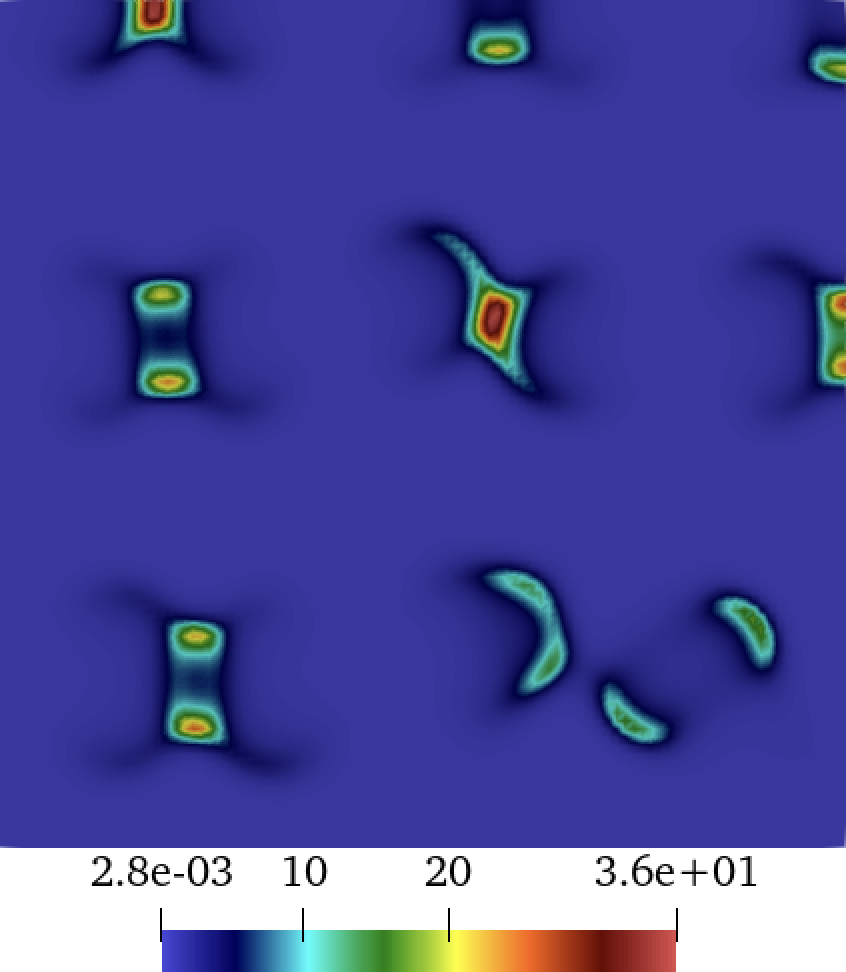}
  \put(23,102){\footnotesize{$t=1365$}}
    \end{overpic} 
        \begin{overpic}[percent,width=0.16\textwidth]{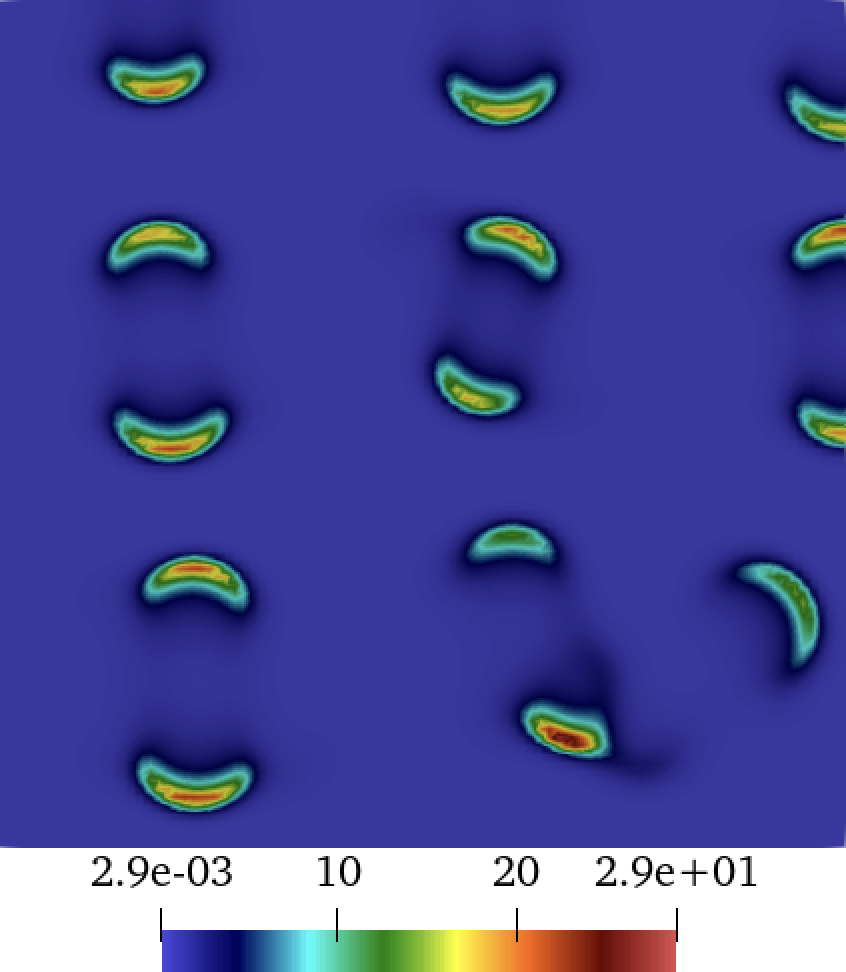}
       \put(23,102){\footnotesize{$t=1367$}}
    \end{overpic} 
      \begin{overpic}[percent,width=0.16\textwidth]{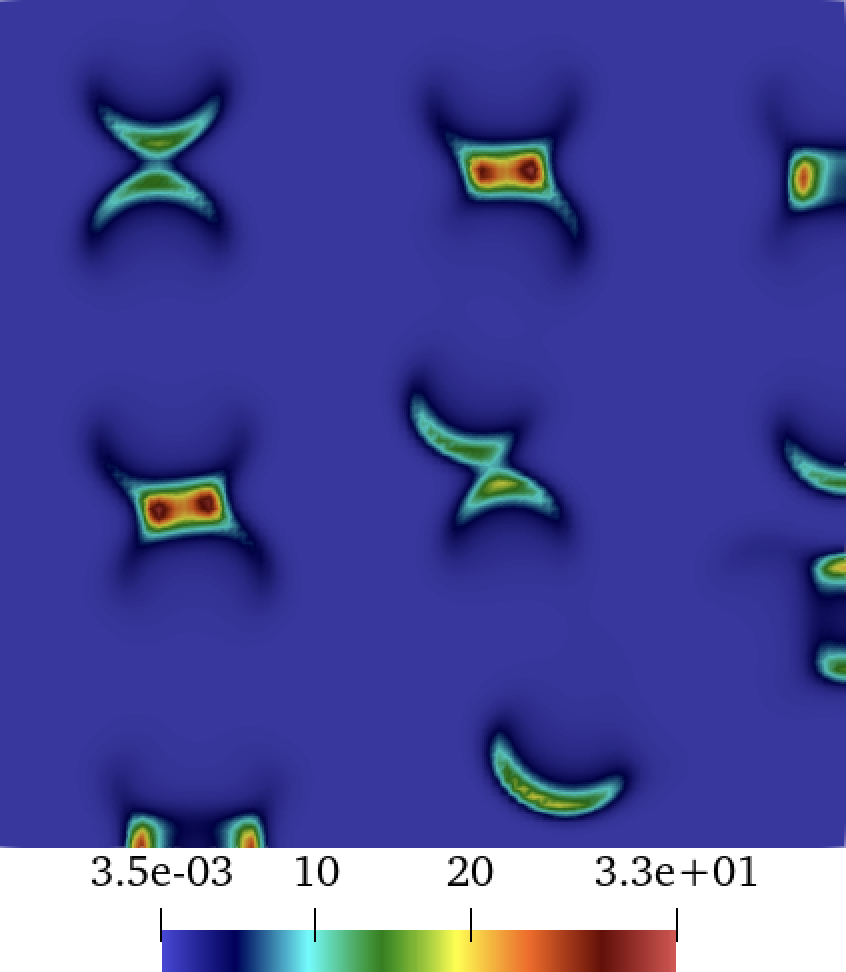}
       \put(23,102){\footnotesize{$t=1369$}}
    \end{overpic} 
  \begin{overpic}[percent,width=0.16\textwidth]{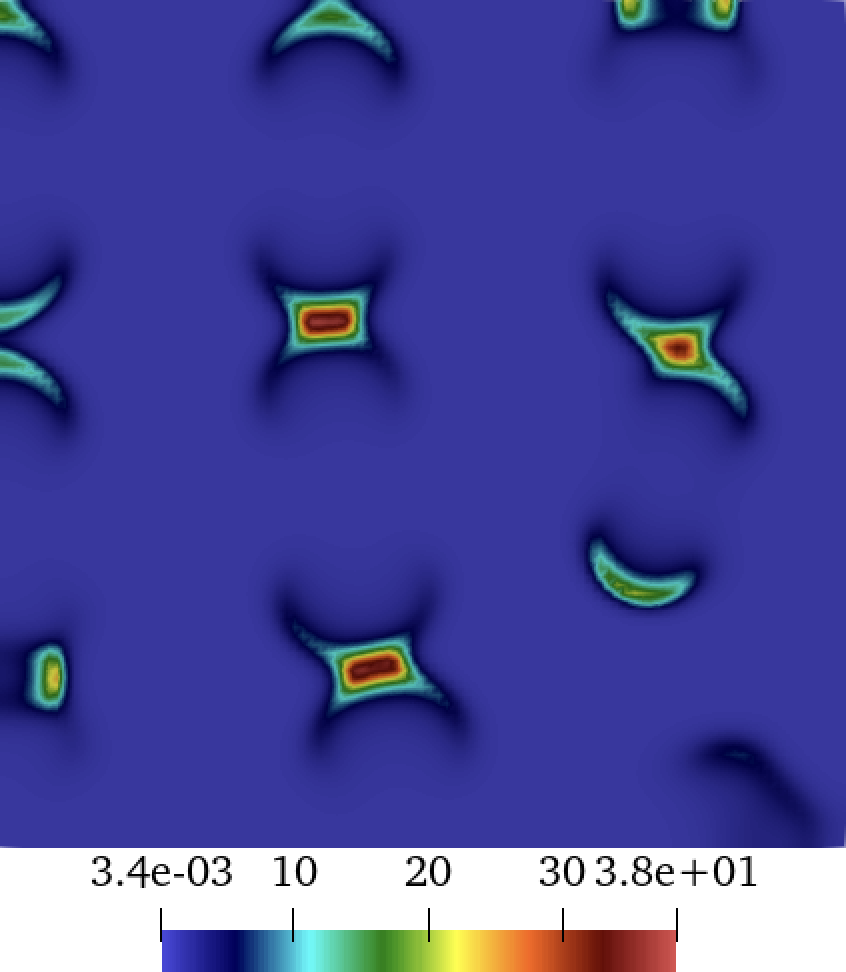}
      \put(23,102){\footnotesize{$t=1378$}}
    \end{overpic} 

    \caption{Case 4: expected number of crimes per unit area $S$ from $t = 1360$ to $t = 1378$. }
    \label{fig:eta0075_snap}
\end{figure}

Figs.~\ref{fig:case2_PDE}-\ref{fig:eta0075_snap} show that, when all other parameters are fixed, $\eta$ plays a key role in the stability of the system, controlling the excited frequencies and hotspot size when the solution is oscillatory. 
Next, we fix $\eta$ to one of the values used in Fig.~\ref{fig:eta}, i.e., $\eta = 0.15$, and vary the value of  
$\Gamma \theta(1-\Sigma)/\omega^2$
with cases 5 and 6, which belong to the oscillatory part of the plane in Fig.~\ref{fig:phase_diagrams} (center). 
Fig.~\ref{fig:bbar} shows the evolution of $\langle S\rangle$ for
cases 3, 5, and 6, and associated power spectra. Smaller $\Gamma \theta(1-\Sigma)/\omega^2$ means fewer criminals (either because they ``regenerate'' less or because they get arrested more) and/or longer time scale for repeat victimization
and/or weaker repeat victimization. 
Thus, it is not surprising to see in 
Fig.~\ref{fig:bbar} that the average 
number of crimes over time lowers as 
$\Gamma \theta(1-\Sigma)/\omega^2$
becomes smaller. However, it is interesting to see how the number of excited frequencies varies as
$\Gamma \theta(1-\Sigma)/\omega^2$ is varied. Note that cases 3 and 6 (i.e., $\Gamma \theta(1-\Sigma)/\omega^2=1.5, 2.5$) have two respective peaks very close to each other, one of them being the dominant peak. Those peaks disappear in case 5 (i.e., $\Gamma \theta(1-\Sigma)/\omega^2=0.5$), which gives more irregular, low amplitude oscillations.
Furthermore, the hotspots in this case, shown in Fig.~\ref{fig:bbar25_snap},
seem to mix the dynamics observed in Figs.~\ref{fig:eta015_snap} and 
\ref{fig:eta0075_snap}.

\begin{figure}[htb!]
    \centering
    \includegraphics[width=1.0\linewidth]{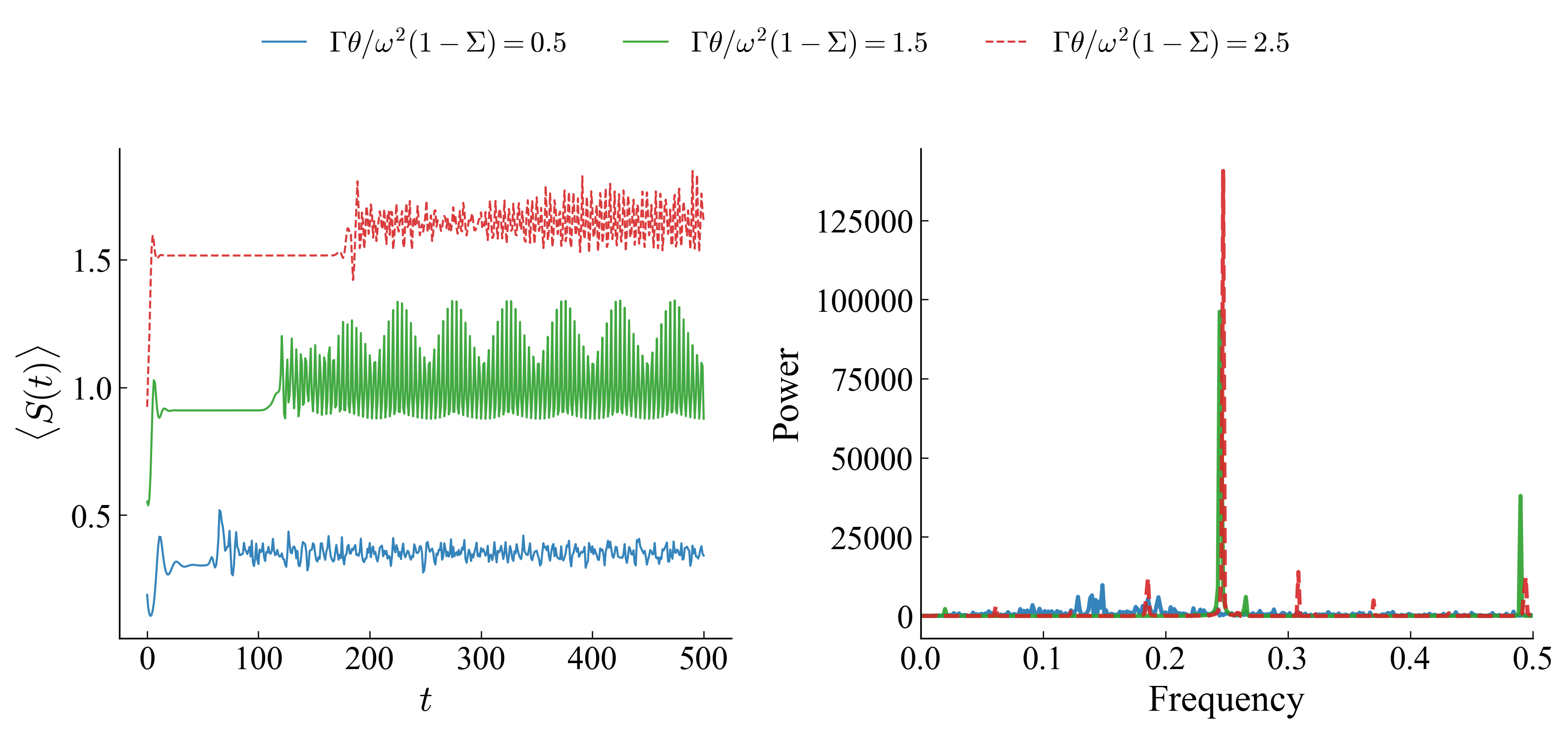}
    \caption{Left: time evolutions of spatially averaged quantities $\langle S \rangle$ for case 5 ($\Gamma \theta(1-\Sigma)/\omega^2 = 0.5$), case 3 ($\Gamma \theta(1-\Sigma)/\omega^2 = 1.5$), and case 6 ($\Gamma \theta(1-\Sigma)/\omega^2 = 2.5$). Right: associated power spectra.}
    \label{fig:bbar}
\end{figure}

\begin{figure}[htb!]
     \centering
         \begin{overpic}[percent,width=0.16\textwidth]{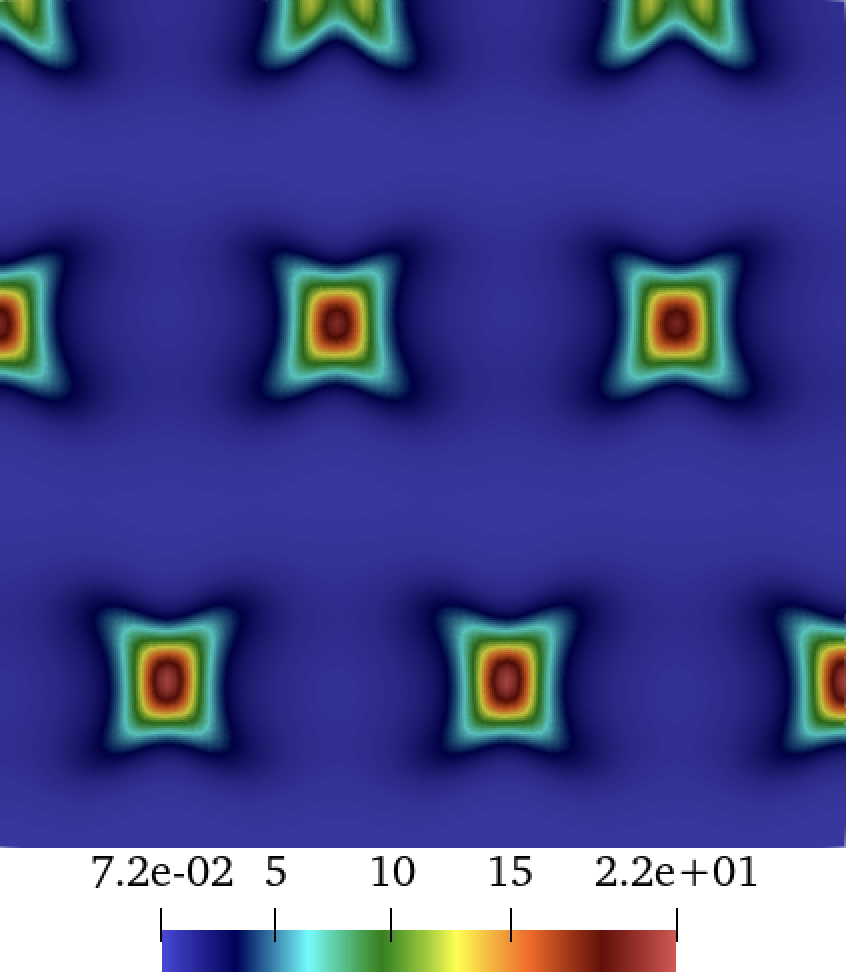}
     \put(23,102){\footnotesize{$t=673$}}
    \end{overpic} 
     \begin{overpic}[percent,width=0.16\textwidth, grid=false]{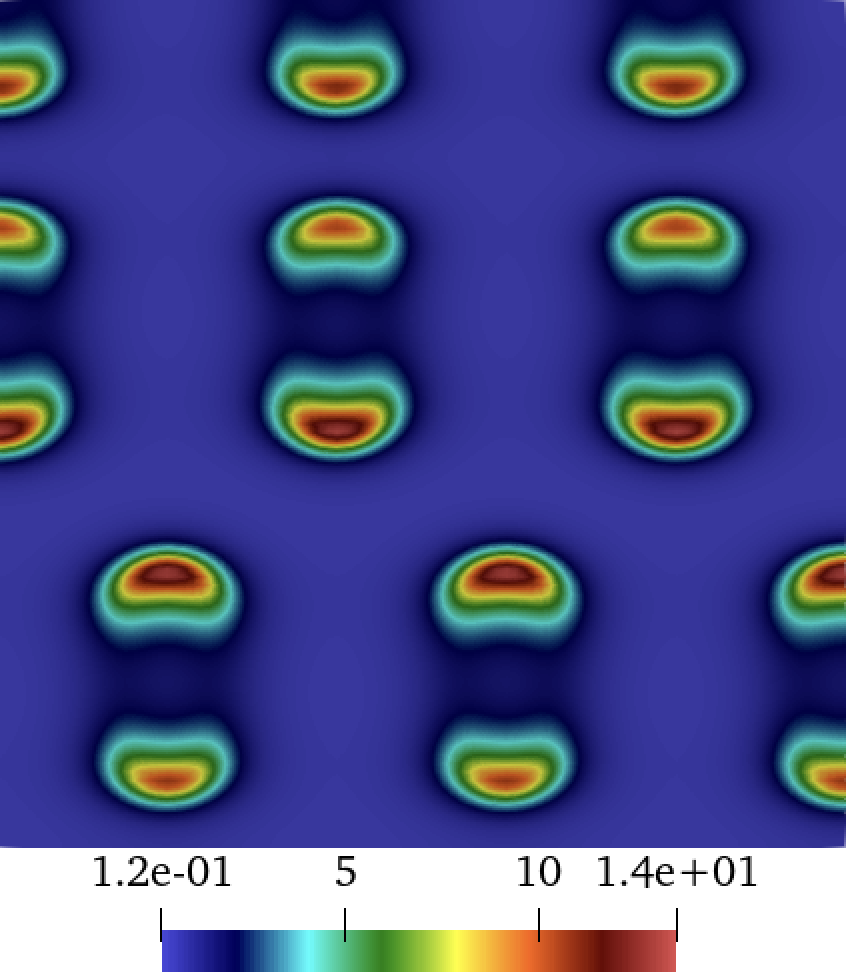}
     \put(23,102){\footnotesize{$t=675$}}
    \end{overpic}
    \begin{overpic}[percent,width=0.16\textwidth]{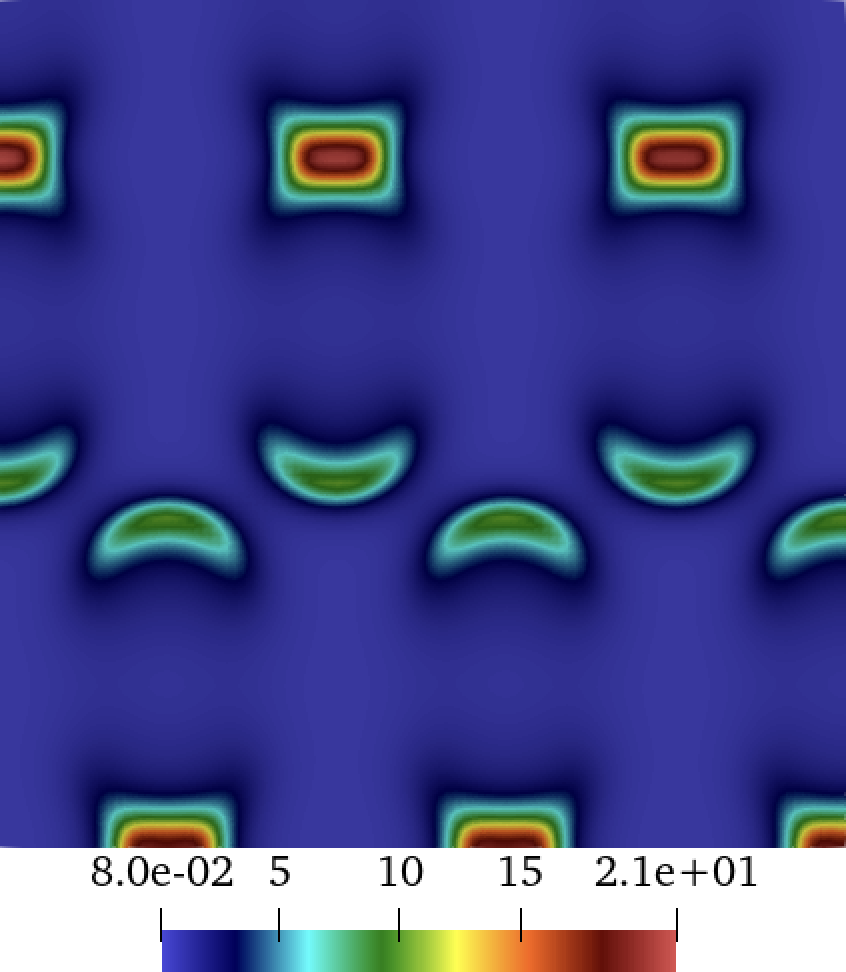}
  \put(23,102){\footnotesize{$t=677$}}
    \end{overpic} 
        \begin{overpic}[percent,width=0.16\textwidth]{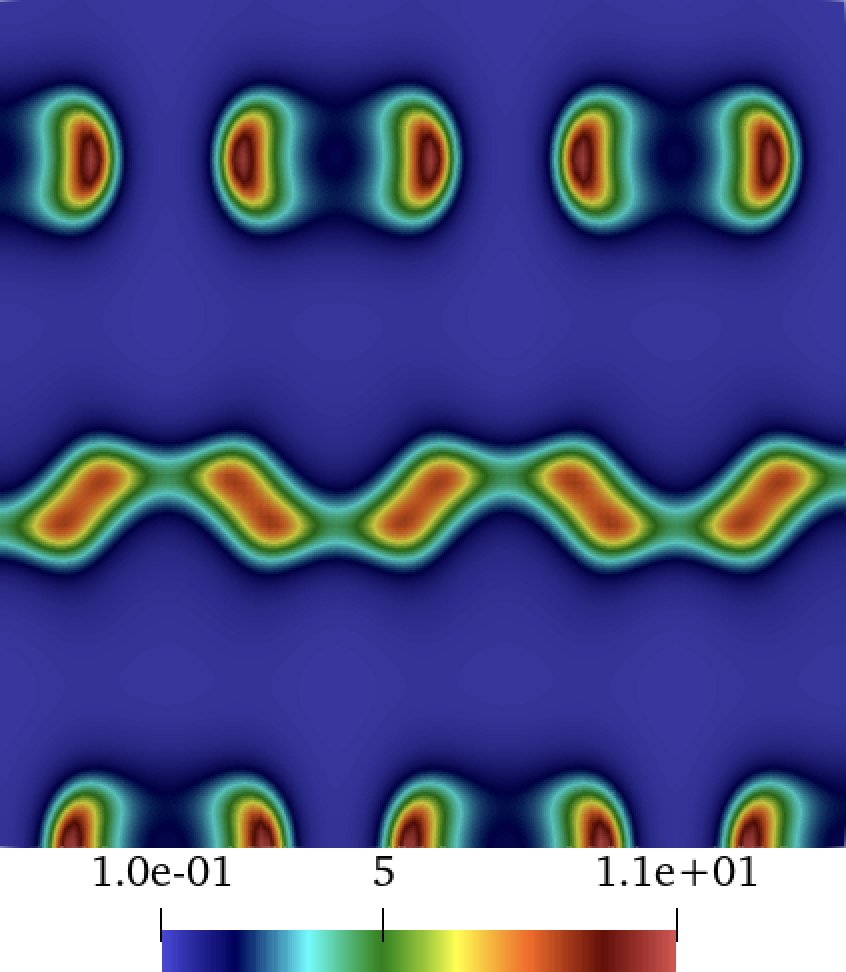}
       \put(23,102){\footnotesize{$t=679$}}
    \end{overpic} 
      \begin{overpic}[percent,width=0.16\textwidth]{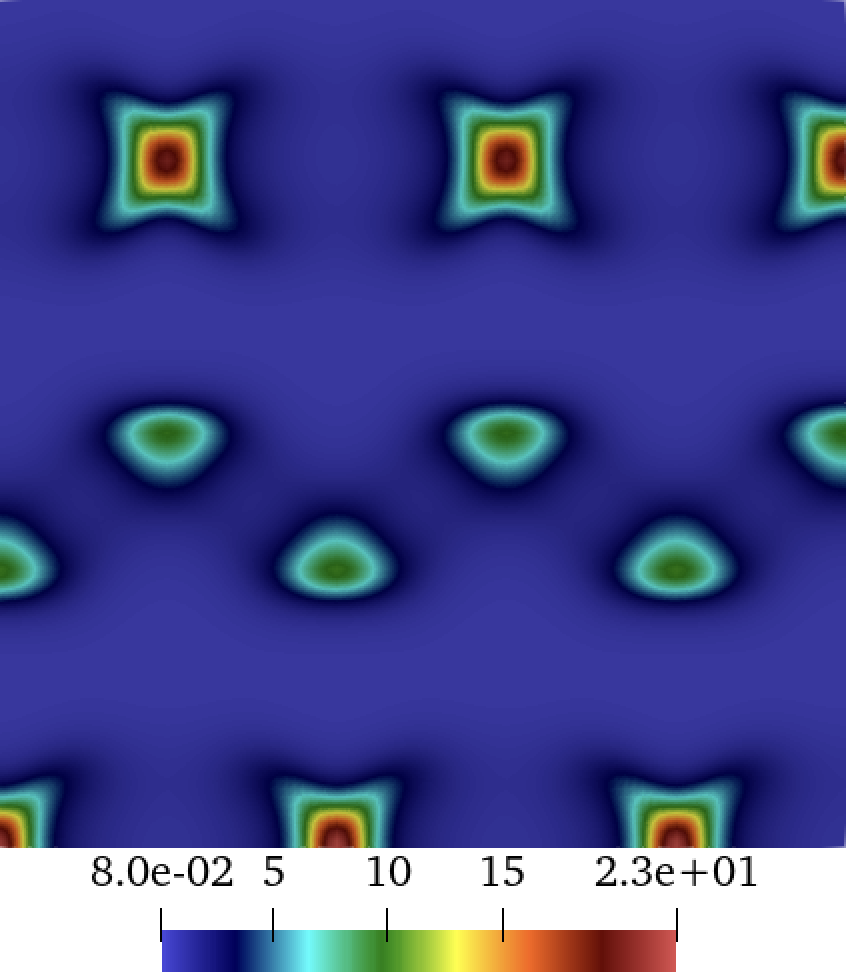}
       \put(23,102){\footnotesize{$t=681$}}
    \end{overpic} 
  \begin{overpic}[percent,width=0.16\textwidth]{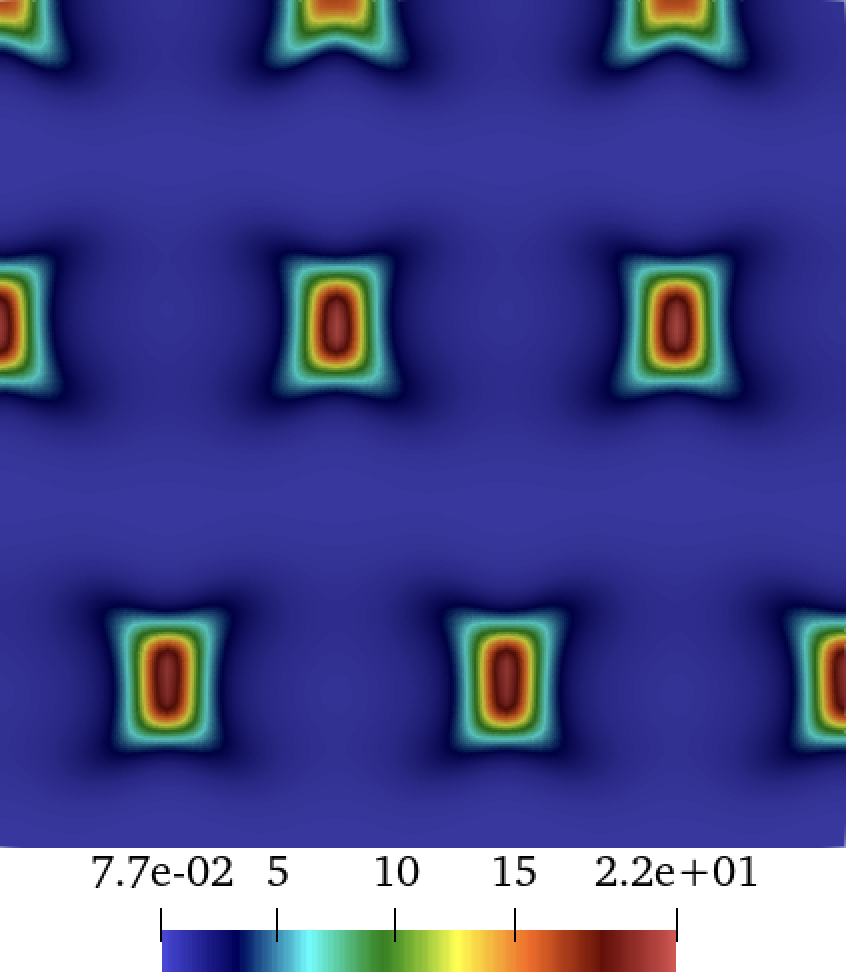}
      \put(23,102){\footnotesize{$t=689$}}
    \end{overpic} 

    \caption{Case 6: expected number of crimes per unit area $S$ from $t = 673$ to $t = 689$.}
    \label{fig:bbar25_snap}
\end{figure}

Let us now fix $\eta = 0.15$ and 
$\Gamma \theta(1-\Sigma)/\omega^2 = 1.5$, and vary $\tau$, 
which represents the time lag between crime occurrence and the availability of crime data used by police to plan patrols.
We consider cases 7 and 8, which lie in the oscillatory part of the plane in Fig.~\ref{fig:phase_diagrams} (left). 
Fig.~\ref{fig:tau}
shows the evolution of $\langle S\rangle$ for cases 3, 7, and 8, and associated power spectra. For very short delays ($\tau = 0.5$), $\langle S\rangle$ rapidly converges to a constant value, which is homogeneous in space (not shown for the sake of brevity). This indicates that when police have prompt access to crime data, their response is effective at stabilizing crime overall. For moderate delays ($\tau = 5$), we see sustained oscillations with relatively large amplitude. The delayed response causes 
the policing strategy to systematically lag behind current crime patterns, leading to cycles of over- and under-enforcement. As a result, crime levels periodically increase and decrease, rather than settling to a steady equilibrium. 
For long delays ($\tau = 50$), the long-term dynamics show irregular, low-amplitude fluctuations around the mean. In this regime, policing decisions are based on outdated information, effectively decoupling enforcement from current crime patterns. 
This is reflected in the 
hotspots, which are irregular in shape and display an erratic motion.
See Fig.~\ref{fig:tau50_snap}.

\begin{figure}[htb!]
    \centering
    \includegraphics[width=1.0\linewidth]{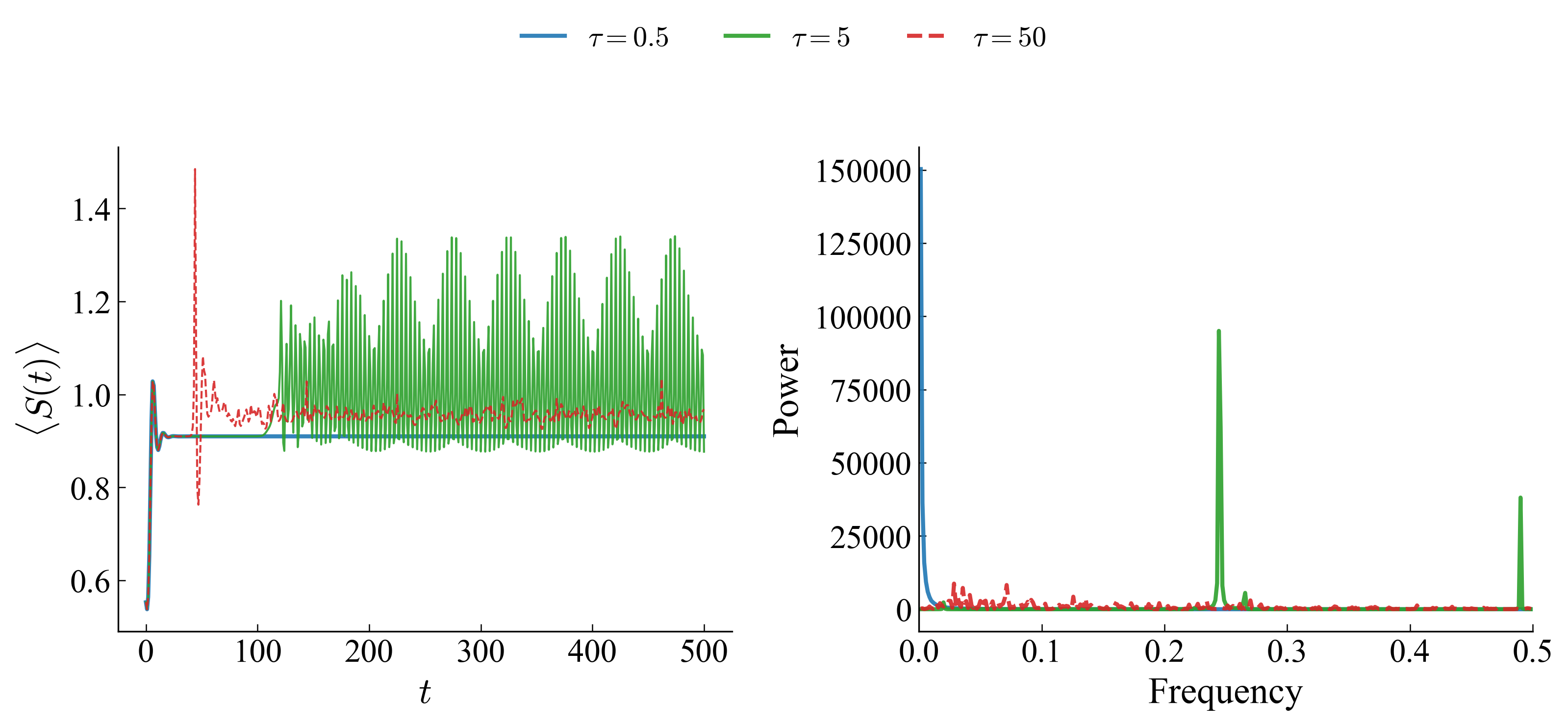}
    \caption{Left: time evolutions of spatially averaged quantities $\langle S \rangle$ for case 7 ($\tau = 0.5$), case 3 ($\tau = 5$), and case 8 ($\tau = 50$). Right: associated power spectra.
    }
    \label{fig:tau}
\end{figure}

\begin{figure}[htb!]
     \centering
         \begin{overpic}[percent,width=0.16\textwidth]{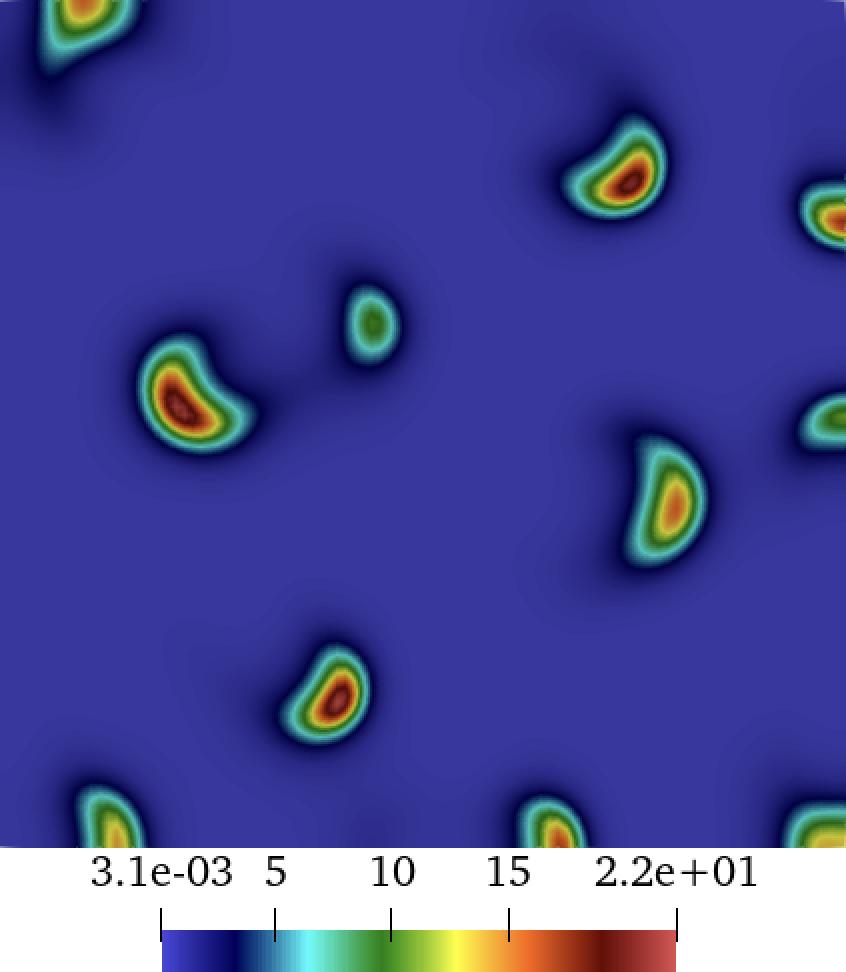}
     \put(23,102){\footnotesize{$t=765$}}
    \end{overpic} 
     \begin{overpic}[percent,width=0.16\textwidth, grid=false]{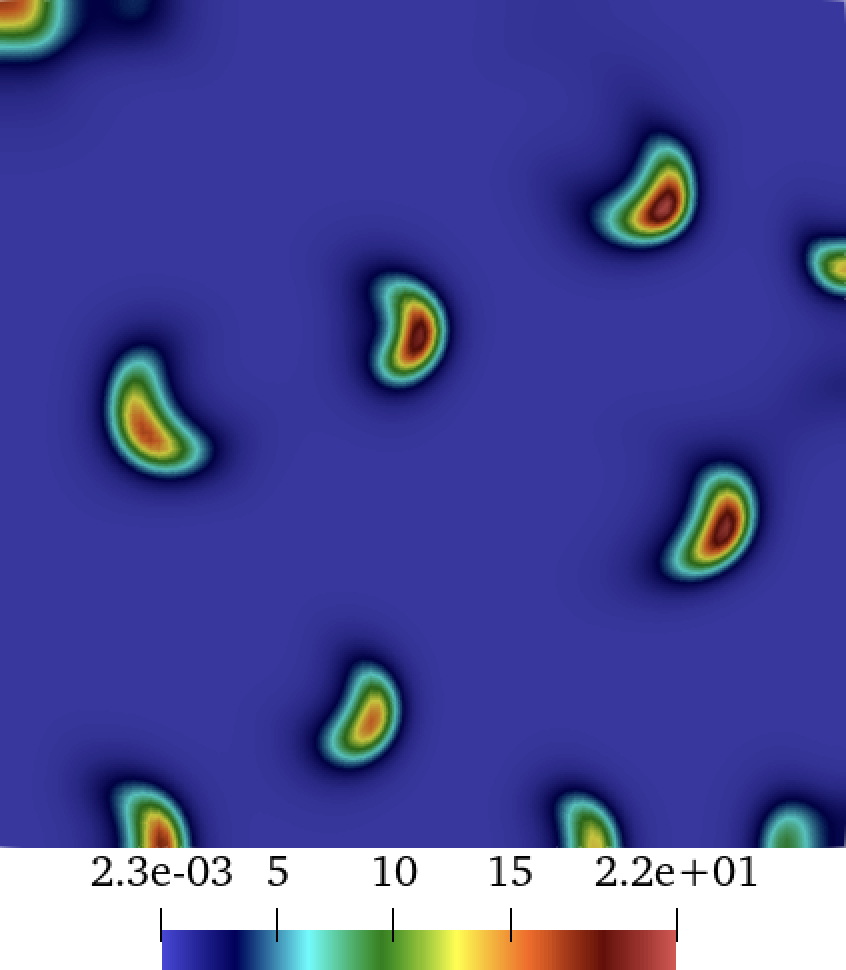}
     \put(23,102){\footnotesize{$t=768$}}
    \end{overpic}
    \begin{overpic}[percent,width=0.16\textwidth]{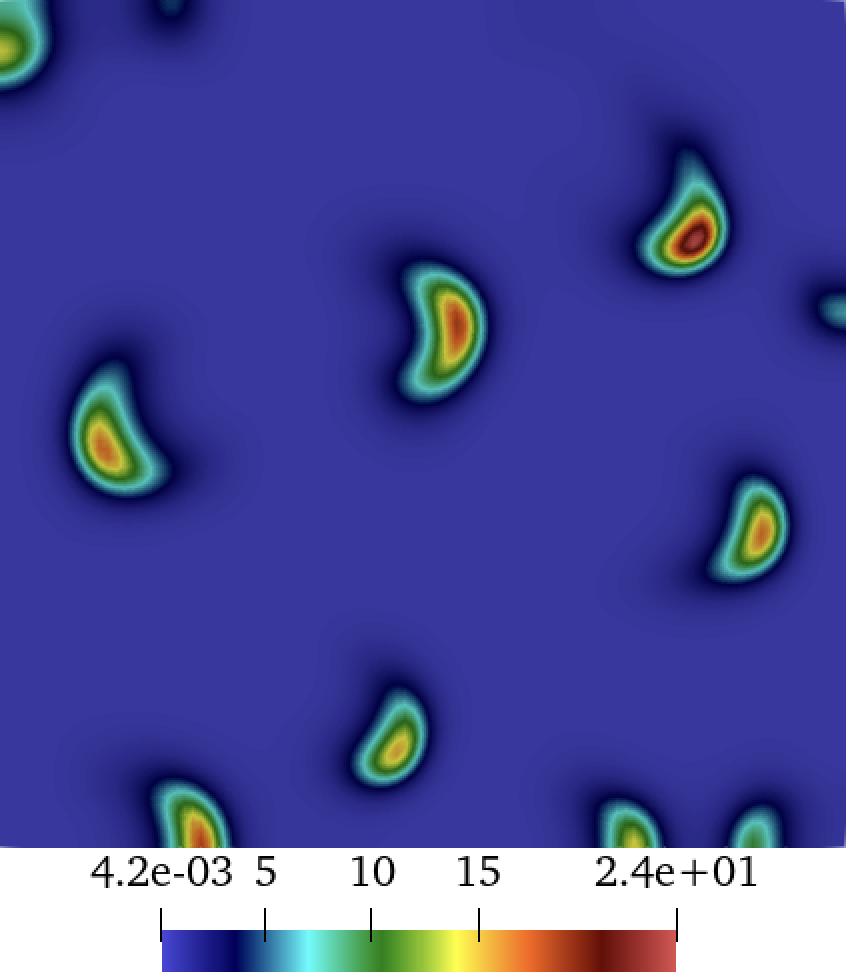}
  \put(23,102){\footnotesize{$t=771$}}
    \end{overpic} 
        \begin{overpic}[percent,width=0.16\textwidth]{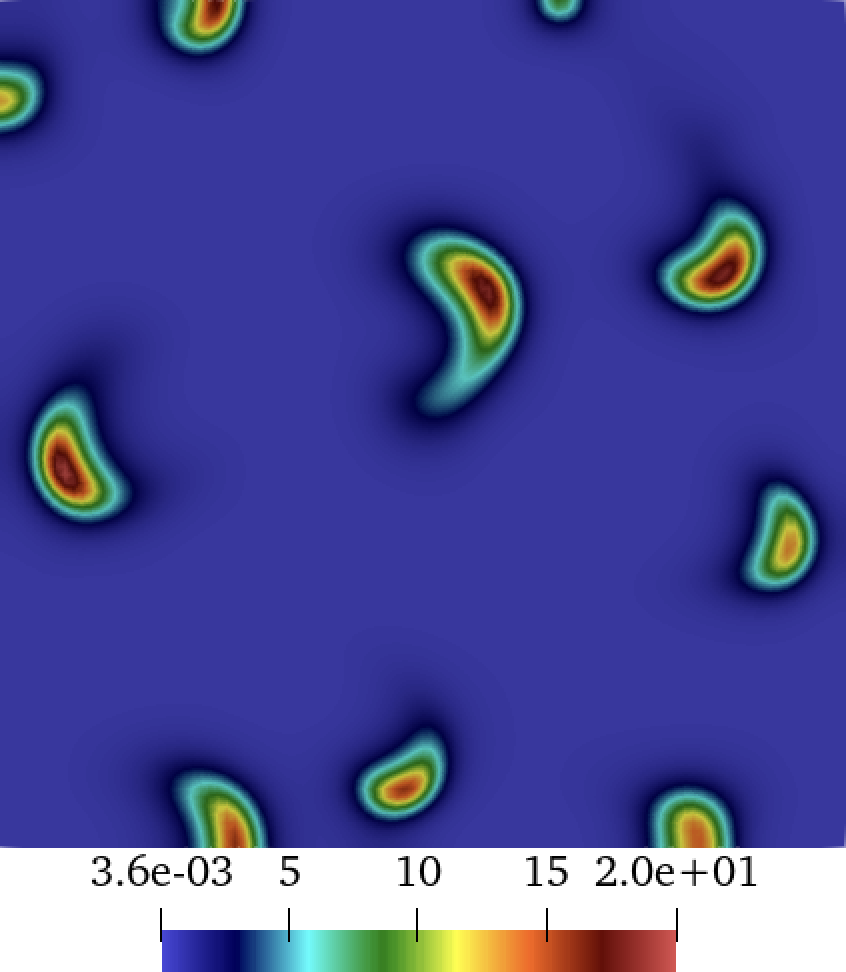}
       \put(23,102){\footnotesize{$t=774$}}
    \end{overpic} 
      \begin{overpic}[percent,width=0.16\textwidth]{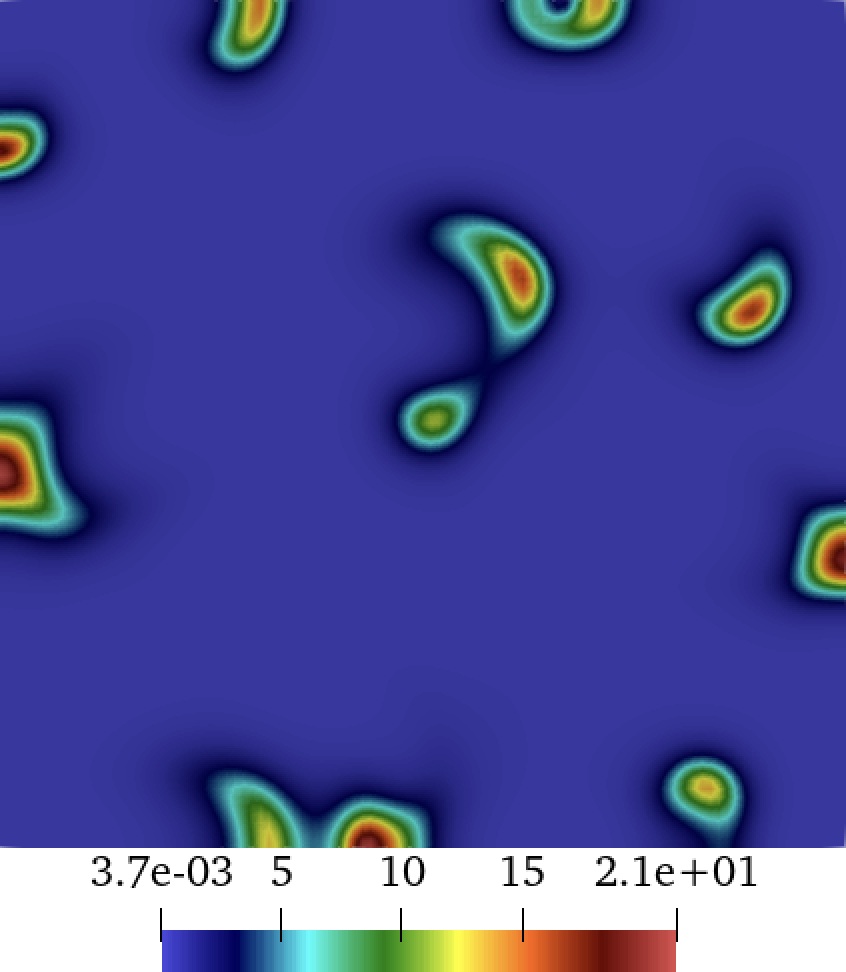}
       \put(23,102){\footnotesize{$t=777$}}
    \end{overpic} 
  \begin{overpic}[percent,width=0.16\textwidth]{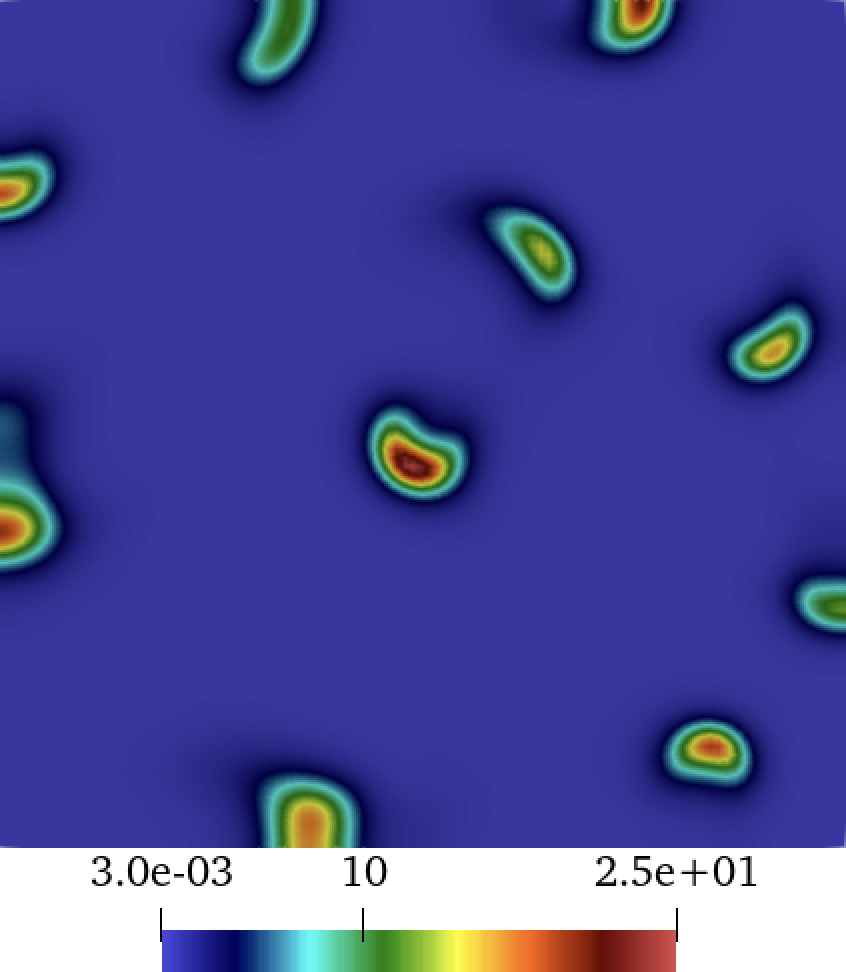}
      \put(23,102){\footnotesize{$t=780$}}
    \end{overpic} 

    \caption{Case 8: expected number of crimes per unit area $S$ from $t = 765$ to $t = 780$.}
    \label{fig:tau50_snap}
\end{figure}

For comparison with Fig.~\ref{fig:tau50_snap}, 
Fig.~\ref{fig:tau50_snapABM} shows
the evolution of $S$ for case 8 given by the agent-based solver after non-dimensionalization. We set 
$\Gamma = 0.0285$, $\theta = 0.2339$, $\Sigma = 0$ and $\omega = 1/15$ to obtain $\Gamma \theta (1 - \Sigma)/\omega^2 \approx 1.5$.
We see that the magnitudes of 
$S$ in Fig.~\ref{fig:tau50_snap}
and \ref{fig:tau50_snapABM} match
and the size and number of hotspots
are comparable. 

\begin{figure}[htb!]
     \centering
         \begin{overpic}[percent,width=0.16\textwidth]{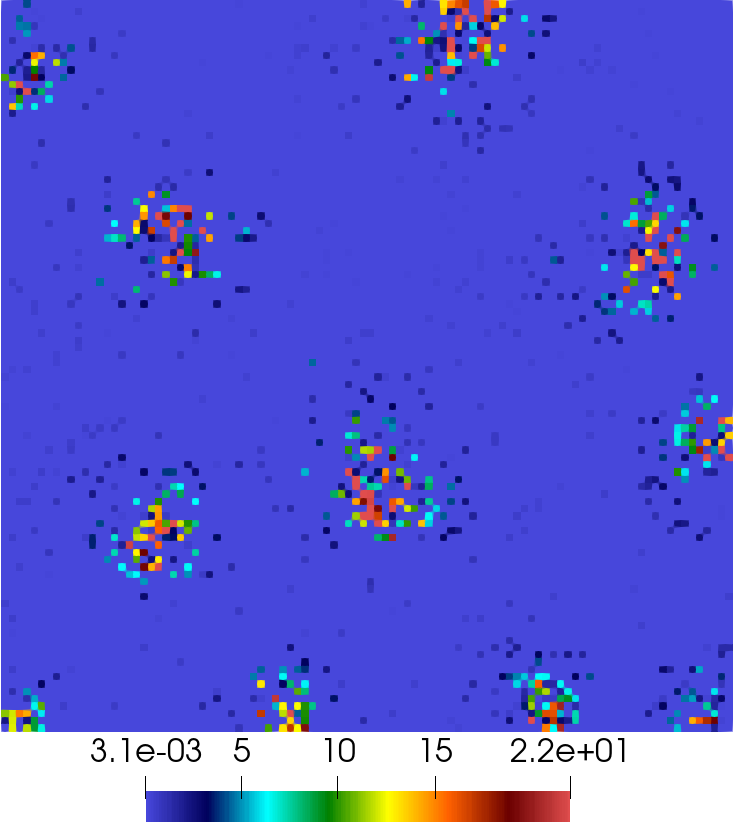}
     \put(23,102){\footnotesize{$t=765$}}
    \end{overpic} 
     \begin{overpic}[percent,width=0.16\textwidth, grid=false]{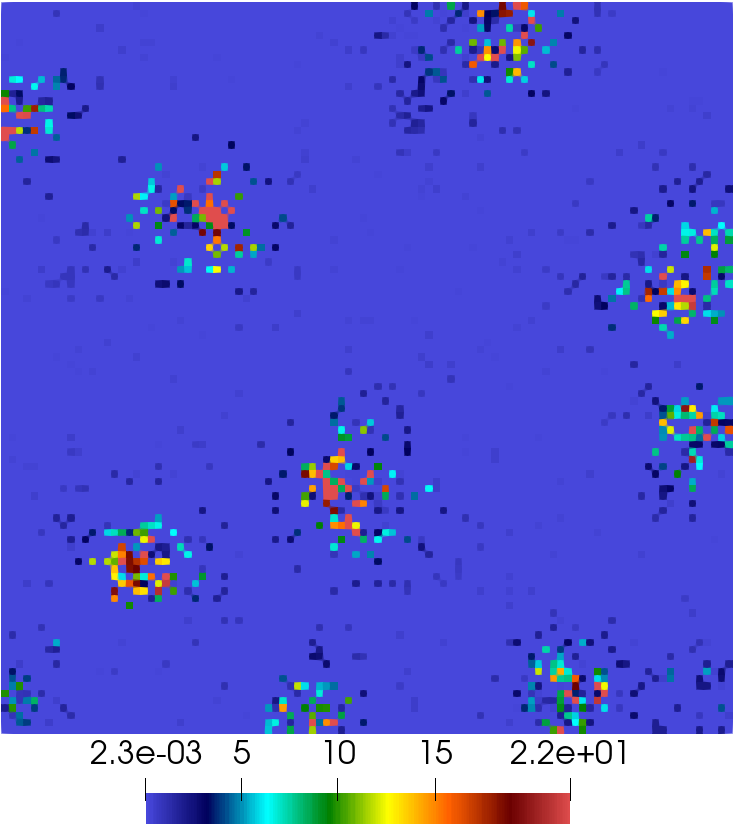}
     \put(23,102){\footnotesize{$t=768$}}
    \end{overpic}
    \begin{overpic}[percent,width=0.16\textwidth]{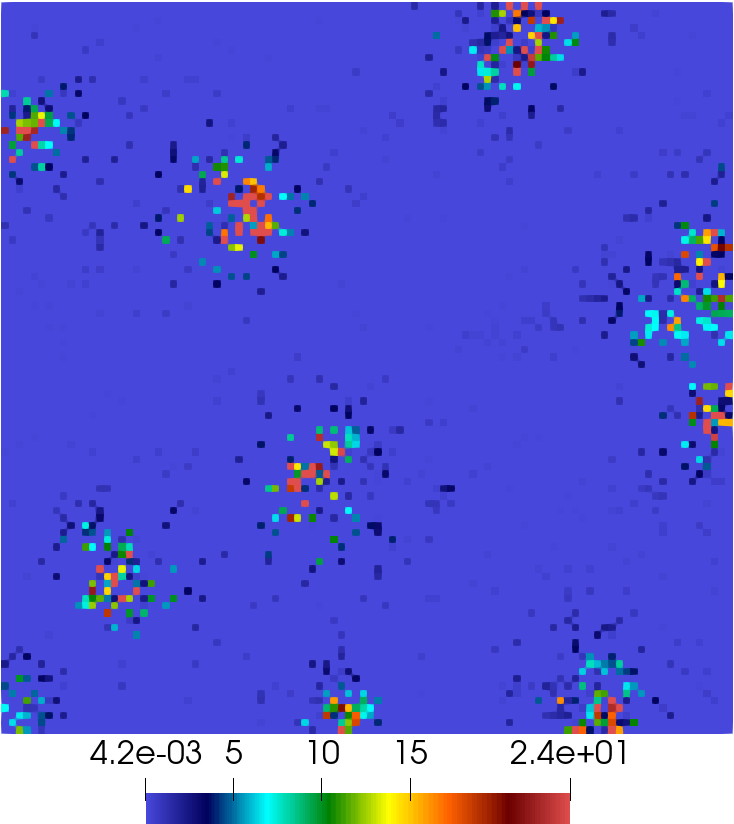}
  \put(23,102){\footnotesize{$t=771$}}
    \end{overpic} 
        \begin{overpic}[percent,width=0.16\textwidth]{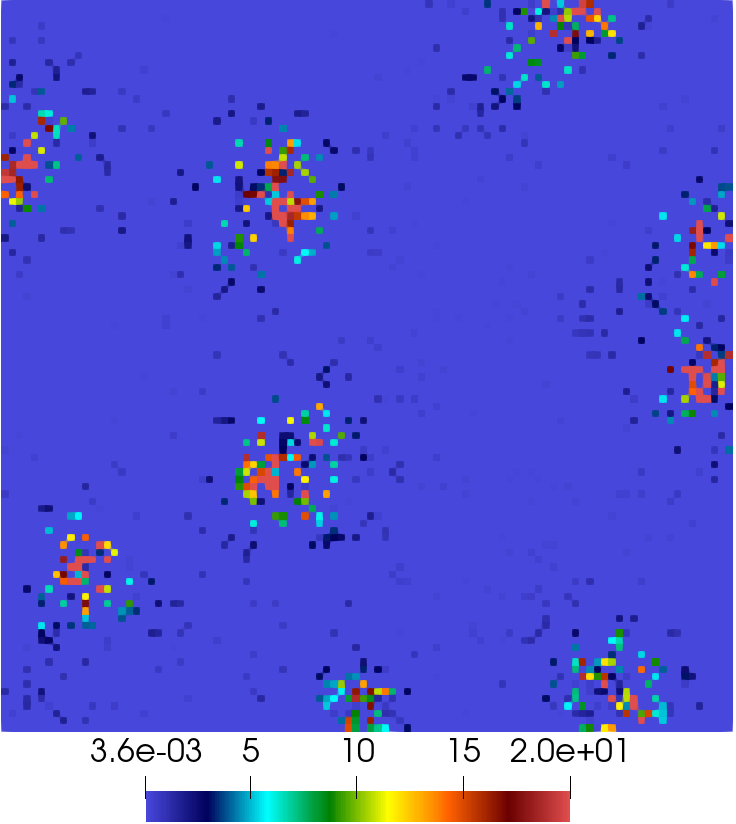}
       \put(23,102){\footnotesize{$t=774$}}
    \end{overpic} 
      \begin{overpic}[percent,width=0.16\textwidth]{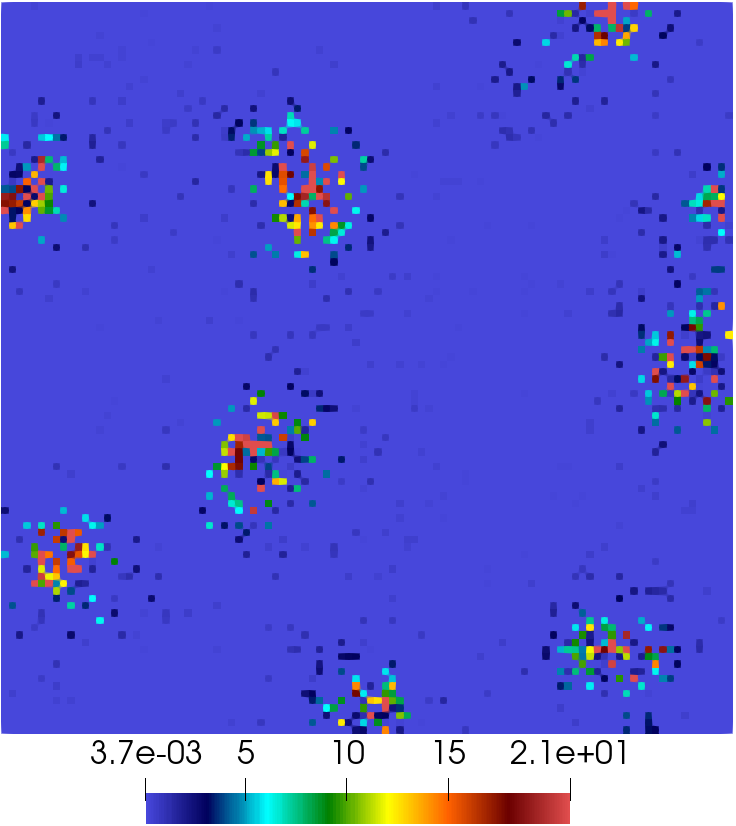}
       \put(23,102){\footnotesize{$t=777$}}
    \end{overpic} 
  \begin{overpic}[percent,width=0.16\textwidth]{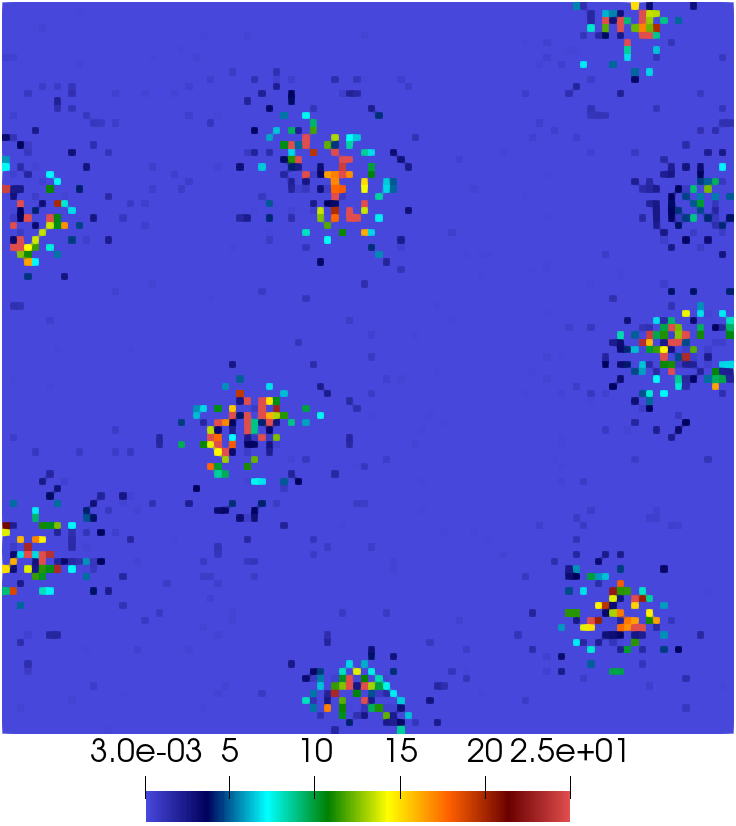}
      \put(23,102){\footnotesize{$t=780$}}
    \end{overpic} 

    \caption{Case 8 - agent-based model: expected number of crimes per unit area $S$ from $t = 765$ to $t = 780$.}
    \label{fig:tau50_snapABM}
\end{figure}

So far, we have kept $\pi_0 = \langle \pi \rangle = 0.5$. Let us now investigate the effect on crime of increasing or decreasing police density with cases 9 and 10. 
Fig.~\ref{fig:pi} shows the evolution of $\langle S\rangle$ for cases 3, 9, and 10
and associated power spectra.
For low police density ($\pi_0 = 0.1$), $\langle S\rangle$ quickly converges to a high, nearly constant value, with only small fluctuations around the mean. 
This indicates that policing is insufficient to significantly suppress crime. Due to weak crime-police feedback, crime patterns evolve largely unhindered, leading to an elevated crime level and large crime hotspots (see Fig.~\ref{fig:pi01_snap}). At intermediate police density ($\pi_0 = 0.5$), $\langle S\rangle$ exhibits sustained oscillations, i.e., the police-crime
interaction produces cyclical dynamics of periods of effective suppression followed by rebounds in crime. Higher police density ($\pi_0 = 1$), reduces
the average number of expected crimes, but temporal fluctuations persist and are more irregular. In fact, in Fig.~\ref{fig:pi} (right) we see more excited frequencies
in the spectrum for $\pi_0 = 1$ than for 
$\pi_0 = 0.5$. This suggests that, while high enforcement suppresses crime on average, it can also
induce complex transient dynamics
due to rapid displacement and fragmentation of hotspots (see Fig.~\ref{fig:pi1_snap}).

\begin{figure}[htb!]
    \centering
    \includegraphics[width=\linewidth]{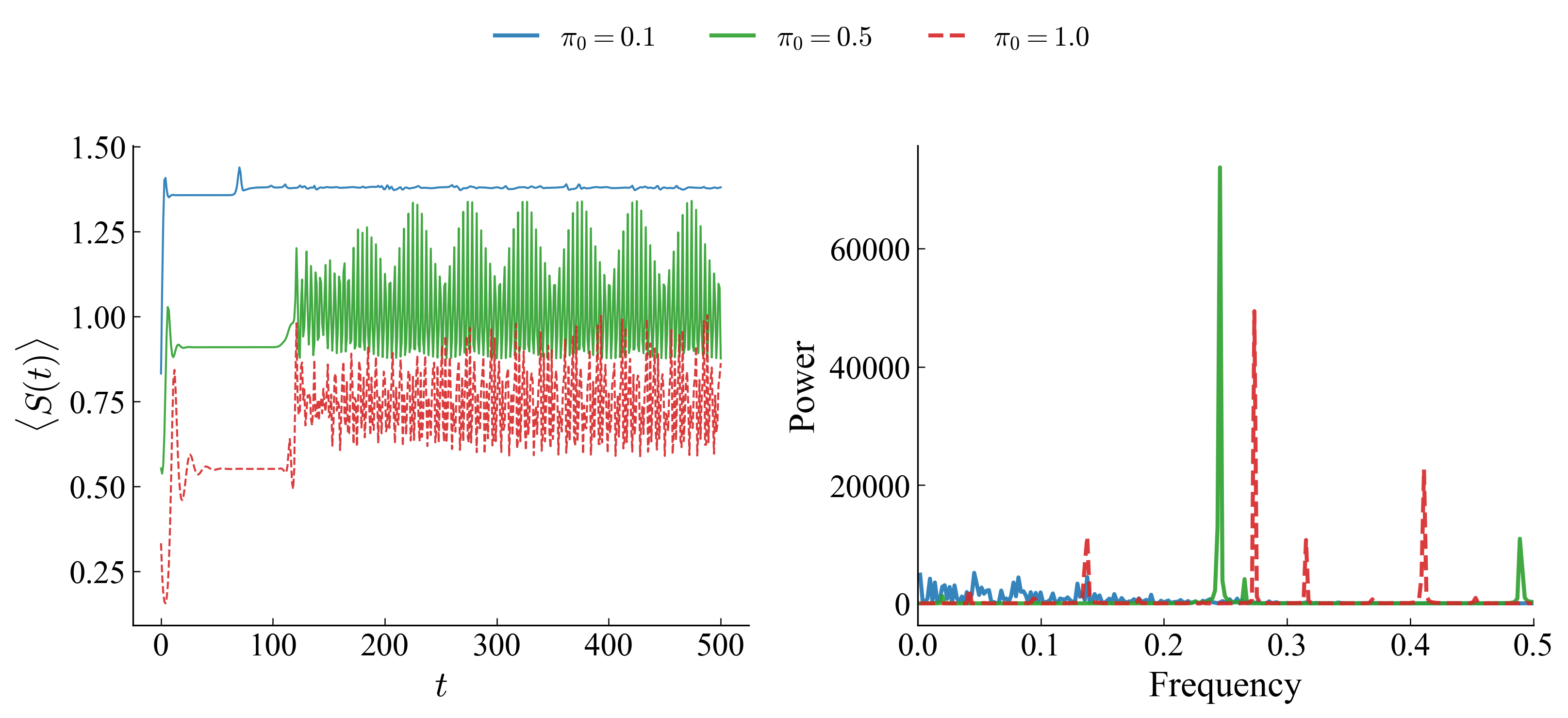}
    \caption{Left:
    time evolutions of spatially averaged quantities $\langle S \rangle$ for case 9 ($\pi_0 = 0.1$), case 3 ($\pi_0 = 0.5$), and case 10 ($\pi_0 = 1$). Right: associated power spectra.
    }
    \label{fig:pi}
\end{figure}

\begin{figure}[htb!]
     \centering
         \begin{overpic}[percent,width=0.16\textwidth]{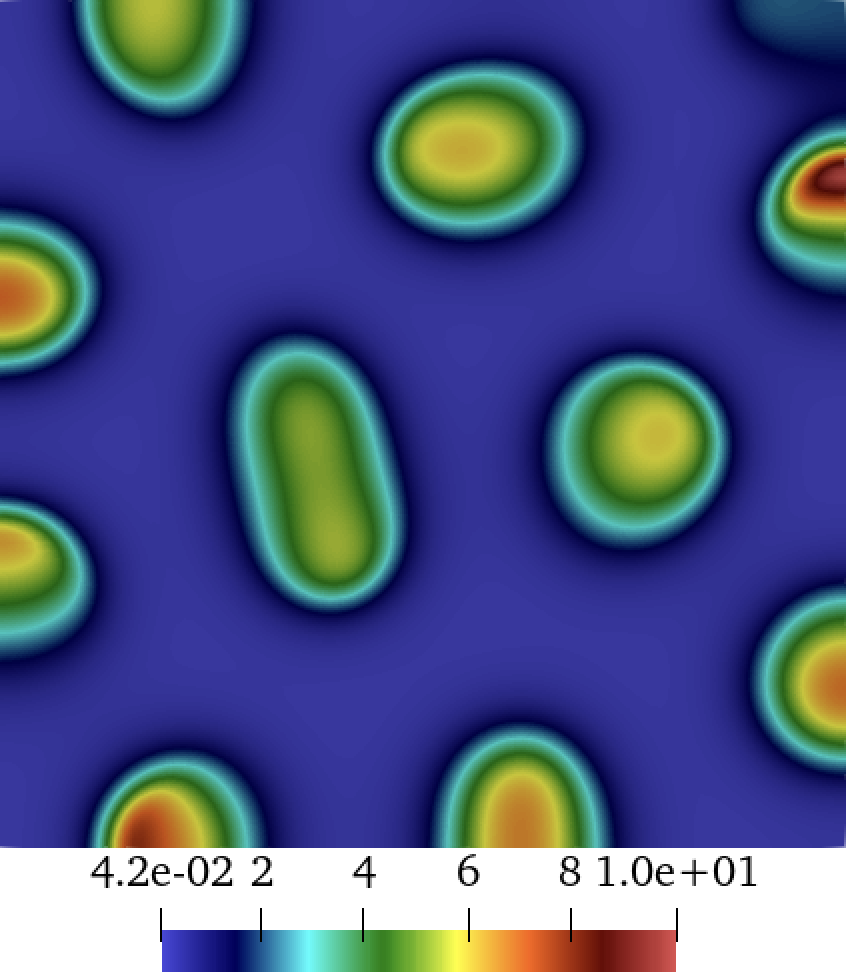}
     \put(23,102){\footnotesize{$t=744$}}
    \end{overpic} 
     \begin{overpic}[percent,width=0.16\textwidth, grid=false]{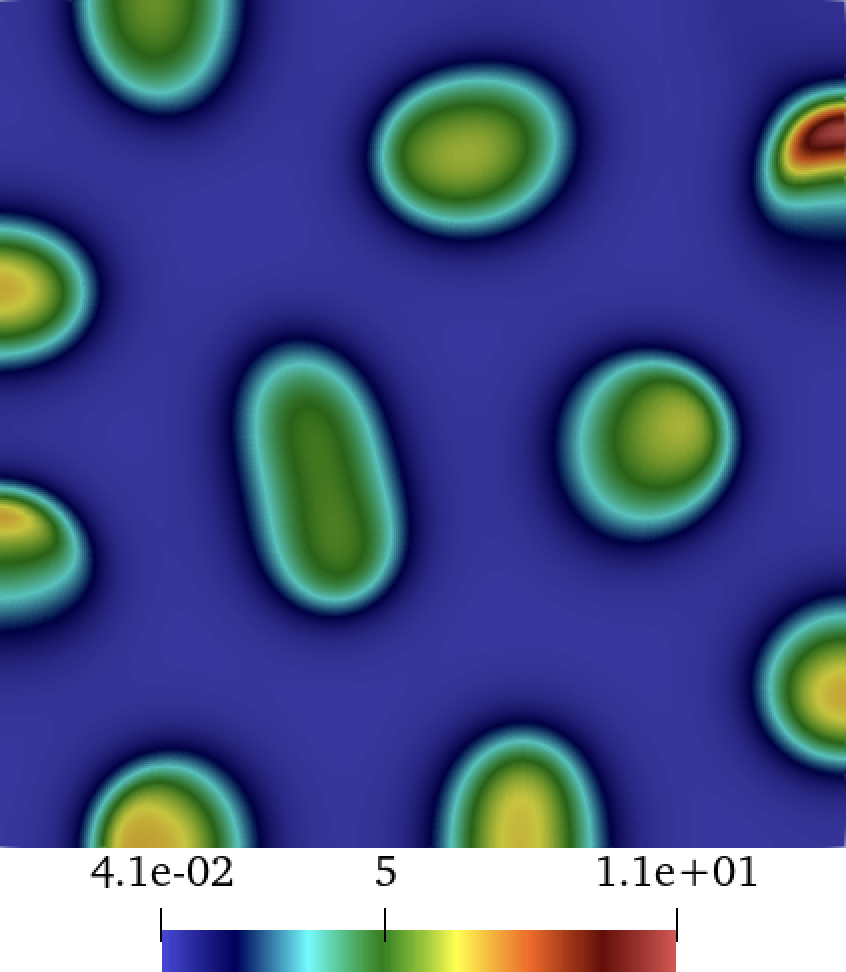}
     \put(23,102){\footnotesize{$t=747$}}
    \end{overpic}
    \begin{overpic}[percent,width=0.16\textwidth]{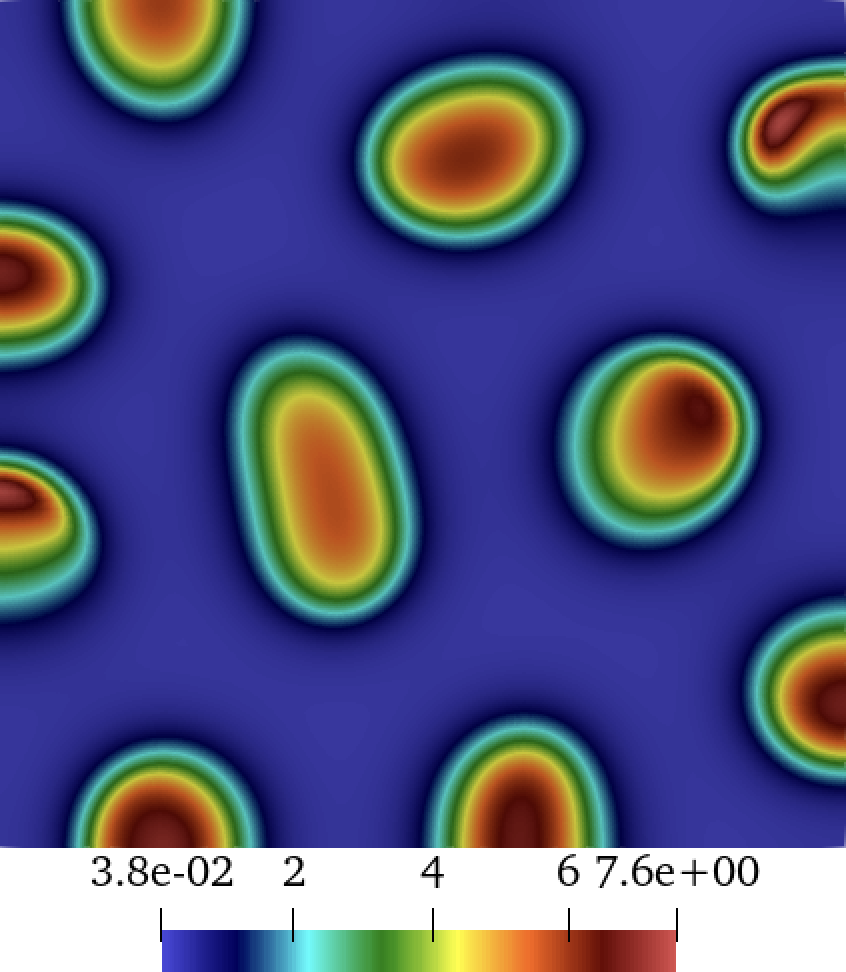}
  \put(23,102){\footnotesize{$t=750$}}
    \end{overpic} 
        \begin{overpic}[percent,width=0.16\textwidth]{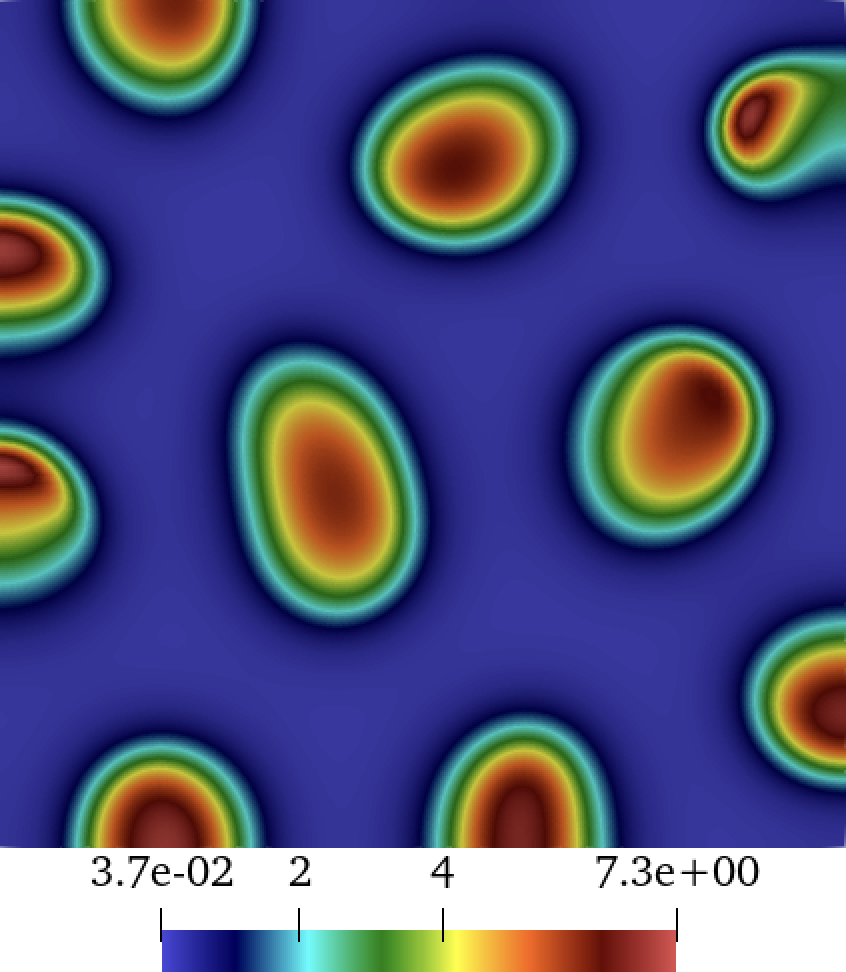}
       \put(23,102){\footnotesize{$t=753$}}
    \end{overpic} 
      \begin{overpic}[percent,width=0.16\textwidth]{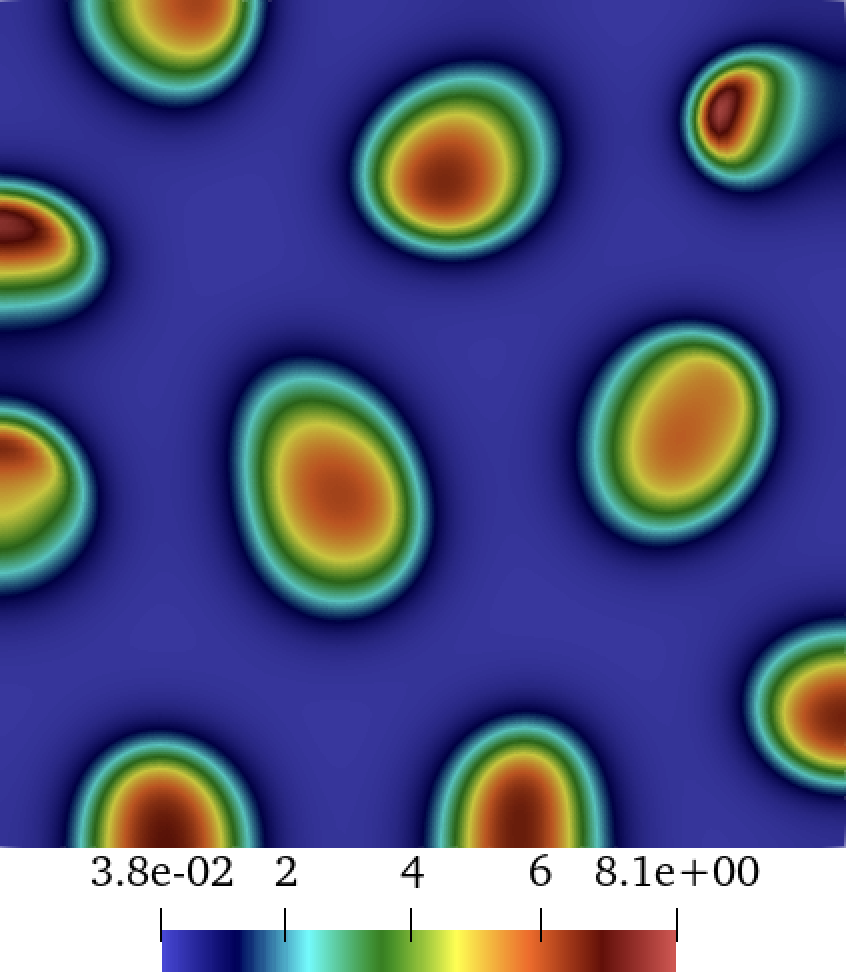}
       \put(23,102){\footnotesize{$t=756$}}
    \end{overpic} 
  \begin{overpic}[percent,width=0.16\textwidth]{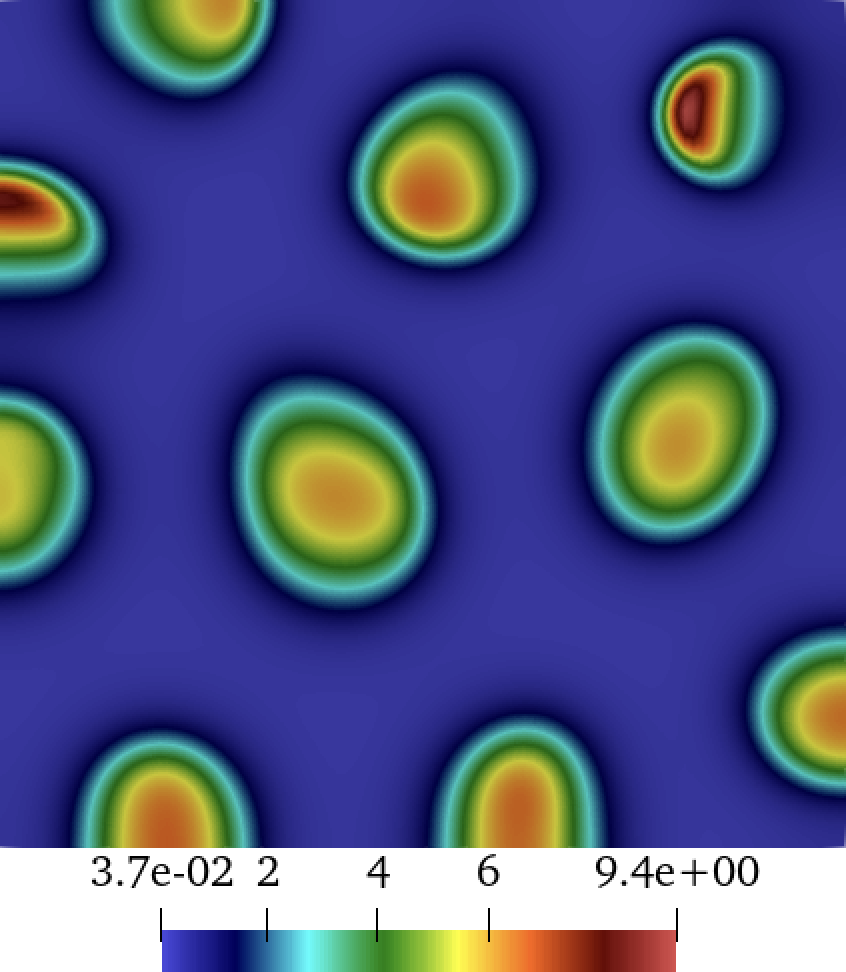}
      \put(23,102){\footnotesize{$t=759$}}
    \end{overpic} 

    \caption{Case 9: expected number of crimes per unit area $S$ from $t = 744$ to $t = 759$.}
    \label{fig:pi01_snap}
\end{figure}

\begin{figure}[htb!]
     \centering
         \begin{overpic}[percent,width=0.16\textwidth]{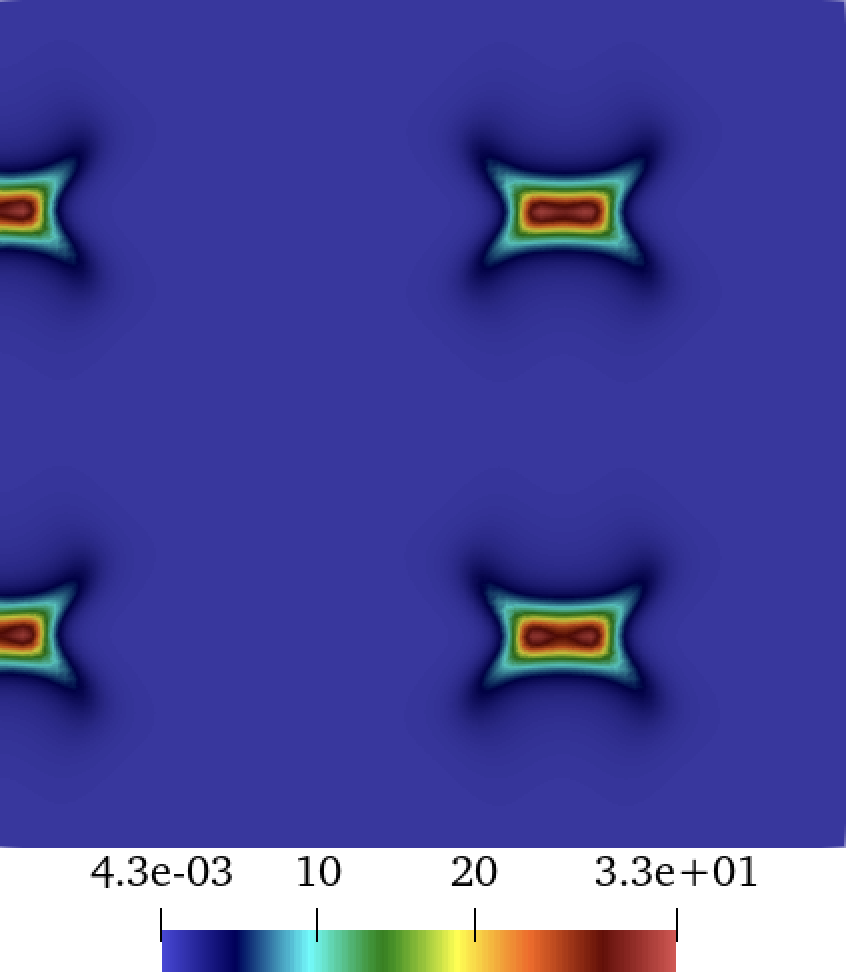}
     \put(23,102){\footnotesize{$t=773$}}
    \end{overpic} 
     \begin{overpic}[percent,width=0.16\textwidth, grid=false]{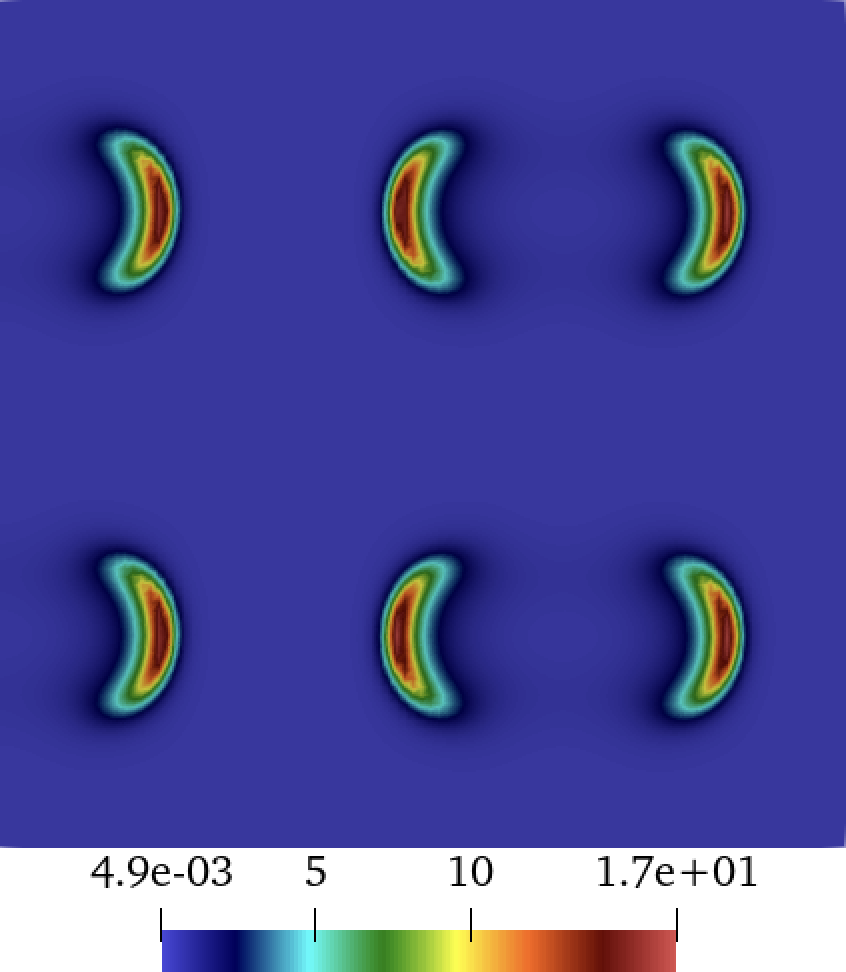}
     \put(23,102){\footnotesize{$t=775$}}
    \end{overpic}
    \begin{overpic}[percent,width=0.16\textwidth]{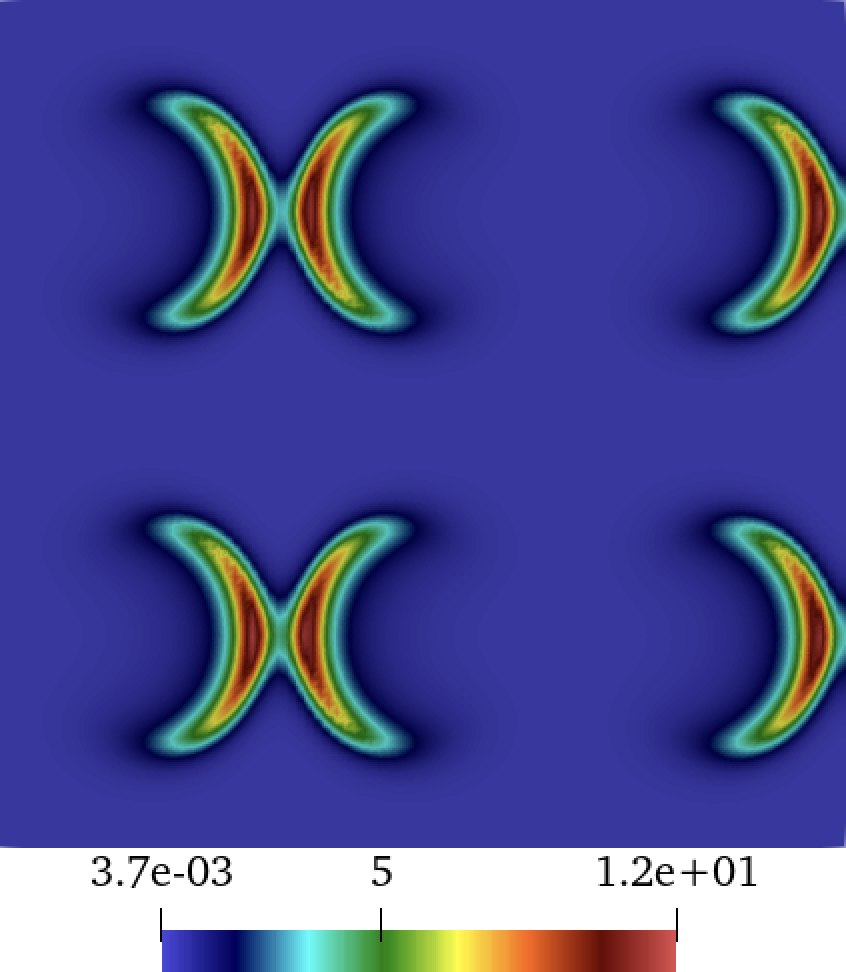}
  \put(23,102){\footnotesize{$t=777$}}
    \end{overpic} 
        \begin{overpic}[percent,width=0.16\textwidth]{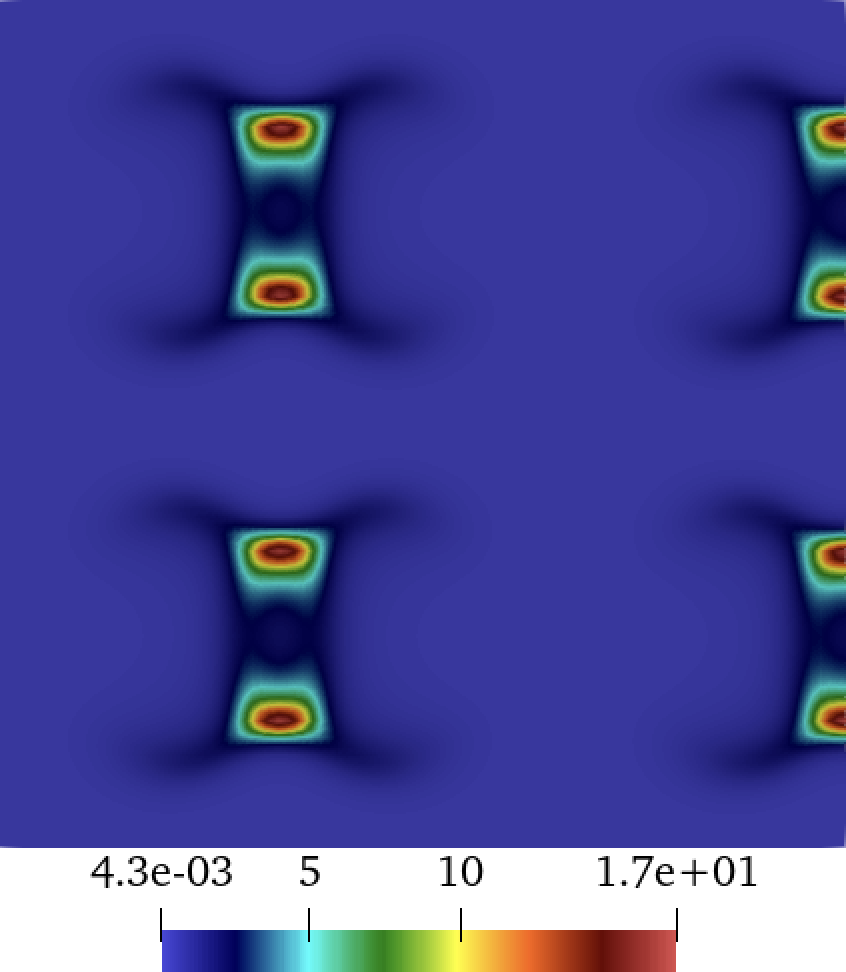}
       \put(23,102){\footnotesize{$t=778$}}
    \end{overpic} 
      \begin{overpic}[percent,width=0.16\textwidth]{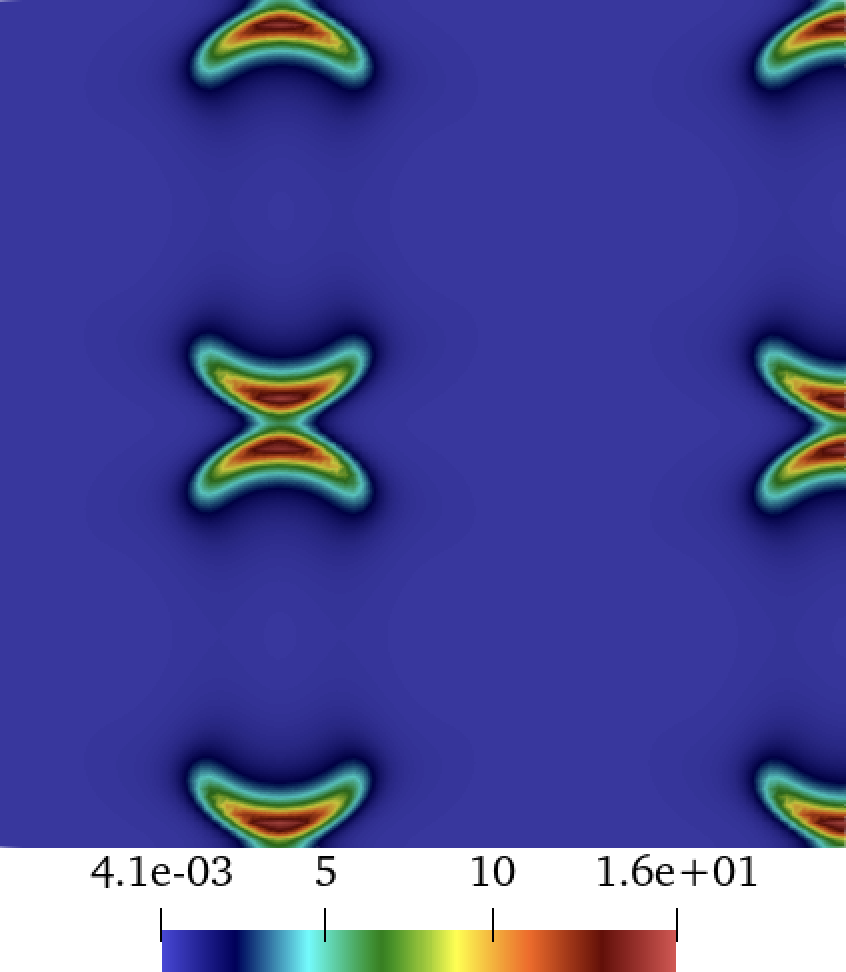}
       \put(23,102){\footnotesize{$t=780$}}
    \end{overpic} 
  \begin{overpic}[percent,width=0.16\textwidth]{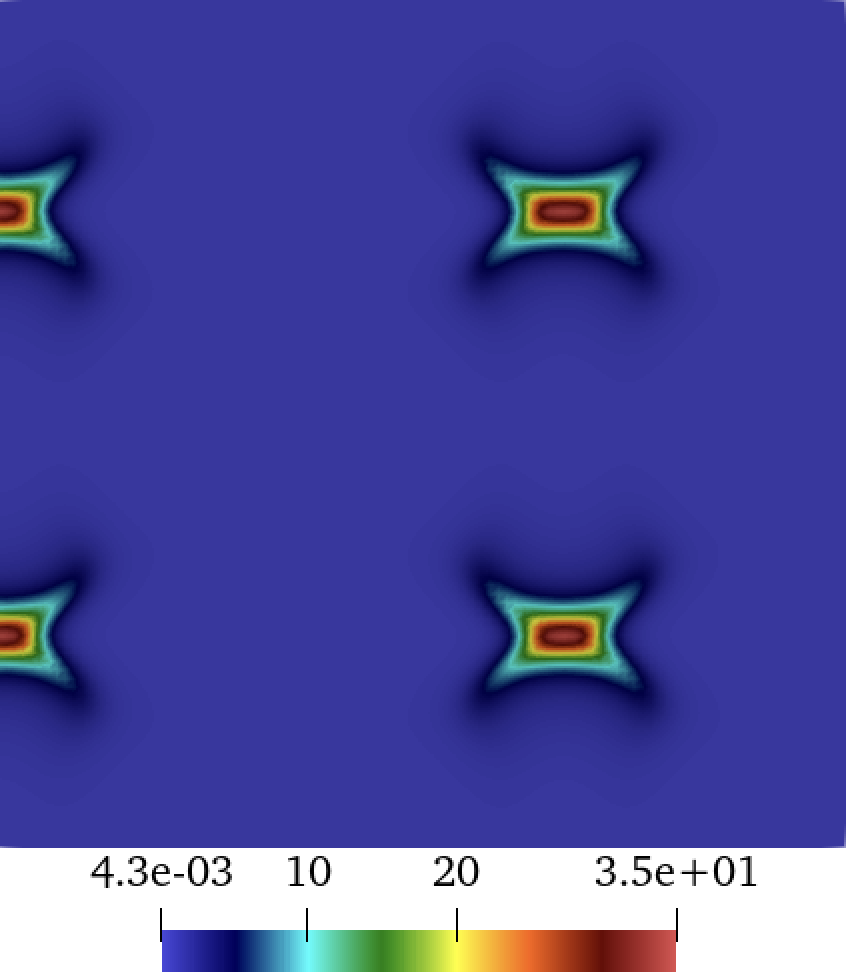}
      \put(23,102){\footnotesize{$t=788$}}
    \end{overpic} 

    \caption{Case 10: expected number of crimes per unit area $S$  from $t = 773$ to $t = 788$.}
    \label{fig:pi1_snap}
\end{figure}

For one final comparison of the results by the PDE model vs agent-based 
model,  
Fig.~\ref{fig:pi01_snapABM} shows
the evolution of $S$ in case 9 given by the agent-based solver after non-dimensionalization. Again, we set 
$\Gamma = 0.0285$, $\theta = 0.2339$, $\Sigma = 0$ and $\omega = 1/15$ to obtain $\Gamma \theta (1 - \Sigma)/\omega^2 \approx 1.5$.
Fig.~\ref{fig:pi01_snapABM} shows 
hotspots that are in good agreement
with those in Fig.~\ref{fig:pi01_snap}
in terms of number, size, and magnitude. 
From the results of cases 2, 8, and 9, we can conclude that the 
solutions given by the PDE solver and agent-based solver match well, especially when 
the densities of criminals 
and police are large (i.e., cases 2 and 8).

\begin{figure}[htb!]
     \centering
         \begin{overpic}[percent,width=0.16\textwidth]{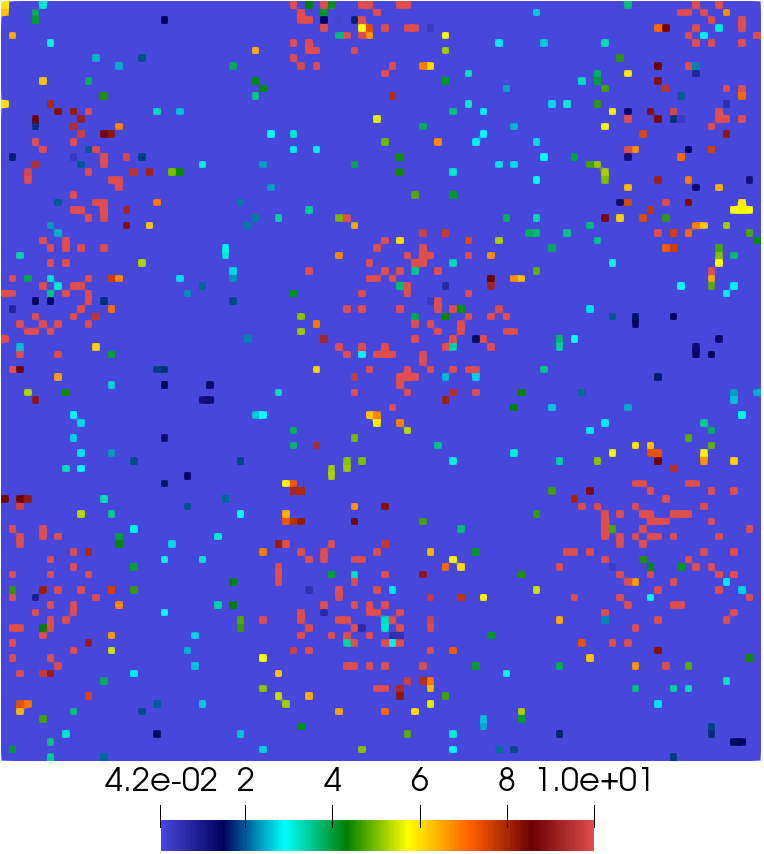}
     \put(23,102){\footnotesize{$t=744$}}
    \end{overpic} 
     \begin{overpic}[percent,width=0.16\textwidth, grid=false]{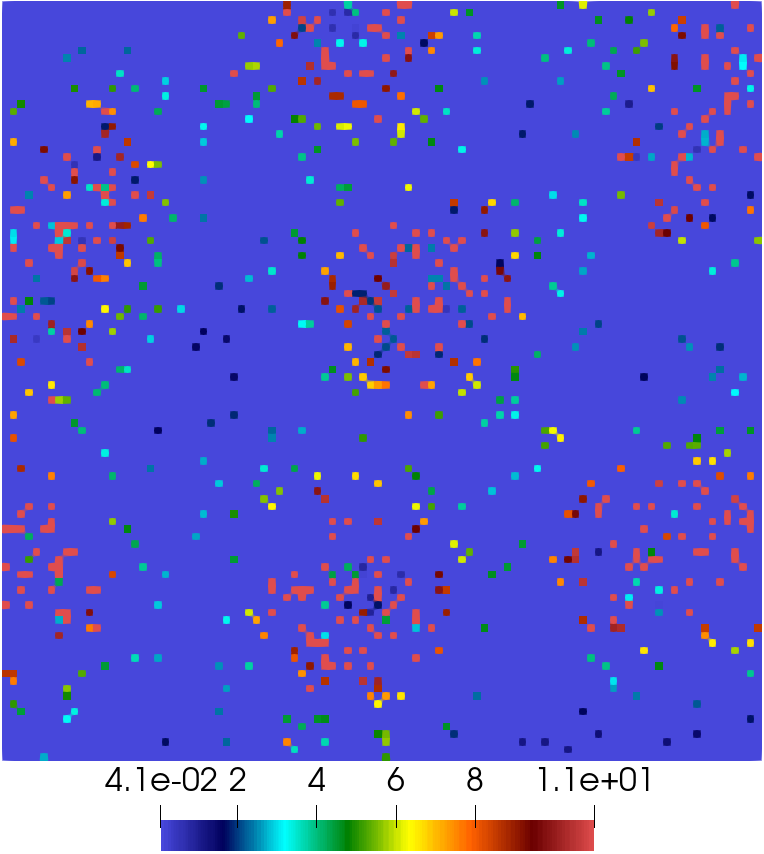}
     \put(23,102){\footnotesize{$t=747$}}
    \end{overpic}
    \begin{overpic}[percent,width=0.16\textwidth]{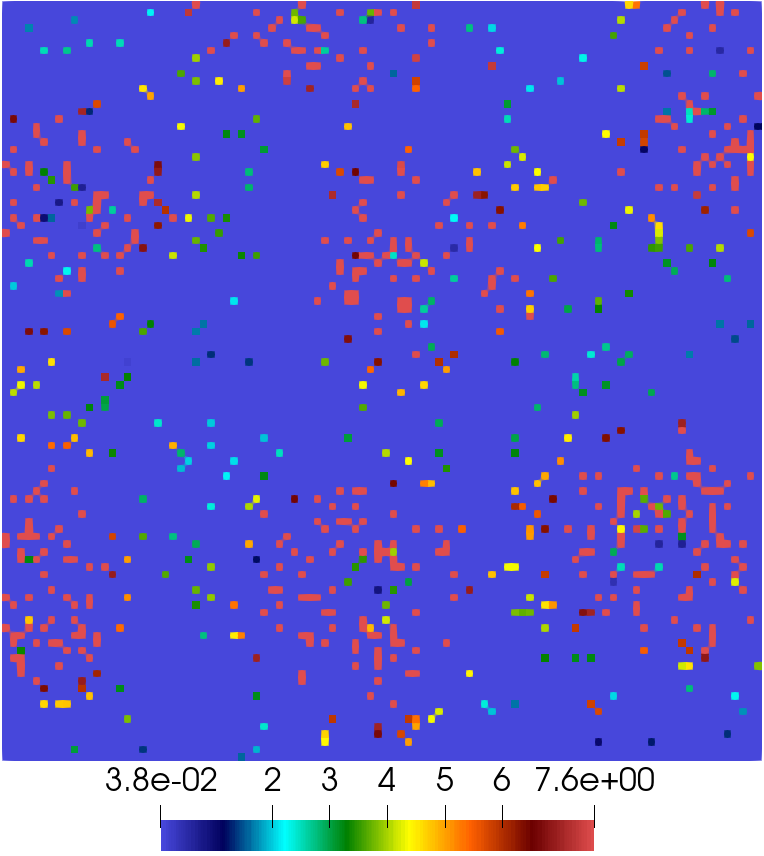}
  \put(23,102){\footnotesize{$t=750$}}
    \end{overpic} 
        \begin{overpic}[percent,width=0.16\textwidth]{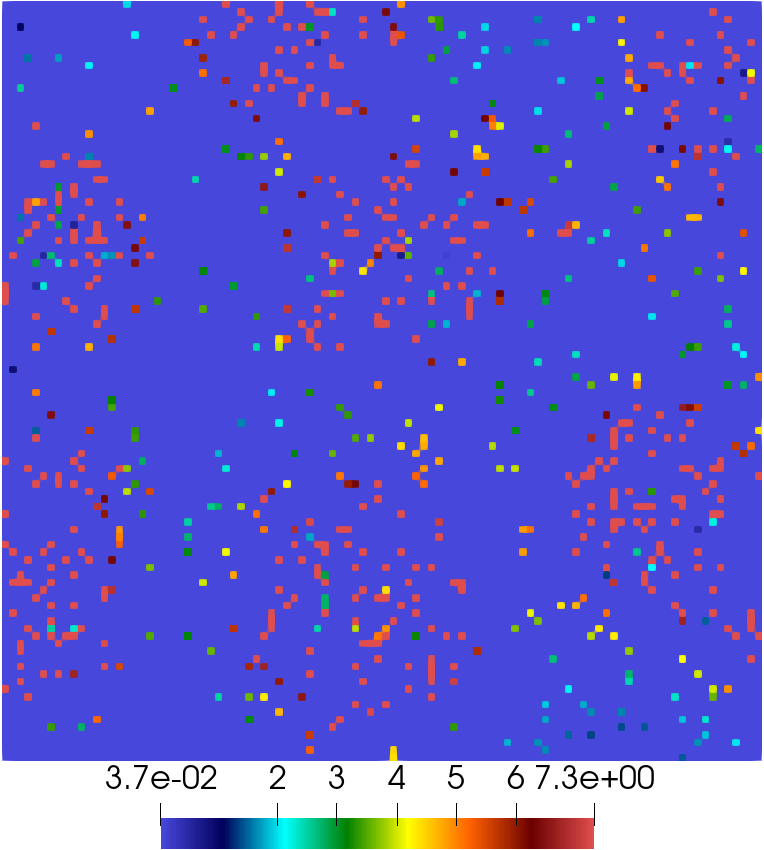}
       \put(23,102){\footnotesize{$t=753$}}
    \end{overpic} 
      \begin{overpic}[percent,width=0.16\textwidth]{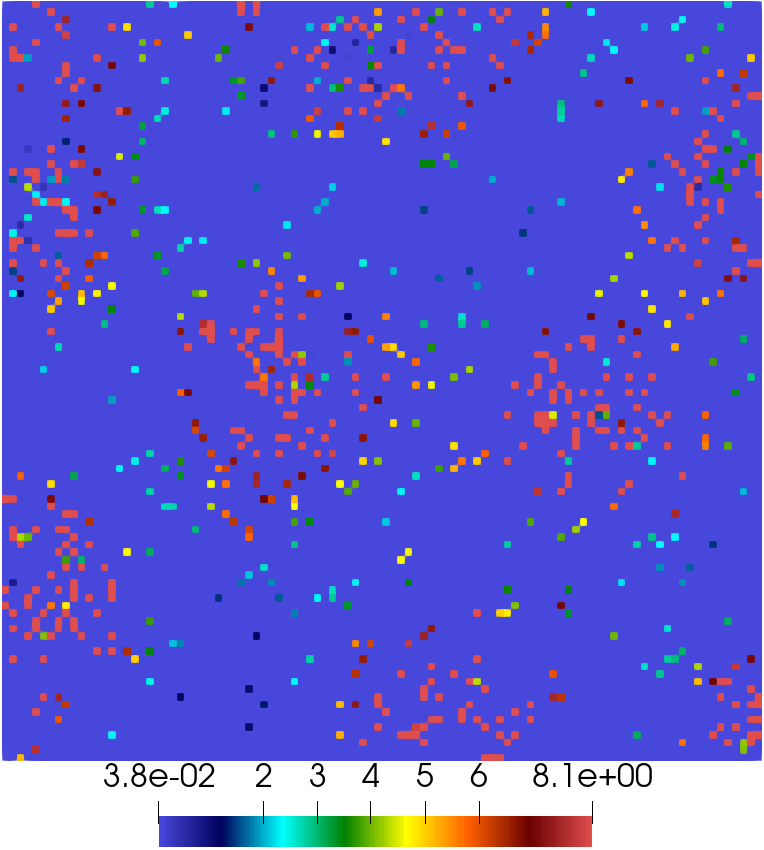}
       \put(23,102){\footnotesize{$t=756$}}
    \end{overpic} 
  \begin{overpic}[percent,width=0.16\textwidth]{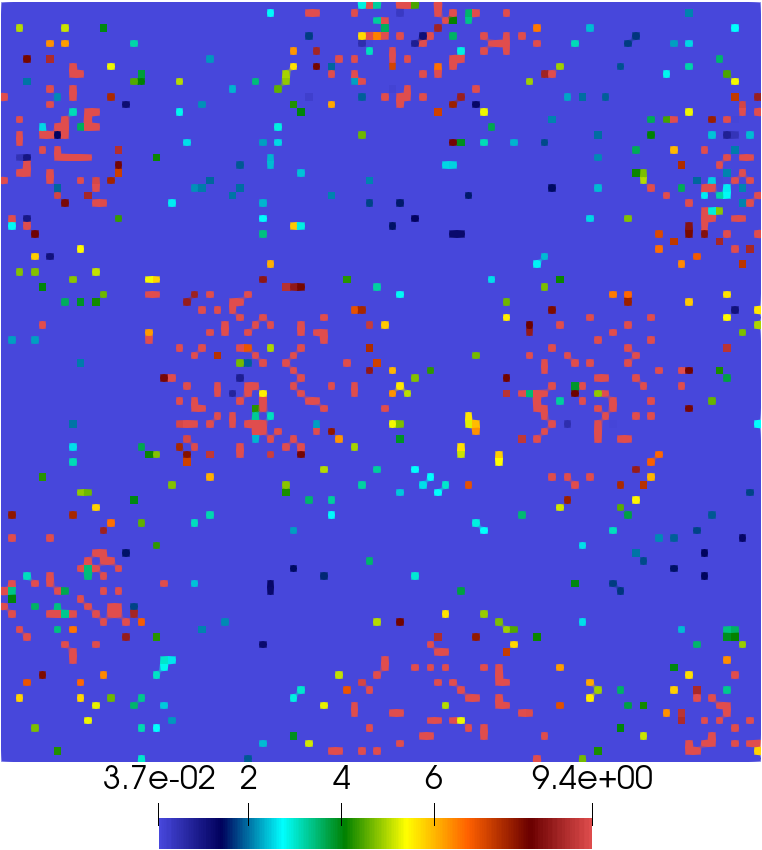}
      \put(23,102){\footnotesize{$t=759$}}
    \end{overpic} 

    \caption{Case 9 - agent-based model: expected number of crimes per unit area $S$ from $t = 744$ to $t = 759$.}
    \label{fig:pi01_snapABM}
\end{figure}


To summarize cases 2-10, we display 
return maps of $\langle S \rangle_{n+1}$ versus $\langle S \rangle_n$ in Fig.~\ref{fig:return_maps}.
Increasing $\eta$ (top row in Fig.~\ref{fig:return_maps}) regularizes the dynamics, replacing a scattered cloud at small $\eta$ with well-defined closed curves at larger values, indicating that larger $\eta$ promotes coherent oscillations and suppresses spatial heterogeneity. 
Increasing
$\Gamma \theta(1-\Sigma)/\omega^2$
(second row in Fig.~\ref{fig:return_maps}) drives the system from a fixed point to periodic and then to higher-dimensional attractors in the return map, consistent with chaotic dynamics. The third row in 
Fig.~\ref{fig:return_maps} confirms
the delay-like effect of $\tau$: small $\tau$ yields a stable equilibrium, intermediate values produce sustained oscillations, and large $\tau$ leads to irregular dynamics.
Finally, variations in $\pi_0$ (bottom row in 
Fig.~\ref{fig:return_maps}) produce a transition similar to variations in 
$\Gamma \theta(1-\Sigma)/\omega^2$:
from equilibrium to periodic and more complex attractors. 

In conclusion, the results in this section so far suggest that the PDE system~\eqref{eq:continous_1}--\eqref{eq:continous_4} supports a rich variety of dynamical regimes.

\begin{figure}[htb!]
    \centering
    \begin{minipage}{1.0\textwidth}
        \centering
        \includegraphics[width=\linewidth]{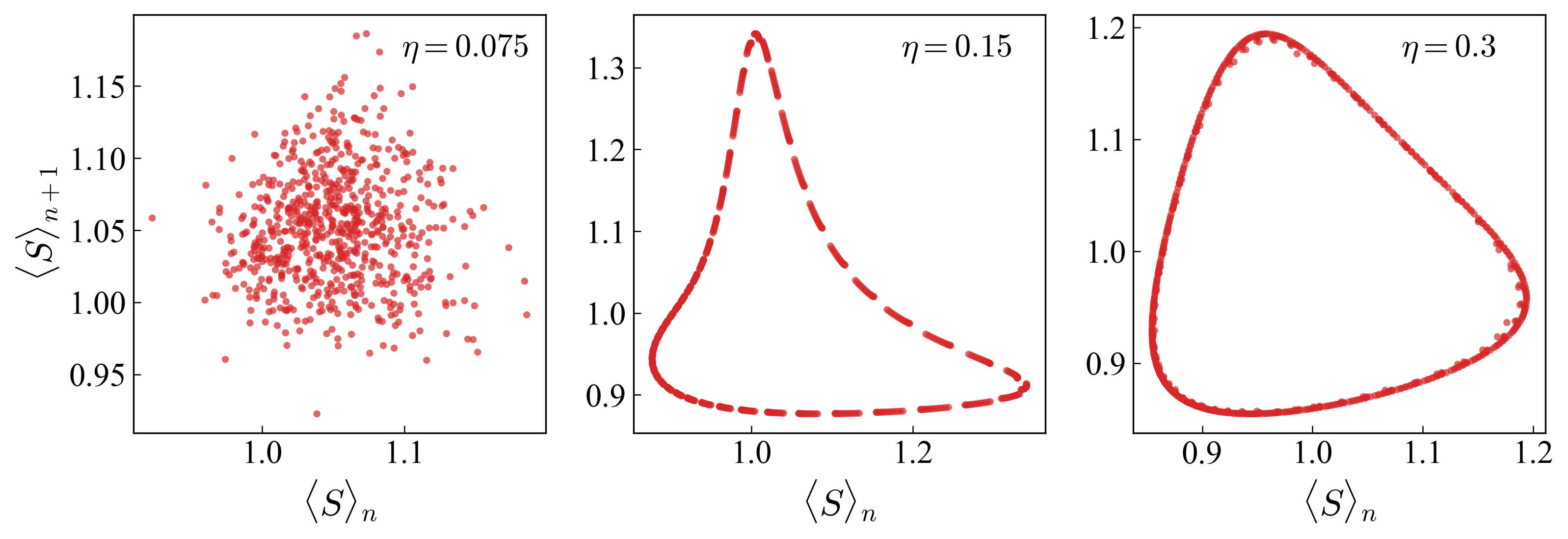}
    \end{minipage}
    \hfill
    \\
    \begin{minipage}{1.0\textwidth}
        \centering
        \includegraphics[width=\linewidth]{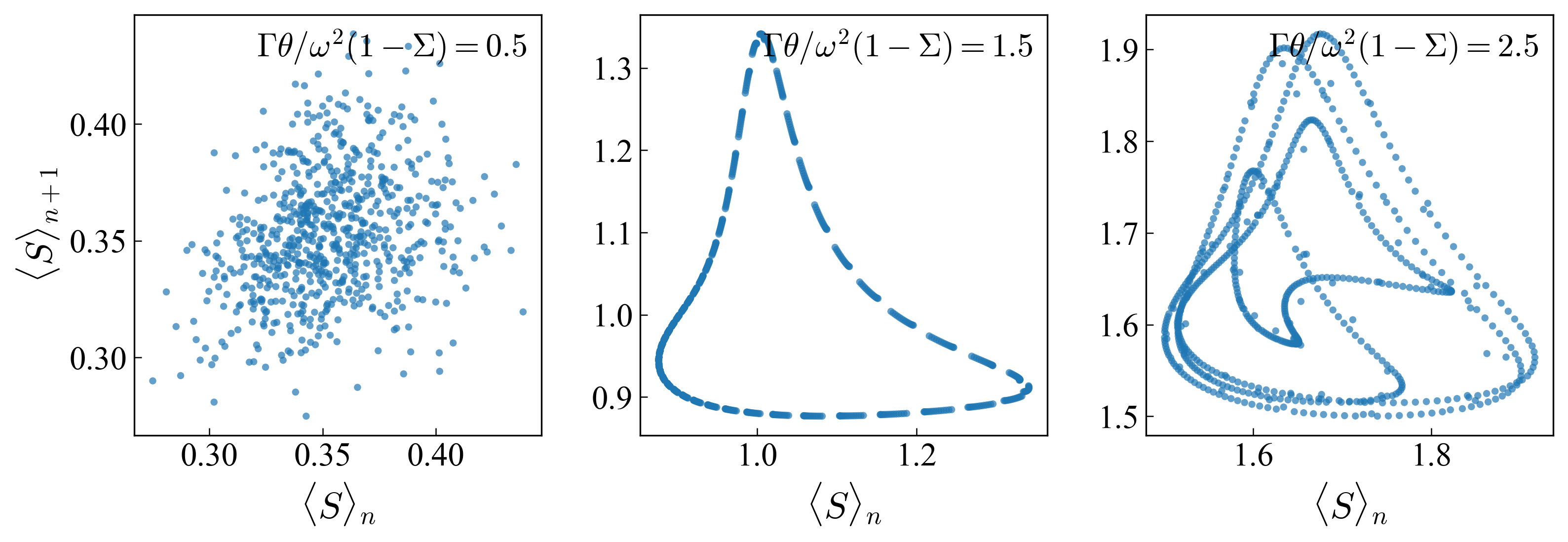}
    \end{minipage}
\\
        \begin{minipage}{1.0\textwidth}
        \centering
        \includegraphics[width=\linewidth]{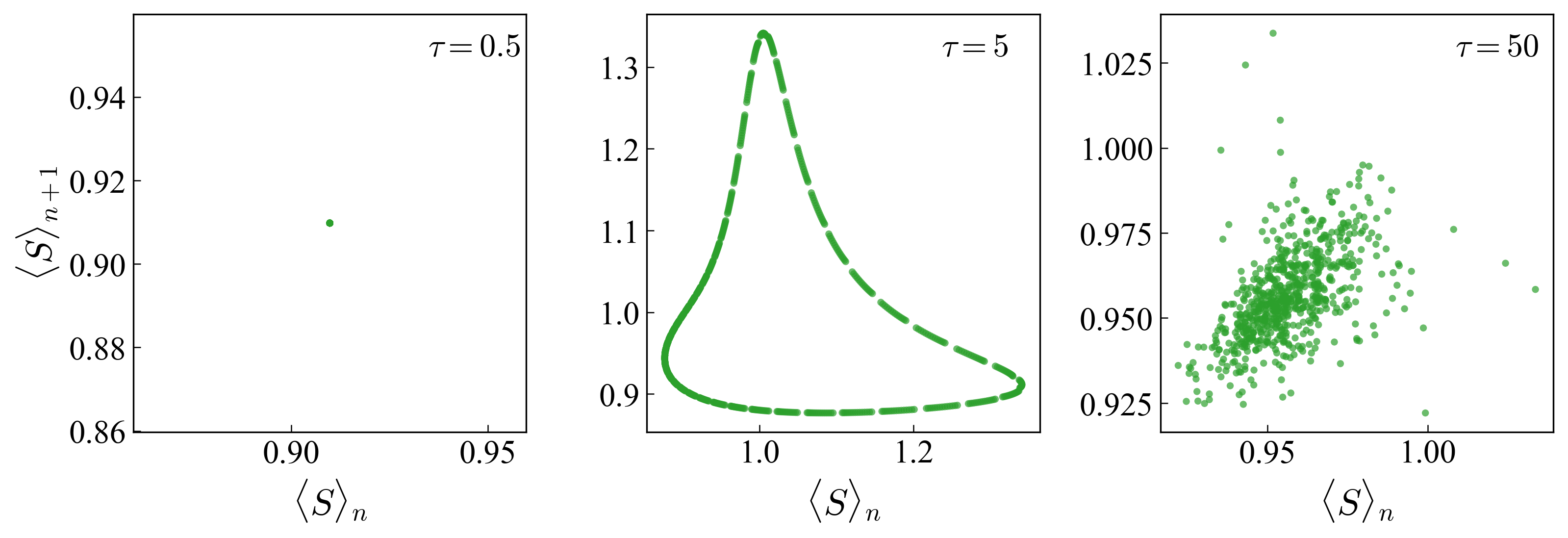}
    \end{minipage}
    \\
        \begin{minipage}{1.0\textwidth}
        \centering
        \includegraphics[width=\linewidth]{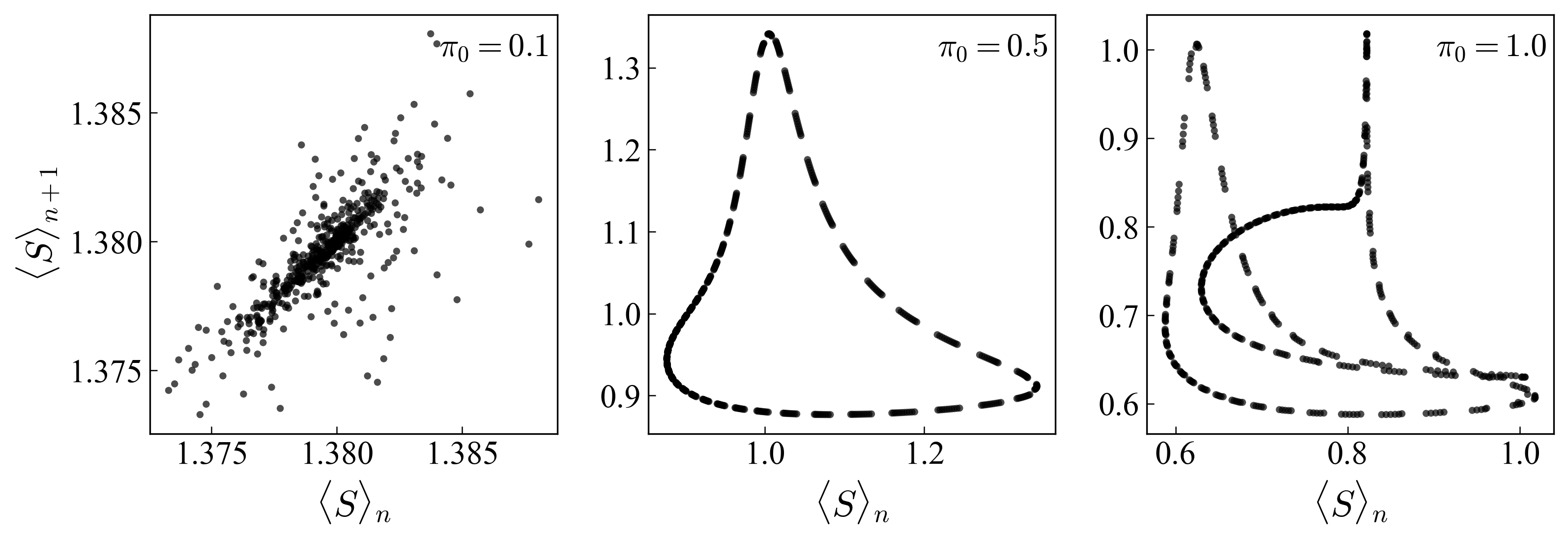}
    \end{minipage}
    \caption{Return maps of $\langle S \rangle_{n+1}$ versus $\langle S \rangle_n$
    for cases 2-10.}
    \label{fig:return_maps}
\end{figure}

\noindent {\bf Verification of one Hopf bifurcation}.
We verify the analytical predictions of Sec.~\ref{sec:lin_stab} 
through direct eigenvalue computation and nonlinear simulations when the parameters are as in case 2 and $\tau$ varies. 
For different values of $\tau$, 
Tab.~\ref{tab:hopf_verify} reports the real and imaginary parts of the dominant 
eigenvalue of~\eqref{eq:eig} at the critical mode $\mu^* = 1.28$, confirming a purely imaginary pair 
(up to numerical error)
at $\tau_c^* \approx 2.48$. 

\begin{table}[htb!]
\centering
\renewcommand{\arraystretch}{1.2}
\begin{tabular}{rcccccc}
\toprule
$\tau$ 
  & $0.5$ 
  & $2.48$ 
  & $2.5$ 
  & $3.0$ 
  & $5.0$ 
  & $10.0$ \\
\midrule
$\mathrm{Re}(\lambda_{\max})$
  & $-0.897$
  & $\approx 0$
  & $+0.003$
  & $+0.040$
  & $+0.114$
  & $+0.266$ \\
$\mathrm{Im}(\lambda_{\max})$
  & $\pm\,0.376$
  & $\pm\,0.452$
  & $\pm\,0.451$
  & $\pm\,0.408$
  & $\pm\,0.296$
  & $\pm\,0.157$ \\
\bottomrule
\end{tabular}
\caption{Real and imaginary parts of the dominant 
eigenvalue of~\eqref{eq:eig} at the critical mode $\mu^* = 1.28$ for different values of $\tau$.}
\label{tab:hopf_verify}
\end{table}


Fig.~\ref{fig:hopf_validation} (left) shows the fundamental 
oscillation frequency $f_{\rm comp}$ extracted from nonlinear PDE simulations as $\tau$ varies, together with the value $f_{\rm lin} = \omega_0/(2\pi) \approx 0.0717$ from the analysis in Sec.~\ref{sec:lin_stab}.
The fundamental 
oscillation frequency computed from the simulations is obtained by dividing the dominant FFT peak of 
$\langle S \rangle$ by four, since each fundamental 
cycle contains four sub-phases with similar spatial 
averages (see Fig.~\ref{fig:case2_snap}).
The distance between $f_{\rm comp}$ and $f_{\rm lin}$ 
increases with $\tau$, showing that 
the further $\tau$ is from $\tau^*_c$, the less accurate the prediction from the linear stability analysis becomes
due to the fact that the problem is highly nonlinear.
Fig.~\ref{fig:hopf_validation} (left)
plots the square of the amplitude of the oscillations in 
$\langle S \rangle$ 
vs $\tau - \tau_c^*$, with linear regression
for the points with $\tau -\tau^*_c \le 5$.
{We see that the squared oscillation amplitude scales linearly with $\tau - \tau_c^*$ near criticality ($R^2 = 0.98$), in agreement with the $\mathcal{O}(\sqrt{\tau - \tau_c^*})$ amplitude scaling predicted by the Hopf normal form. Moreover, for $\tau$ slightly above $\tau_c^*$ the system settles onto a small-amplitude limit cycle whose amplitude grows continuously from zero, with no finite-time blow-up observed as $\tau$ crosses $\tau_c^*$. These observations are consistent with a supercritical Hopf bifurcation.
A rigorous proof of the nature of the Hopf bifurcation
is beyond the scope of this work. Deviations from the linear fit at larger $\tau - \tau_c^*$ reflect higher-order effects outside the weakly nonlinear regime.}

\begin{figure}[htb!]
\centering
\includegraphics[width=\textwidth]{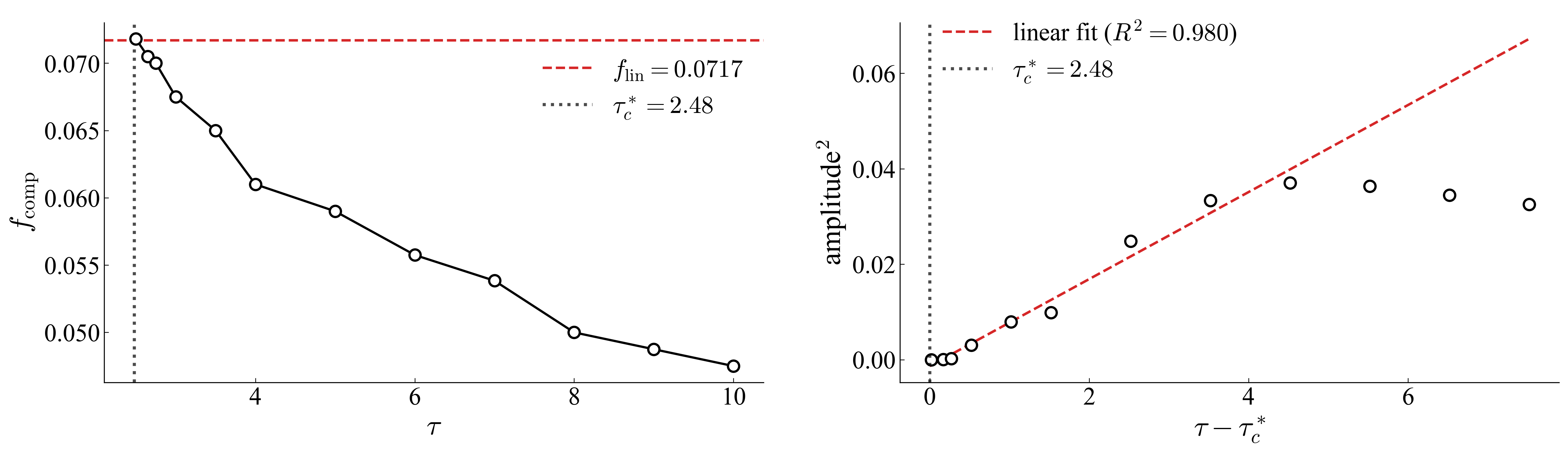}
\caption{Left: 
fundamental 
oscillation frequency $f_{\rm comp}$ extracted from nonlinear PDE simulations vs $\tau$ (black curve), together with the value $f_{\rm lin}$ (red dashed line) from the analysis in Sec.~\ref{sec:lin_stab}. Right: square of the amplitude of the oscillations in 
$\langle S \rangle$ 
vs $\tau - \tau_c^*$ (black curve), together with linear regression (red dashed line) for low $\tau - \tau_c^*$ values.
}
\label{fig:hopf_validation}
\end{figure}

Finally, 
Fig.~\ref{fig:bifurcation} shows 
the real part of the dominant eigenvalue, i.e., $\max_{\mu_j}\mathrm{Re}(\lambda)$, as $\tau$ is varied over a larger interval than in Fig.~\ref{fig:hopf_validation} (left). We clearly see that
the horizontal axis is crossed at $\tau_c^*\approx 2.48$ (marked with a red dashed line) and we are able to approximate the  
transversality $\partial_\tau\mathrm{Re}(\lambda)|_{\tau_c^*} 
\approx 0.091 > 0$. 

\begin{figure}[htb!]
\centering
\includegraphics[width=.9\linewidth]{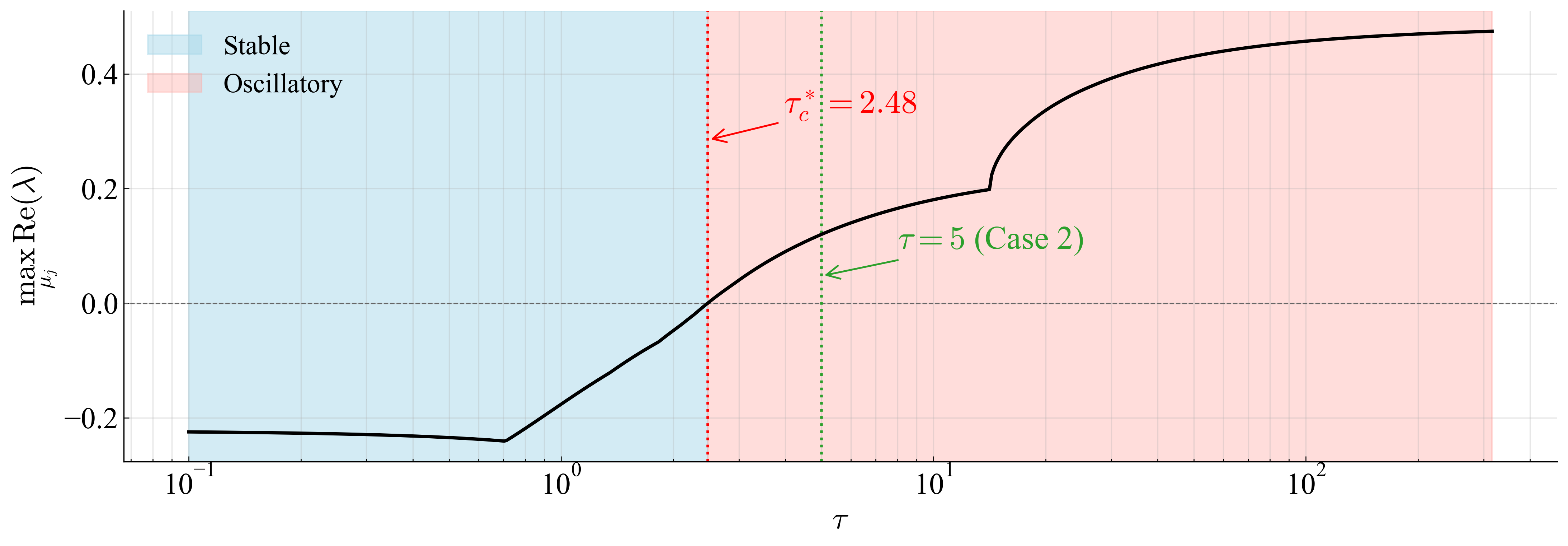}
\caption{
Real part of the dominant eigenvalue vs $\tau$, with highlighted the value of $\tau^*_c$ (red dotted line) and the value of $\tau$ in case 2 (green dotted line). The stable part of the plane is colored in blue, while the oscillatory part is colored in pink.}
\label{fig:bifurcation}
\end{figure}

\subsection{Comparison of policing strategies}\label{sec:comparison}

The goal of this section is to compare the model proposed in this paper with
two models available in the literature:
\begin{itemize}
    \item[-] {\emph{Fixed policing}}: 
    the criminal activity is allowed to evolve with no policing till
    $t = t_s$, at which time 
    police density is 
    computed from the current attractiveness and 
    held fixed thereafter~\cite{short2010nonlinear}. See Remark \ref{rem:pol_strategies}.
    \item[-] {\emph{Optimal policing}}:   at each time step, the police field minimizes the expected number of crimes at the given time
    subject to the constraint that the total police population is kept constant~\cite{zipkin2014}. See Remark \ref{rem:pol_strategies}.
\end{itemize}
Neither of these two strategies is realistic. We will call our model with delay $\tau = 5$ \emph{realistic policing.}

To compare the three policing strategies, we set $\eta = 0.15$ and the other parameters as in Case~3. We first run the simulation with no police, i.e., using the model from 
\cite{short2008statistical} till 
$t_s = 200$. At this time, the crime hotspots are fully developed. 
Fig.~\ref{fig:strategy_comparison}
reports $\langle S \rangle$ over time (black curve) and
shows the hotspots at $t_s = 200$. At $t_s = 200$, police are deployed according to the three strategies mentioned above, assuming 
$M = \int_\Omega \pi\, d\mathbf{x} = 50$. In the case of fixed policing, 
we set $\mu = 5$ and $A_c = 1.5$. We slightly modify \eqref{eq:fixed_pol} as follows
\begin{align*}
    \pi (\bm{x}) = -\log\left( \max \left\{ \frac{1}{2} \left( 1 - \tanh{\mu \left( A(\bm{x}, t_s) - A_c\right) } \right),10^{-12} \right\} \right), 
\end{align*}
However, there is no guarantee that this $\pi$ would be such that $\int_\Omega \pi\, d\mathbf{x} = 50$. Thus, if $ \int_\Omega \pi\, d\mathbf{x} \neq 50$, we rescale it as: 
\begin{equation*}
    \pi(\bm{x}) \leftarrow \pi(\bm{x}) \frac{M}{\int_\Omega \pi(\bm{x}) d\bm{x}}.
\end{equation*}






Fig.~\ref{fig:strategy_comparison} shows $\langle S(t) \rangle$ for each strategy. Optimal policing achieves the lowest crime level ($\langle S \rangle \approx 0.92$) 
and eliminates all spatial heterogeneity, producing a uniform crime 
landscape (also reported in Fig.~\ref{fig:strategy_comparison}). However, this policing strategy requires solving a global 
optimization problem with complete real-time crime data at every 
time step. 
Fixed policing reduces crime to $\langle S \rangle \approx 1.2$ 
but retains displaced hotspots, also shown in Fig.~\ref{fig:strategy_comparison}, in regions not targeted by the frozen-in-time police patrolling. 
Our realistic policing model produces an 
intermediate time-averaged crime level ($\langle S \rangle \approx 1.03$) 
with sustained oscillations consistent with the linear stability analysis in Sec.~\ref{sec:lin_stab}. In this scenario, the moving hotspots 
(shown in Fig.~\ref{fig:strategy_comparison}) reflect the delayed feedback between crime 
and police: police respond to lagged crime data, suppress the 
current hotspot, but arrive too late to prevent its re-emergence 
elsewhere.

\begin{figure}[htb!]
\centering
\includegraphics[width=0.97\linewidth]{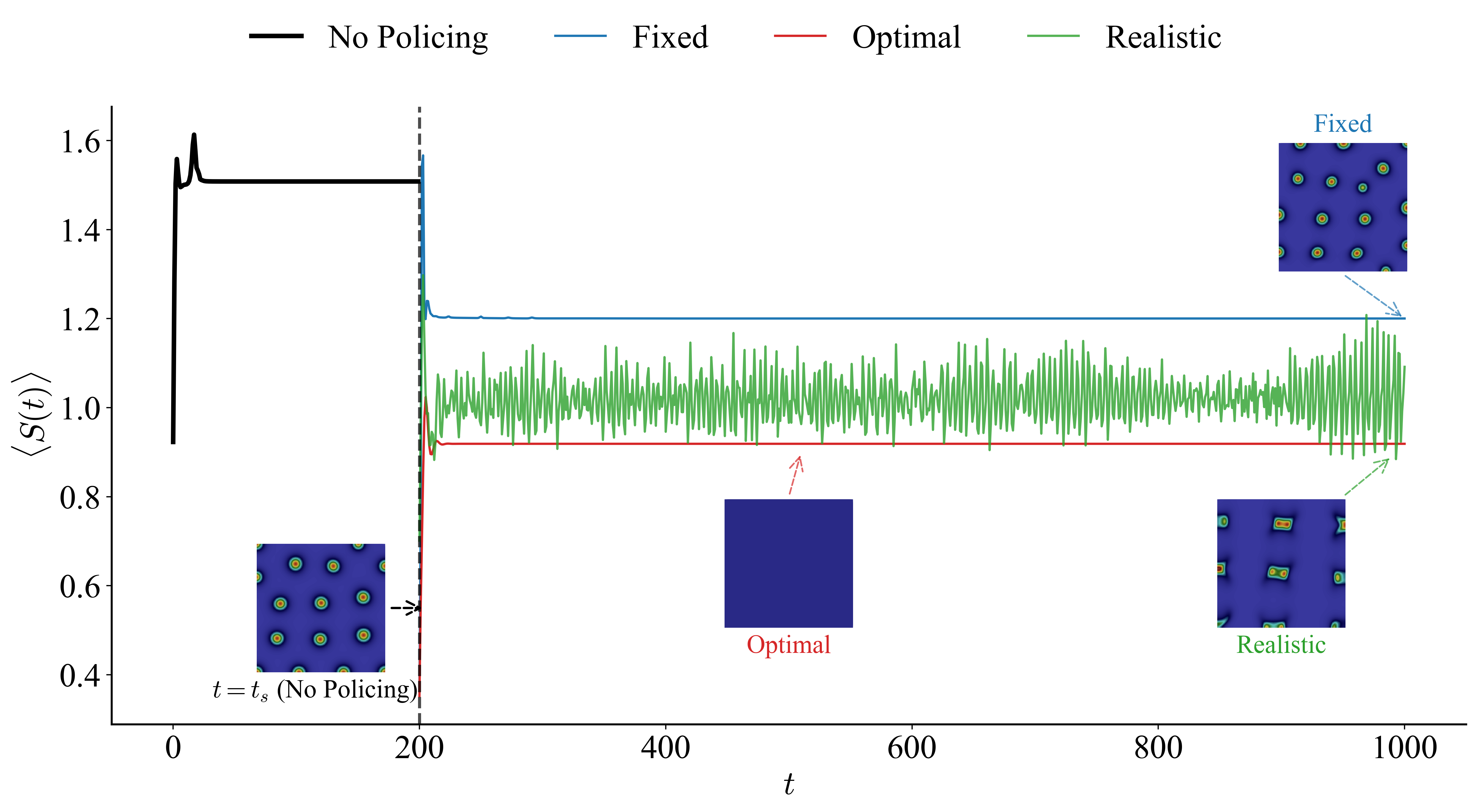}
\caption{Spatially averaged crime rate $\langle S\rangle$ for three policing strategies starting at 
$t_s = 200$. For $t < 200$, we let crime hotspots fully form in absence of law enforcement. The hotspots in $S$ at $t_s = 200$ 
are shown, as well as representative $S$
for each policing strategy after 
$t_s = 200$. 
All spatial plots use a common color scale ranging from 0 (blue) to 25 (red).}
\label{fig:strategy_comparison}
\end{figure}

These results illustrate how the degree of information availability 
shapes not only the level but also the spatial and temporal structure 
of crime. Optimal policing, while theoretically superior, assumes 
instantaneous global coordination that is operationally unrealistic. 
Fixed policing requires no real-time data but cannot adapt to crime 
displacement. Our model occupies the realistic middle ground: police 
adapt continuously but with a finite delay $\tau$, producing dynamics 
that are qualitatively distinct from both extremes.


\section{Conclusions}\label{sec:concl}

We extended a well known agent-based model for residential burglary
to include police deployment through a delayed feedback mechanism. 
The model assumes that if police are not busy at a certain location, they move according to 
a biased random walk towards areas with higher reported crime data. Such crime data is assumed to become available and used for decision-making with some delay. 
On the basis of the agent-based model, we derived a continuum PDE system and showed analytically that time delays in police response can destabilize the homogeneous equilibrium through a Hopf bifurcation. Linear stability analysis via the Routh-Hurwitz criterion identified a critical delay threshold for given model parameters, and numerical eigenvalue computations provided evidence for the supercritical nature of the bifurcation.

Simulations of both the PDE and agent-based models corroborate the theoretical predictions and reveal a rich range of spatio-temporal dynamics, including persistent oscillations, splitting and merging of hotspots, and multi‑frequency regimes. A parametric study showed that neighborhood effects, police density, and delay jointly control stability, hotspot size, and oscillatory behavior, with intermediate delays and police densities being most conducive to sustained cycles. 
Our results indicated that timely access to crime data is crucial for stabilizing crime levels. Small delays allow police to respond effectively and suppress large-scale temporal fluctuations, while intermediate delays promote strong oscillations due to delayed feedback. Very large delays weaken the feedback mechanism, leading to less structured and potentially less predictable crime dynamics.
Our results suggested that increasing police density can qualitatively change the system dynamics. Intermediate enforcement levels may produce predictable crime cycles, whereas very strong enforcement suppresses crime more effectively but generates irregular temporal patterns rather than steady suppression.

A comparison with two police strategies in the mathematical literature, called fixed and optimal policing, highlighted that our PDE model occupies a realistic middle ground: it outperforms fixed policing (i.e., frozen-in-time police patrolling) while avoiding the unrealistic assumptions of instantaneous global optimization of the police force.

\section*{Acknowledgments}
We thank Jacobe West for searching spatio-temporal patrolling and crime data that led to Figs.~\ref{fig:empirical_correspondence_q3_2025} and 
\ref{fig:empirical_correspondence_q4_2025}. Zhong is supported by NSF-AoF grant $\#2225507$. Zhong also gratefully acknowledges funding provided by the Oak Ridge Associated Universities Ralph E. Powe Junior Faculty Enhancement Award for FY$2024$.

\section*{Appendix A: Data sources and processing for Fig.~\ref{fig:empirical_correspondence_q3_2025}-\ref{fig:empirical_correspondence_q4_2025}}

Reported burglary incidents are
taken from a publicly available
dataset of the City of Chicago \cite{Chicago_crime}
and filtered to primary type ``Burglary''. Investigatory Stop Reports
(ISRs) are obtained from 
the Chicago Police Department (CPD) public releases \cite{ISR_Chicago} and filtered to
non-vehicle ``investigatory stop'' with ``proximity to crime'' or matching
an ``offender description''. CPD sworn staffing is taken from the Office of Inspector General
for the City of Chicago
\cite{OIG_Chicago}
and restricted to geographic CPD beats. 
All beat-level quantities use the
official CPD beat polygons from the City of Chicago boundary data \cite{Chicago_beats}.

The burglary and ISR hotspot panels in Figs.~\ref{fig:empirical_correspondence_q3_2025} and 
\ref{fig:empirical_correspondence_q4_2025} are constructed from point-level records
using kernel density estimation on a common map extent, grid resolution, and
smoothing bandwidth. For a given quarter, burglary and ISR counts are divided by three
to obtain average monthly event densities. 
For police staffing, we average the
monthly staffing counts within the quarter. 
Dividing by beat area
gives police density in officers per \(\mathrm{km}^2\).

\section*{Appendix B: Explicit form of the 
coefficients in the characteristic polynomial \eqref{eq:char_poly}
}\label{app:coefficients}
 
This appendix provides the derivation of the
coefficients $a_0, a_1, a_2, a_3$ in the characteristic
polynomial~\eqref{eq:char_poly}. 
For compactness of presentation, let us set
\begin{equation}\label{eq:shorthand}
\alpha = \bar{A}\,e^{-\bar{\pi}},
\qquad
\zeta  = \bar{\rho}\,e^{-\bar{\pi}},
\end{equation}
so that $\bar{H} = \bar{\rho}\,\alpha = \zeta\,\bar{A}$. 
%

We find $P(\lambda)
=\det(\lambda I - J)$, with $J$ given in~\eqref{eq:eig}, 
by cofactor expansion along the third row:
\begin{equation}\label{eq:factorization}
P(\lambda)
= (\lambda + \mu_j)\!\left(\lambda + \tfrac{1}{\tau}\right)
  \!\big[\lambda^2 + b_1\lambda + b_0\big]
+ \frac{2\mu_j\bar{\pi}}{\tau}\,
  (\lambda + \eta\mu_j + 1)(\lambda + \mu_j + \alpha),
\end{equation}
where
\begin{align}
b_1 &= (1+\eta)\mu_j + 1 + \alpha - \zeta,
      \label{eq:b1_app}\\
b_0 &= (\eta\mu_j + 1)(\mu_j + \alpha) - 3\mu_j\zeta.
      \label{eq:b0_app}
\end{align}
The quadratic factor $\lambda^2 + b_1\lambda + b_0$ depends
only on $\eta$, $\mu_j$, and the equilibrium
values in~\eqref{eq:hom_eq_sol}, which
do not depend on $\tau$.

 
Expanding~\eqref{eq:factorization} and collecting powers of
$\lambda$ 
yields
\begin{align}
a_3 &= (2+\eta)\mu_j + 1
       + \bar{A}e^{-\bar{\pi}} - \bar{\rho}e^{-\bar{\pi}}
       + \frac{1}{\tau},
      \label{eq:a3_app}\\
a_2 &= (1+2\eta)\mu_j^2 + 2\mu_j
       + (1+\eta)\mu_j\,\bar{A}e^{-\bar{\pi}}
       + \bar{A}e^{-\bar{\pi}}
       - 4\mu_j\,\bar{\rho}e^{-\bar{\pi}}
       \notag\\
    &\quad
     + \frac{1}{\tau}\Bigl[
         (2+\eta)\mu_j + 1
         + \bar{A}e^{-\bar{\pi}} - \bar{\rho}e^{-\bar{\pi}}
         + 2\mu_j\bar{\pi}
       \Bigr],
      \label{eq:a2_app}\\
a_1 &= \mu_j\,(\eta\mu_j + 1)
         \bigl(\mu_j + \bar{A}e^{-\bar{\pi}}\bigr)
       - 3\mu_j^2\,\bar{\rho}e^{-\bar{\pi}}
       \notag\\
    &\quad
     + \frac{1}{\tau}\Bigl[
         (1+2\eta)\mu_j^2 + 2\mu_j
         + (1+\eta)\mu_j\,\bar{A}e^{-\bar{\pi}}
         + \bar{A}e^{-\bar{\pi}}
         - 4\mu_j\,\bar{\rho}e^{-\bar{\pi}}
         \notag\\
    &\qquad\qquad
         + 2\mu_j\bar{\pi}\bigl(
             (1+\eta)\mu_j + 1 + \bar{A}e^{-\bar{\pi}}
           \bigr)
       \Bigr],
      \label{eq:a1_app}\\
a_0 &= \frac{\mu_j}{\tau}
       \Bigl[
         (1+2\bar{\pi})
         (\eta\mu_j + 1)
         \bigl(\mu_j + \bar{A}e^{-\bar{\pi}}\bigr)
         - 3\mu_j\,\bar{\rho}e^{-\bar{\pi}}
       \Bigr].
      \label{eq:a0_app}
\end{align}
 
%

We note that all four coefficients 
\eqref{eq:a3_app}-\eqref{eq:a0_app}
have a $\tau$-independent part plus a $1/\tau$ contribution. 
Thus, one way to simplify the coefficients is to consider the
$\tau \to \infty$ limit and have the $1/\tau$ terms vanish. 
In this limit, the coefficient $a_0 \propto 1/\tau \to  0$, so $P(\lambda)$ acquires a zero root. This corresponds to $H$ becoming frozen: it
retains its initial value indefinitely and ceases to respond to changes in crime intensity. 
The remaining coefficients $a_1, a_2, a_3$ reduce to the coefficients of the
$(A,\rho,\pi)$ subsystem.
The police effectively receives no information, and no periodic orbit can bifurcate.

To verify this numerically, we
set $\tau = 10^{12}$ and all the other parameters as in case 8. Fig.~\ref{fig:tau-1e12} (left) shows the evolution of $\langle S\rangle$, which 
rapidly converges to a constant value as confirmed by the power spectrum in 
Fig.~\ref{fig:tau-1e12} (right).
Notice that $\langle S\rangle$ is just slight higher for $\tau = 10^{12}$ than for
$\tau = 0.5$ (see Fig.~\ref{fig:tau}), However, there is a big difference: $S$ is homogeneously low in space for $\tau = 0.5$, while it shows steady-state hotspots that police cannot control for $\tau = 10^{12}$ (see Fig.~\ref{fig:tau1e12_snap}). 
Very large delays in the use of crime data for patrolling weaken the feedback mechanism, leading to persistent hotspots that cannot get eradicated.

\begin{figure}[htb!]
    \centering
    \includegraphics[width=.9\linewidth]{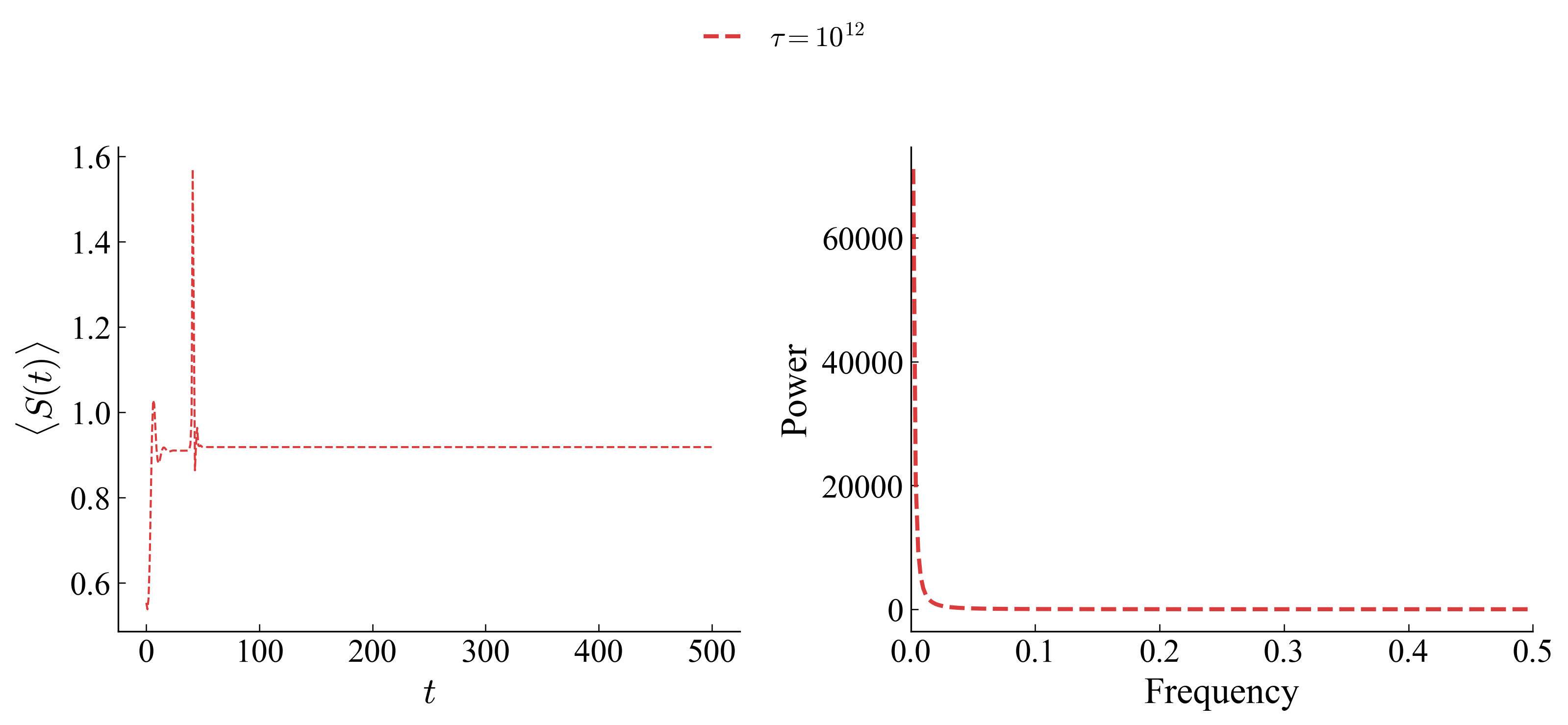}
    \caption{Time evolutions of spatially averaged quantities $\langle S \rangle$ for $\tau = 10^{12}$ (left) and associated power spectrum (right).}
    \label{fig:tau-1e12}
\end{figure}

\begin{figure}[htb!]
     \centering
         \begin{overpic}[percent,width=0.16\textwidth]{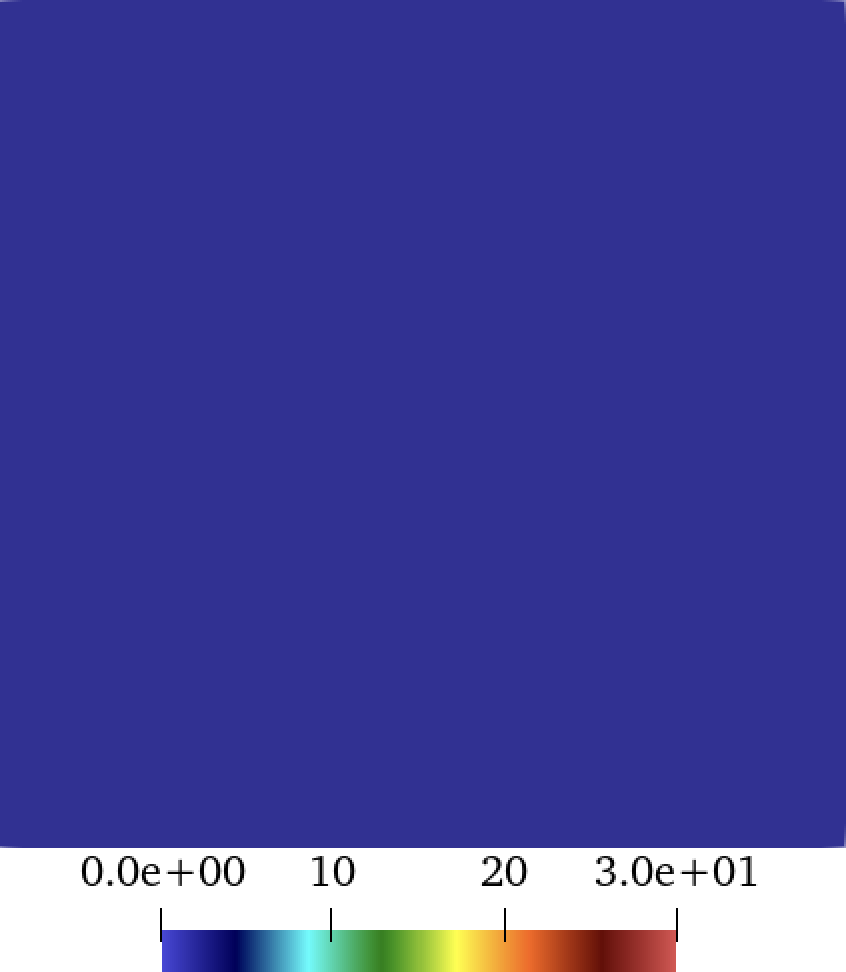}
     \put(23,102){\footnotesize{$t=0$}}
    \end{overpic} 
     \begin{overpic}[percent,width=0.16\textwidth, grid=false]{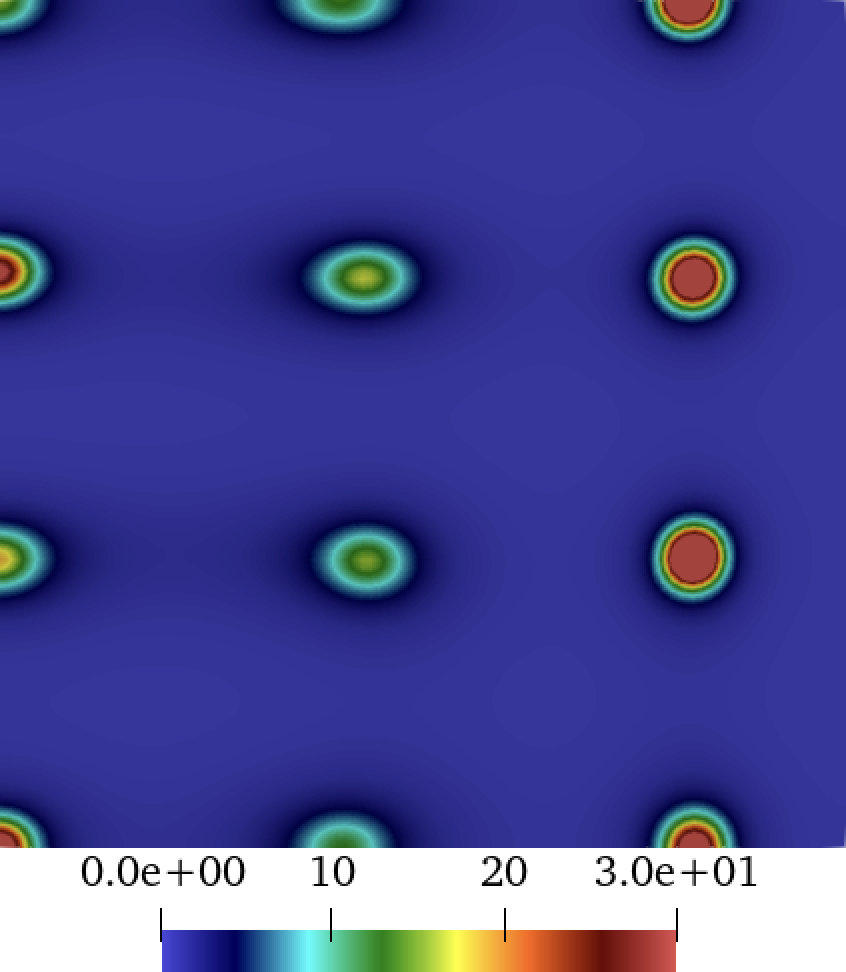}
     \put(23,102){\footnotesize{$t=40$}}
    \end{overpic}
    \begin{overpic}[percent,width=0.16\textwidth]{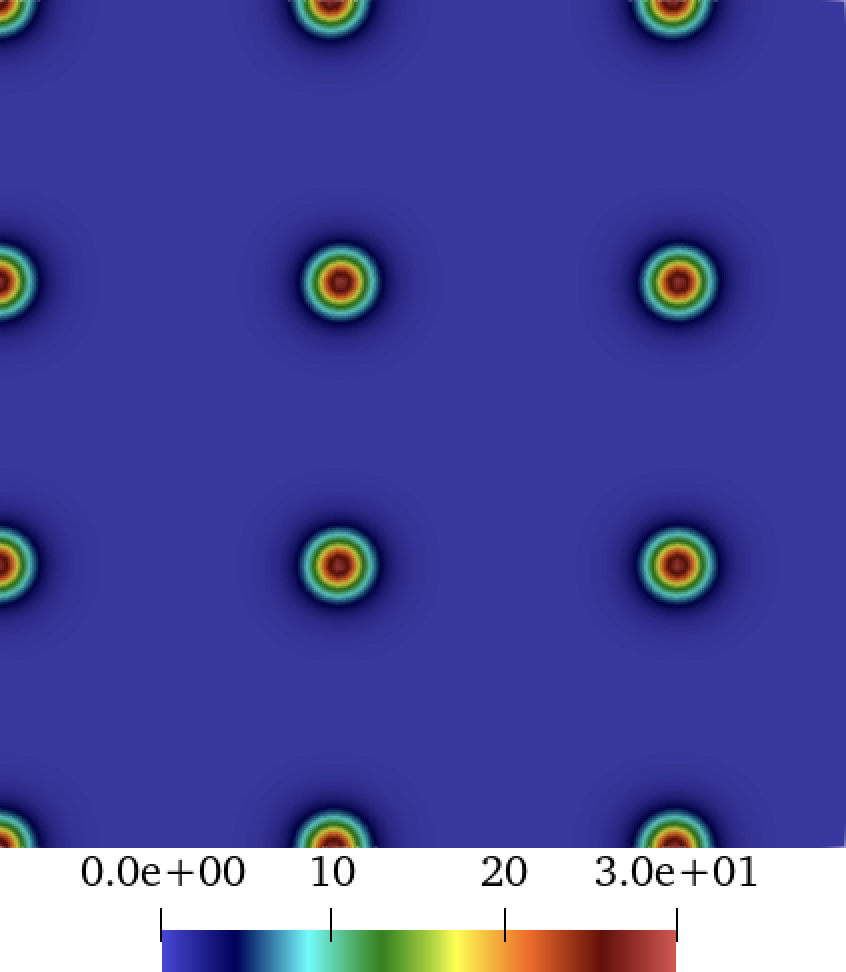}
  \put(23,102){\footnotesize{$t=250$}}
    \end{overpic} 
        \begin{overpic}[percent,width=0.16\textwidth]{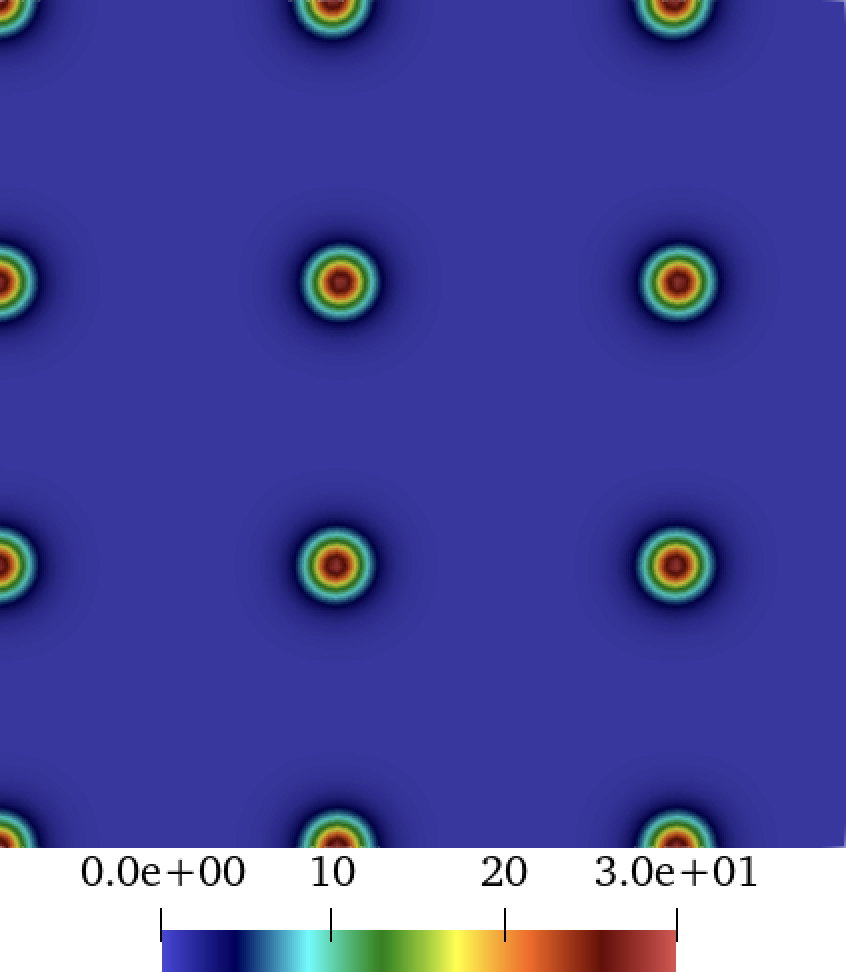}
       \put(23,102){\footnotesize{$t=500$}}
    \end{overpic} 
      \begin{overpic}[percent,width=0.16\textwidth]{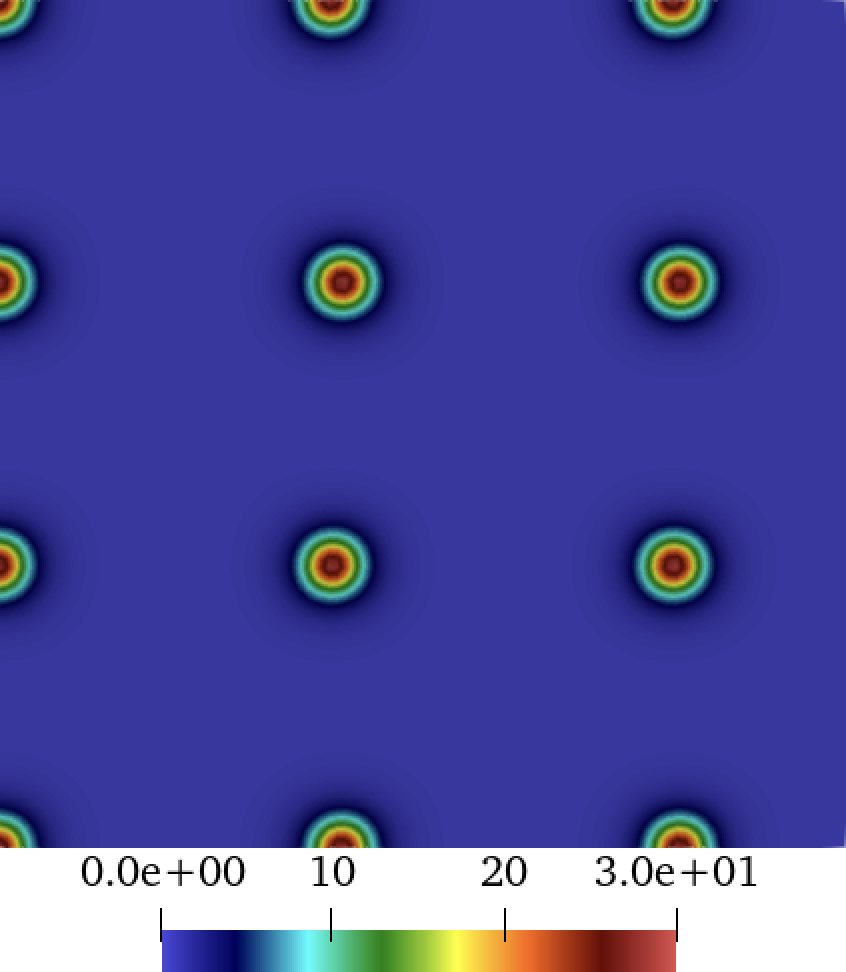}
       \put(23,102){\footnotesize{$t=750$}}
    \end{overpic} 
  \begin{overpic}[percent,width=0.16\textwidth]{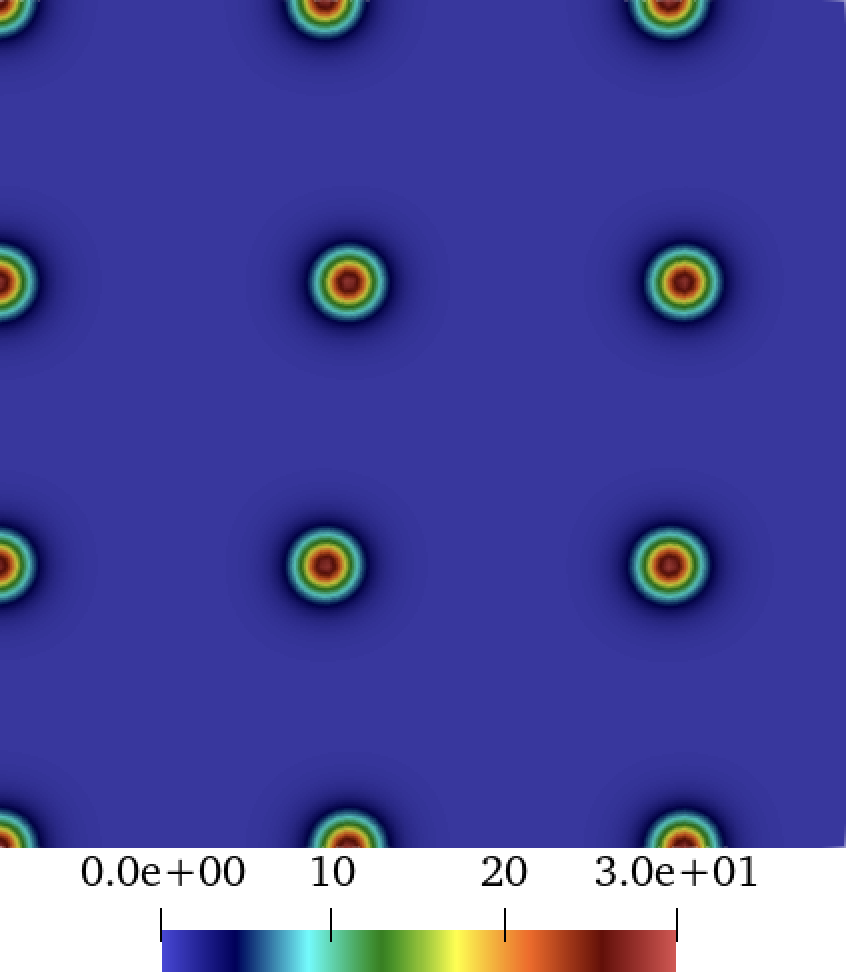}
      \put(23,102){\footnotesize{$t=1000$}}
    \end{overpic} 

    \caption{Expected number of crimes per unit area $S$ for $\tau = 10^{12}$  from $t = 0$ to $t = 1000$.}
    \label{fig:tau1e12_snap}
\end{figure}

The second limit in which the delayed feedback loop breaks down is $\bar{\pi} \to 0$, i.e., no police. When $\bar{\pi}=0$,
$2\nabla H/H$ in~\eqref{eq:continous_3} vanishes
identically since $H$ is spatially uniform at equilibrium (thus,
$\nabla H = 0$), so \eqref{eq:continous_3} reduces to
a pure diffusion equation that decouples from the remaining system of eq.~\eqref{eq:continous_1}, \eqref{eq:continous_2}, and 
\eqref{eq:continous_4}. 
This means that the delayed variable $H$ becomes inert, severing the feedback loop between crime
data and police deployment.
Similarly, 
setting $\bar{\pi} = 0$ in the linearized
system~\eqref{eq:eig}, entry $2\mu_j \bar{\pi}/\bar{H}$
in the third row vanishes,
decoupling the $\pi$-equation from the
$(A,\rho,H)$ block. Then,
the characteristic polynomial 
\eqref{eq:factorization} reduces
to
\begin{equation}\label{eq:char_poly_nopolice}
P(\lambda)\big|_{\bar{\pi}=0}
= (\lambda + \mu_j)\!\left(\lambda + \tfrac{1}{\tau}\right)
  \!\left[\lambda^2 + b_1\lambda + b_0\right],
\end{equation}
where $b_1$ and $b_0$ given by 
\eqref{eq:b1_app} and 
\eqref{eq:b0_app}.
%
%
Polynomial \eqref{eq:char_poly_nopolice} has two roots whose values are
fully explicit: $\lambda = -\mu_j \le 0$ and
$\lambda = -1/\tau < 0$ for every finite $\tau > 0$. Both roots are
real and strictly negative and therefore cannot cross the
imaginary axis as $\tau$ is varied. 
As for the other two roots, 
they do not depend on $\tau$ because 
neither
$b_0$ nor $b_1$ depends on $\tau$.  These 
$\tau$-independent roots
coincide with the eigenvalues of the linearized $(A,\rho)$
subsystem of the model with no police~\cite{short2008statistical},
whose stability properties are governed by $\eta$ and
$\Gamma\theta/\omega^2$ alone.
No root of~\eqref{eq:char_poly_nopolice} can cross the imaginary axis as $\tau$ is varied, and a Hopf
bifurcation is impossible. 
This confirms analytically that
police presence ($\bar{\pi} > 0$) is a necessary
condition for delay-induced oscillations: without enforcement,
the system reverts to the steady hotspot regime
of~\cite{short2008statistical,Hao2026}.

Finally, we note that 
although both the $\tau \to \infty$ and the $\bar{\pi} \to 0$ limits suppress Hopf bifurcations, the
mechanisms are structurally distinct. The infinite-delay
limit zeroes out the constant term of $P(\lambda)$, while
the no-police limit factors $P(\lambda)$ into decoupled
blocks with $\tau$-free roots.

\end{document}